\documentclass[10pt,a4paper]{article}
\usepackage[top=70pt,bottom=90pt,left=64pt,right=64pt]{geometry}
\usepackage{amsmath} %For its many new commands
\usepackage{amsthm}
\usepackage{amssymb} %Some extra symbols
\usepackage{amsfonts}
\usepackage{mathrsfs}
\usepackage{graphicx} %If you want to include postscript graphics
\usepackage{inputenc}
\usepackage{upref} %use \ref
\usepackage{float}
\usepackage[all]{xy}
\usepackage{bm}
\usepackage{calc}
\usepackage{hyperref}
\usepackage{tikz-cd}
\usepackage[all]{xy}
\usepackage{hyperref}
\usepackage[title]{appendix}

\def\makeautorefname#1#2{\expandafter\def\csname#1autorefname\endcsname{#2}}
\def\equationautorefname~#1\null{(#1)\null}

\makeautorefname{theorem}{Theorem}
\makeautorefname{lemma}{Lemma}
\makeautorefname{proposition}{Proposition}
\makeautorefname{remark}{Remark}
\makeautorefname{definition}{Definition}
\makeautorefname{notation}{Notation}
\makeautorefname{conjecture}{Conjecture}
\makeautorefname{construction}{Construction}
\makeautorefname{corollary}{Corollary}
\makeautorefname{question}{Question}

\newtheorem{theorem}{Theorem}[section]

\newtheorem{lemma}{Lemma}[section]
\newtheorem{definition}{Definition}[section]

\newtheorem{proposition}{Proposition}[section]
\newtheorem{notation}{Notation}[section]

\newtheorem{remark}{Remark}[section]

\newtheorem{question}{Question}[section]

\makeatletter
\let\c@lemma=\c@theorem
\let\c@proposition=\c@theorem
\let\c@remark=\c@theorem
\let\c@definition=\c@theorem
\let\c@conjecture=\c@theorem
\let\c@construction=\c@theorem
\let\c@corollary=\c@theorem
\let\c@question=\c@theorem
\let\c@notation=\c@theorem
\numberwithin{equation}{section}
\makeatother

\newcommand{\F}{\mathbb{F}}
\newcommand{\R}{\mathbb{R}}

\newcommand{\Q}{\mathbb{Q}}

\newcommand{\Z}{\mathbb{Z}}

\title{Computing Equivariant Homotopy with a Splitting Method}
\author{Yutao Liu \\ University of Chicago
\\ Email: \href{mailto:yutao492@math.uchicago.edu}{yutao492@math.uchicago.edu}
}

\begin{document}
\maketitle

\paragraph{Abstract:} We develop a new method in the computation of equivariant homotopy, which is based on the splitting of cofiber sequences associated to universal spaces in the category of equivariant spectra. In particular, we use this method to compute the homotopy of $H\underline{\Z}$ for $G=D_{2p}$ and $A_5$.
\medskip

\tableofcontents

\section{Introduction}

Let $G$ be a compact Lie group. In \cite{LMM}, Lewis, May and McClure defined Eilenberg-Maclane $RO(G)$-graded cohomology with Mackey functor coefficients, which 
gave rise to a fundamental but hard question in equivariant stable homotopy: computing $RO(G)$-graded cohomology (or homology) of a point, which is equivalent to the equivariant homotopy of Eilenberg-Maclane spectra. We make this question more precise:

\begin{question}
Let $M$ be an arbitrary $G$-Mackey functor. Compute 
$$\pi_V^G(HM):=[S^V,HM]^G$$
for any virtual $G$-representation $V$. Use $\bigstar$ to denote the $RO(G)$-grading. Then we want to compute $\pi_\bigstar^G(HM)$ as an $RO(G)$-graded abelian group.

Moreover, we have the Mackey functor valued homotopy $\underline{\pi}_\bigstar^G(HM)$ given by
$$\underline{\pi}_V^G(HM)(G/H):=[G/H_+\wedge S^V,HM]^G\cong[S^V,HM]^H.$$
A more general question is to compute $\underline{\pi}_\bigstar^G(HM)$ as an $RO(G)$-graded Mackey functor.

In addition, when $M$ is chosen as a Green functor, which is a monoid in the category of Mackey functors, there will be extra multiplicative structures on the homotopy groups of $HM$. A final question is to compute $\pi_\bigstar^G(HM)$ or $\underline{\pi}_\bigstar^G(HM)$ as an $RO(G)$-graded ring or Green functor.
\end{question}

Common choices of $M$ include constant Mackey functors $\underline{\Z}$, $\underline{\F_p}$, $\underline{\Q}$ and the Burnside ring Mackey functor $A_G$, which is the unit object in the category of Mackey functors. All these choices admit natural Green functor structures.
\bigskip

When $G$ is a cyclic group with prime order, the question is completely solved in multiple ways. Partial computations have been done for larger groups like $C_{p^2}$, $C_2^n$, dihedral group $D_{2p}$, and quaternion group $Q_8$, for which we refer to \cite{Ell}, \cite{Geo}, \cite{HK}, \cite{KL}, \cite{Lu}, \cite{Zen}, and \cite{Zou}. In this paper, we will develop a new method which computes equivariant stable homotopy by decomposition, which works for general equivariant spectra:

\begin{theorem} \label{catmain}
(\textbf{Main Theorem}) Assume that $G$ is a finite group whose order contains multiple prime factors. For any $G$-spectrum $X$ and $G$-virtual representation $V$, $\pi_V^G(X)$ can be expressed as the limit of a diagram with objects $\pi_{V|_H}^H(X)$ and $\pi_{(V|_L)^K}^{L/K}(\Phi^KX)$ for proper $H\subset G$ and $K\lhd L\subset G$ (with proper localizations) such that $|H|,|L/K|$ are prime powers, and maps induced by restriction, conjugacy and localization.

In other words, the whole equivariant stable homotopy theory with finite acting groups is encoded in the cases of $p$-groups.
\end{theorem}

A more detailed statement is given in \autoref{machineinput}.
\medskip

The idea of our method is motivated by the study of rational $G$-spectra for finite $G$, which is first discussed by Greenlees and May:

\begin{theorem} \label{idemsplit}
\cite{GM} There is an orthogonal basis $\{e_H:H\subset G\}$ of the rational Burnside ring, which only contains idempotent elements. For any rational $G$-spectrum $X$, define $e_HX$ as
$$e_HX:=colim(X\xrightarrow{e_H}X\xrightarrow{e_H}X\xrightarrow{e_H}...)$$
Then we have
$$X\simeq\bigvee_H e_HX$$
and
$$[X,Y]^G\cong\prod_H[e_HX,e_HY]^G$$
with one $H$ chosen from each conjugacy class of subgroups.
\end{theorem}

This theorem is further studied by Barnes in \cite{Bar}, which reproves the splitting above by some topological constructions and sets up an algebraic model which completely describes the behavior of rational $G$-spectra.

One important idea in Barnes' proof is to use the universal $G$-space $E\mathscr{F}_H$, which is characterized by the fixed point subspaces:
$$(E\mathscr{F}_H)^K\simeq
\begin{cases}
S^0,\text{ if } K \text{ is conjugate to } H,\\
*, \text{ otherwise.}
\end{cases}$$

\begin{theorem} \label{computesplit}
\cite{Bar} For rational $G$-spectra $X,Y$, we have an isomorphism
$$[X,Y]^G\cong\prod_H [E\mathscr{F}_H\wedge X,E\mathscr{F}_H\wedge Y]^G$$
with one $H$ chosen from each conjugacy class. Moreover, the functor
$$X\mapsto\bigvee_H E\mathscr{F}_H\wedge X$$
is symmetric monoidal.
\end{theorem}

The computation in the rational world is not hard since the category of rational $G$-spectra splits into small and simple pieces (with one piece corresponding to each conjugacy class of subgroups). In fact, in order to obtain this full splitting, it suffices to invert all prime factors of $|G|$ instead of applying rationalization.
\bigskip

Our idea comes from a weaker splitting. If we just invert some prime factors of $|G|$, sometimes the category of $G$-spectra will still split, but into some less simple pieces. If we choose the inverted primes properly, the homotopy of the split pieces may still remain computable. Thus we can compute the localized homotopy by collecting all the pieces.
\smallskip

Each localization loses torsion information at the inverted primes, but that information is retained at other localizations and no torsion information is lost when we glue together different localizations. That is, when $|G|$ contains multiple prime factors, we can apply the idea above multiple times with different prime factors inverted each time. The unlocalized homotopy can be recovered by collecting different localized data.
\medskip

For example, let $G=D_{2p}$. When inverting $2$ or $p$, we get two different splittings in the category of $G$-spectra. Both make the equivariant homotopy computable. For any $G$-spectrum $X$, we can compute $\pi_\bigstar^G(X)[1/2]$ and $\pi_\bigstar^G(X)[1/p]$. The unlocalized $\pi_\bigstar^G(X)$ can be computed by the following pullback diagram:
$$\xymatrix{
\pi_\bigstar^G(X) \ar[r] \ar[d] & \pi_\bigstar^G(X)[1/p] \ar[d] \\
\pi_\bigstar^G(X)[1/2] \ar[r] & \pi_\bigstar^G(X)[1/2,1/p]
}$$
\bigskip

The key element in our splitting method is the universal space $E\mathscr{F}$, where $\mathscr{F}$ is a family of subgroups of $G$. There is a natural cofiber sequence
$$E\mathscr{F}_+\rightarrow S^0\rightarrow\widetilde{E\mathscr{F}},$$
where $\widetilde{E\mathscr{F}}$ is the unreduced suspension of $E\mathscr{F}$. When certain prime factors of $|G|$ are inverted, the cofiber sequence above will split after taking $\Sigma^\infty$. Thus we have
$$X\simeq (E\mathscr{F}_+\wedge X)\vee(\widetilde{E\mathscr{F}}\wedge X)$$
for any $G$-spectrum $X$. 

The inverted prime factors depend on the choice of $\mathscr{F}$. When we invert all prime factors of $|G|$, we can choose all families of subgroups together, which will result in a full splitting of $G$-spectra and recover \autoref{computesplit}. In our computation, we invert all but one prime factor and make special choices of families. Details will be provided in sections 3 and 4.
\bigskip

\paragraph{Structure of the paper:}

In section 2, we will introduce the universal space $E\mathscr{F}$. We are especially interested in the equivariant homology of $E\mathscr{F}$ with coefficients in the Burnside ring Mackey functor $A_G$, which is the starting point of our splitting method. To be precise, we will prove:

\begin{proposition} \label{homefall}
When $|G|$ is inverted, the Mackey functor valued homology
$$\underline{H}_*^G(E\mathscr{F};A_G)$$
concentrates in degree $0$. Moreover, the trivial map $E\mathscr{F}\rightarrow{*}$ induces an inclusion
$$\underline{H}_0^G(E\mathscr{F};A_G)\hookrightarrow\underline{H}_0^G(*;A_G)\cong A_G$$
which makes $\underline{H}_0^G(E\mathscr{F};A_G)$ into a direct summand of $A_G$.
\end{proposition}

In section 3, we will study the cofiber sequence
$$E\mathscr{F}_+\rightarrow S^0\rightarrow\widetilde{E\mathscr{F}}.$$
We will give a necessary and sufficient criterion about which primes we need to invert in order to split the cofiber sequence above after applying $\Sigma^\infty$.

Section 4 is the core of the paper, in which our splitting method is fully explained. We will explain the choice of families and how to decompose the split pieces into the data of smaller subgroups of $G$. The splitting method itself is not a computation, but a machine which uses the computations on smaller groups to obtain the data on larger groups. The input of the machine consists of equivariant homotopy on $p$-subgroups of $G$ with restriction and conjugacy maps. We will also discuss the computability of $H\underline{\Z}$ and $MU_G$ with our splitting method.

In sections 5 and 6, we will apply our splitting method to compute the homotopy of $H\underline{\Z}$ for $G=D_{2p}$ and $A_5$.

In section 7, we will explain how to compute the Mackey functor valued homotopy of $H\underline{\Z}$. The input of our machine becomes the Mackey functor valued $\underline{\pi}_\bigstar^H(H\underline{\Z})$ for proper subgroups $H$ of $G$. The structure maps in $\underline{\pi}_\bigstar^G(H\underline{\Z})$ can be expressed by the structure maps in $\underline{\pi}_\bigstar^H(H\underline{\Z})$. We will explicitly compute the Mackey functor valued $\underline{\pi}_V^{A_5}(H\underline{\Z})$ for some special choices of $V$.
\bigskip

\paragraph{List of computations:}
$\pi_\bigstar^{D_{2p}}(H\underline{\Z})$: \autoref{d2phz}.

$\pi_\bigstar^{A_5}(H\underline{\Z})$: Theorems \ref{3local}, \ref{5local}, \ref{2local}, \ref{positivez}, and \ref{negativez}.

General computations on $D_{2p},A_5$-spectra: Theorems \ref{d2pgeneral}, \ref{a5generalpre}, and \ref{a5general}.
\bigskip

\paragraph{Notations:} In this paper, we use $*$ when the homotopy or homology is graded over $\Z$, and use $\bigstar$ when graded over $RO(G)$. We add an underline to express the Mackey functor valued homotopy or homology: $\underline{H}$, $\underline{\pi}$.

To be more precise, let $X$ be an unbased $G$-space and $M$ be a Mackey functor. Define the Mackey functor valued equivariant homology of $X$ with coefficients in $M$ as
$$\underline{H}_*^G(X;M):=\underline{\pi}_*^G(X_+\wedge HM).$$
So we have
$$\underline{H}_*^G(X;M)(G/H)=[\Sigma^*G/H_+,HM\wedge X_+]^G,$$
$$H_*^G(X;M)=\underline{H}_*^G(X;M)(G/G).$$
Here $HM$ is the equivariant Eilenberg-Maclane spectrum corresponding to $M$.
\smallskip

In this paper, $G$ will always be a finite group.
\bigskip

\paragraph{Acknowledgment:} The author would like to thank his PhD advisor, Peter May, for numerous suggestions on the organization and editing of this paper.

\section{Universal spaces and homological properties}

We introduce the universal space $E\mathscr{F}$ and study its equivariant homology in this section. The most important properties are given in \autoref{homsplit} and \autoref{torsion}.

\subsection{Universal $G$-spaces}

\begin{definition} \label{clspace}
A \textbf{family} $\mathscr{F}$ is a collection of subgroups of $G$ which is closed under conjugation and taking subgroups. The corresponding \textbf{universal space} $E\mathscr{F}$ is an unbased $G$-space such that $(E\mathscr{F})^H$ is contractible for all $H\in\mathscr{F}$ and empty otherwise.
\end{definition}

\begin{remark}
When $\mathscr{F}$ only contains the trivial subgroup, $E\mathscr{F}$ becomes $EG$, which is a contractible space with free $G$-action.
\end{remark} 
\smallskip

One construction of $E\mathscr{F}$ is given in \cite{Dieck} as a join of spaces. We will explain this idea later and use it to prove some important properties. An alternative construction is given by \cite{Elm} as a categorical bar construction.
\smallskip

We can fully characterize $E\mathscr{F}$ by fixed point subspaces:

\begin{lemma} \label{unique}
The universal space $E\mathscr{F}$ is unique up to weak equivalence.
\end{lemma}

\paragraph{Proof:} Assume that $X,Y$ are $G$-spaces with the same conditions on fixed point subspaces as $E\mathscr{F}$.

If we have a $G$-map $X\rightarrow Y$, its restriction on each fixed point subspace is a weak equivalence since $X,Y$ have the same types of fixed point subspaces as either $*$ or $\emptyset$. Thus this map is a weak $G$-equivalence.

If we do not have a map between $X,Y$, consider $X\times Y$, which is another $G$-space with the same conditions on fixed point subspaces. The above argument shows that the projections $X\times Y\rightarrow X$ and $X\times Y\rightarrow Y$ are weak $G$-equivalences. Therefore, $E\mathscr{F}$ is unique up to weak equivalence. $\Box$
\bigskip

We assume $E\mathscr{F}$ to be a $G$-CW complex by applying the $G$-CW approximation.
\begin{lemma} \label{uniqueef}
Any $G$-self map of $E\mathscr{F}$ is $G$-homotopic to the identity.
\end{lemma}

We need the following definition first:

\begin{definition}
Let $\mathcal{O}_G$ be the category of $G$-orbits. Define an $\mathcal{O}_G$-space to be a contravariant functor from $\mathcal{O}_G$ to topological spaces.
\end{definition}

\begin{theorem} \label{elmendorf}
\cite[VI.6]{May} The following pair of functors
$$\Phi: (G-spaces)\rightleftarrows(\mathcal{O}_G-spaces):\Psi$$
defined by $\Phi(X)(G/H):=X^H$ and $\Psi(T):=T(G/\{e\})$ form a Quillen equivalence.
\end{theorem}

\paragraph{Proof of \autoref{uniqueef}:} Construct an $\mathcal{O}_G$-space $T$ by setting 
$$T(G/H)=
\begin{cases}
*, \text{ if }H\in\mathscr{F}\\
\emptyset, \text{ otherwise}
\end{cases}$$

Let $CT$ be a cofibrant approximation of $T$. Then
$$\Psi(CT)^H=\Phi\Psi(CT)(G/H)\cong CT(G/H)\cong
\begin{cases}
*, \text{ if }H\in\mathscr{F}\\
\emptyset, \text{ otherwise}
\end{cases}$$

Recall that
$$\Phi E\mathscr{F}(G/H)=(E\mathscr{F})^H\simeq
\begin{cases}
*, \text{ if }H\in\mathscr{F}\\
\emptyset, \text{ otherwise}
\end{cases}$$

The unique map $\Phi E\mathscr{F}\rightarrow T$ induces $\Phi E\mathscr{F}\rightarrow CT$, whose left adjoint $E\mathscr{F}\rightarrow\Psi(CT)$ becomes a $G$-equivalence.

So we have the following equivalences of homotopy classes on unbased maps:
$$[E\mathscr{F},E\mathscr{F}]^G\cong[E\mathscr{F},\Psi(CT)]^G\cong[\Phi E\mathscr{F},CT]^{\mathcal{O}_G}\cong[\Phi E\mathscr{F},T]^{\mathcal{O}_G},$$
which contains a single element.

Thus the identity map is the only self-map of $E\mathscr{F}$ up to $G
$-homotopy. $\Box$

\bigskip

Notice that $E\mathscr{F}\times E\mathscr{F}$ also appears as a universal space for $\mathscr{F}$. Thus the above two lemmas imply:

\begin{lemma} \label{monoid}
The universal space $E\mathscr{F}$ is a $G$-topological semigroup up to homotopy.
\end{lemma}

Usually $E\mathscr{F}$ is not a Hopf space since there is no unit. But under proper localizations, $\Sigma^\infty E\mathscr{F}_+$ will appear as a homotopy ring spectrum. We will explain this idea in section 3.
\medskip

Notice that the based spaces $E\mathscr{F}_+$ and $\widetilde{E\mathscr{F}}$ are also characterized by their fixed point subspaces:
$$(E\mathscr{F}_+)^H\simeq
\begin{cases}
S^0,\text{ if }H\in\mathscr{F},\\
*, \text{ otherwise.}
\end{cases}$$
$$(\widetilde{E\mathscr{F}})^H\simeq
\begin{cases}
*,\text{ if }H\in\mathscr{F},\\
S^0, \text{ otherwise.}
\end{cases}$$

Thus we can apply the same argument as above, but in a based version:

\begin{proposition} \label{basedmonoid}
Any self $G$-map of $E\mathscr{F}_+$ or $\widetilde{E\mathscr{F}}$ which is a $G$-equivalence is $G$-homotopic to the identity.

The based spaces $E\mathscr{F}_+$ and $\widetilde{E\mathscr{F}}$ are based $G$-topological semigroups up to $G$-homotopy. In addition, the natural map $S^0\rightarrow\widetilde{E\mathscr{F}}$ commutes with the following products up to $G$-homotopy:
$$S^0\wedge S^0\xrightarrow{\simeq}S^0,$$
$$\widetilde{E\mathscr{F}}\wedge\widetilde{E\mathscr{F}}\xrightarrow{\simeq}\widetilde{E\mathscr{F}}.$$
\end{proposition}

\subsection{$\Z$-graded homology with coefficients in $A_G$}

We recall the definition of the Burnside ring and the corresponding Mackey functor:

\begin{definition} \label{ag}

For any finite group $G$, the collection of isomorphism classes of finite $G$-sets forms a commutative monoid, with addition induced by disjoint unions. The \textbf{Burnside ring} $A(G)$ is defined to be the group completion of this monoid. We can express $A(G)$ as a free $\Z$-module, whose basis corresponds to isomorphism classes of $G$-orbits. We use $\{G/H\}$ to denote the basis element of the orbit $G/H$.
\medskip

The \textbf{Burnside ring Mackey functor} $A_G$ is defined by $A_G(G/H):=A(H)$. The Mackey functor structure is given by the following maps:
\smallskip

For any $L\subset H$, the transfer map 
$$T_L^H: A(L)\rightarrow A(H)$$
sends $\{L/K\}$ to $\{H/K\}$ for any $K\subset L$. 

The restriction map
$$R_L^H: A(H)\rightarrow A(L)$$
sends $\{H/N\}$ to itself, but viewed as an $L$-space.

When $L=g^{-1}Hg$, we also have an isomorphism
$$C_L^H:A(L)\rightarrow A(H)$$
sending $\{L/K\}$ to $\{H/gKg^{-1}\}$ for any $K\subset L$.
\medskip

The Cartesian product of finite $G$-sets makes $A(G)$ into a ring. This multiplicative structure also makes $A_G$ into a Green functor, which is a monoid under the box product.
\end{definition}

\begin{proposition} \label{burnsidehom}
For a $G$-CW complex $X$, the integer graded homology of $X$ with coefficients in $A_G$ can be computed as
$$H_*^G(X;A_G)\cong\bigoplus_H H_*(X^H/W_GH;\Z)$$
with one $H$ chosen from each conjugacy class of subgroups of $G$. Here $W_GH$ is the Weyl group of $H$ in $G$.

As a Mackey functor,
$$\underline{H}_*^G(X;A_G)(G/L)\cong H_*^L(X;A_L)\cong\bigoplus_K H_*(X^K/W_LK;\Z)$$
with one $K$ chosen from each conjugacy class of subgroups of $L$.
\end{proposition}

\paragraph{Proof:} It suffices to prove the first equation. The second one can be proved by the same argument.
\medskip

When we view $A_G$ as a covariant coefficient system, the restriction maps are removed. For any $G$-map $G/L_1\rightarrow G/L_2$ induced by the multiplication of $g\in G$ (where we require $g^{-1}L_1g\subset L_2$), we have an induced map
$$A(L_1)\rightarrow A(L_2)$$
sending each $\{L_1/H\}$ to $\{L_2/g^{-1}Hg\}$. Therefore, as a coefficient system, $A_G$ can be decomposed as
$$A_G=\bigoplus_{H}A_G^H$$
with one $H\subset G$ chosen from each conjugacy class. Here $A_G^H$ is the sub-coefficient system of $A_G$ such that each $A_G^H(G/L)\subset A(L)$ is generated by all $\{L/g^{-1}Hg\}$ with $g^{-1}Hg\subset L$. 
\smallskip

Now we have
$$H_*^G(X;A_G)\cong\bigoplus_{H}H_*^G(X;A_G^H).$$
It suffices to prove
$$H_*^G(X;A_G^H)\cong H_*(X^H/W_GH;\Z)$$
for any $H\subset G$.
\bigskip

As described in \cite{Wil}, the equivariant homology can be computed in a cellular way:

Let $\underline{C}_*(X)$ be the chain complex of contravariant coefficient systems such that
$$\underline{C}_*(X)(G/H):=C_*(X^H),$$
where $C_*(X^H)$ is the cellular chain complex of $X^H$ and the orbit maps are sent to conjugacies and inclusions of subcomplexes. For any covariant coefficient system $M$, $H_*^G(X;M)$ is the homology of the chain complex:
$$C_*^G(X;M):=\underline{C}_*(X)\otimes_{\mathcal{O}_G}M,$$
where the tensor product is taken over the category $\mathcal{O}_G$ of $G$-orbits. To be more precise, we have
$$C_*^G(X;M):=\left(\bigoplus_{L\subset G}C_*(X^L)\otimes M(G/L)\right)/\sim.$$

For any $G$-map $f:G/L_1\rightarrow G/L_2$, the equivalence relation identifies $f^*a\otimes b$ and $a\otimes f_*b$ for all $a\in C_*(X^{L_2})$ and $b\in M(G/L_1)$.
\medskip

Now we choose $M=A_G^H$. Since each $M(G/L)$ is freely generated by $\{L/g^{-1}Hg\}$ and the structure maps send these basis elements to each other, $C_*^G(X;M)$ is the free $\Z$-module generated by the equivalence classes of $e\otimes \{L/g^{-1}Hg\}$, for all cells $e$ in $X^L$ and $g^{-1}Hg\subset L$.
\medskip

We can simplify the expression of $C_*^G(X;A_G)$ by the following four facts:

\textbf{(1)} The free $\Z$-module $A_G^H(G/H)$ is generated by a single element $\{H/H\}$.

\textbf{(2)} Each $e\otimes\{L/g^{-1}Hg\}$ is identified with $ge\otimes\{H/H\}$ by the equivalence relation.

\textbf{(3)} If there exists another $g^\prime\in G$ such that
$$ge\otimes\{H/H\}\sim e\otimes\{L/g^{-1}Hg\}\sim g^\prime e\otimes\{H/H\},$$
then $\{L/g^{-1}Hg\}=\{L/(g^{\prime})^{-1}Hg^\prime\}$. Thus
$$g^{-1}Hg=l^{-1}((g^\prime)^{-1}Hg^\prime)l$$
for some $l\in L$. So $g^\prime lg^{-1}\in W_GH$. Since $e$ is $L$-fixed, we have $(g^\prime lg^{-1})(ge)=g^\prime e$. Thus $ge$ and $g^\prime e$ are in the same $W_GH$-orbit.

\textbf{(4)} The converse of \textbf{(3)} is also true: $e_1\otimes\{H/H\}$ and $e_2\otimes\{H/H\}$ are identified if $e_1,e_2$ are in the same $W_GH$-orbit.
\medskip

Now we have a 1-1 correspondence between the basis of $C_*^G(X;A_G^H)$ and the cells in $X^H/W_GH$. Therefore, we get
$$C_*^G(X;A_G^H)\cong C_*(X^H/W_GH).$$
Taking the homology on both sides gives us the required equation. $\Box$
\bigskip

\begin{remark}
\autoref{burnsidehom} follows from the fact that the underlying coefficient system of $A_G$ splits into $A_G^H$. However, such splitting cannot be lifted to the Mackey functor level. Thus the decomposition does not work for cohomology with coefficients in $A_G$.

In fact, since the orbit spectra $(-)/G$ of genuine $G$-spectra can only be defined after passing from a complete universe to the trivial universe,
$$X\mapsto H_*(X^H/W_GH)$$
is a homology theory only for $G$-spaces, and hence cannot be represented by any $G$-spectrum.
\end{remark}

\subsection{Homology of $E\mathscr{F}$ in degree $0$}

Since $E\mathscr{F}$ only has empty and contractible fixed point subspaces, \autoref{burnsidehom} helps us to compute the 0th degree equivariant homology of $E\mathscr{F}$ explicitly:

\begin{proposition} \label{degree0}
The trivial map $E\mathscr{F}\rightarrow\{*\}$ implies an inclusion
$$\underline{H}_0^G(E\mathscr{F};A_G)\hookrightarrow\underline{H}_0^G(*;A_G)\cong A_G.$$

The image of $\underline{H}_0^G(E\mathscr{F},A_G)(G/L)$ in $A(L)$ is generated by all $\{L/K\}\in A(L)$ with $K\in\mathscr{F}$. 
\end{proposition}

We already have a small splitting here:

\begin{proposition} \label{homsplit}
When $|G|$ is inverted, $\underline{H}_0^G(E\mathscr{F};A_G)$ is a direct summand of $A_G$ as Mackey functors.
\end{proposition}

\begin{notation} \label{mfnf}
We use $M_\mathscr{F}$ to denote the Mackey functor $\underline{H}_0^G(E\mathscr{F};A_G)$ and use the same notation for its localizations. 

When $M_\mathscr{F}$ appears to be a direct summand of $A_G$, denote its complement as $N_\mathscr{F}$. A more concrete description of $N_\mathscr{F}$ is given in the proof of \autoref{homsplit} below.
\end{notation}

\paragraph{Proof of \autoref{homsplit}:} Assume that $|G|$ is inverted everywhere. For any $H\subset G$, define a linear map
$$\chi_H:A(G)\rightarrow\Z$$
which sends each $G$-set $S$ to $|S^H|$. 

Let $s_{(K,H)}$ denote $\chi_K(\{G/H\})=|G/H|^K$, which can be computed as the product between $|W_GH|$ and the number of subgroups of $G$ containing $K$ and in the conjugacy class of $H$. Thus $s_{(K,H)}\neq 0$ if and only if $K$ is sub-conjugate to $H$.

When $|G|$ is inverted, all $|W_GH|$ and non-zero $s_{(K,H)}$ become invertible. We choose elements $e_H\in A(G)$ inductively by defining
$$e_H=|W_GH|^{-1}\left(\{G/H\}-\sum_K s_{(K,H)}^{-1}e_K\right)$$
with one $K$ chosen from each conjugacy class that contains a proper subgroup of $H$.

By induction, we have
$$\chi_K(e_H)=
\begin{cases}
1, \text{ if }K\text{ is conjugate to }H\\
0, \text{ otherwise}
\end{cases}$$

Since the dimension of $A(G)$ agrees with the number of conjugacy classes, the collection of $e_H$, with one $H$ chosen from each conjugacy class, forms a basis of $A(G)$.
\medskip

In general, for any $H\subset L\subset G$, define a linear map
$$\chi_H^L:A(L)\rightarrow\Z$$
which sends each $L$-set to the size of its $H$-fixed subset. We can get a similar basis $\{e_H^L\}$, with one $H$ chosen from each conjugacy class of subgroups of $L$, such that
$$\chi_K^L(e_H^L)=
\begin{cases}
1, \text{ if }K\text{ is conjugate to }H\text{ in }L\\
0, \text{ otherwise}
\end{cases}$$
\medskip

Now we discuss how these maps interact with transfer and restriction maps:
\smallskip

For any $K\subset L_1\subset L_2\subset G$, the restriction map $R_{L_1}^{L_2}$ keeps each $L_2$-set but views it as an $L_1$-set, hence does not change its $K$-fixed subset. So we have
$$\chi_K^{L_2}=\chi_K^{L_1}\circ R_{L_2}^{L_1}.$$

On the other hand, for any $K\subset L_1$, the elements in $\{L_1/K\}$ have isotropy groups conjugate to $K$ inside $L_1$. The elements in $T_{L_1}^{L_2}(\{L_1/K\})=\{L_2/K\}$ have isotropy groups conjugate to $K$ inside $L_2$. There may be more isotropy groups. But these additional isotropy groups are still chosen from the conjugacy class of $K$ inside $G$. According to the construction of the basis element $e_H^L$, we have
$$\chi_H^{L_2}(T_{L_1}^{L_2}e_K^{L_1})\neq 0\text{ only if }H\text{ is conjugate to a subgroup of }K\text{ in }G$$

Let
$$N(L):=\bigcap_{H\in\mathscr{F},H\subset L}\ker\chi_H^L$$
In other words, $N(L)$ is generated by $e_H^L$ for all $H\subset L$ and $H\notin\mathscr{F}$. The above discussion tells us that $N(L)$ is closed under transfer and restriction maps. Thus we get a sub-Mackey functor whose value at $G/L$ agrees with $N(L)$. Denote that as $N_\mathscr{F}$.
\smallskip

Notice that for each $L\subset G$, $e_H^L\in N(L)=N_\mathscr{F}(G/L)$ if $H\notin\mathscr{F}$, $e_H^L\in M_\mathscr{F}(G/L)$ if $H\in\mathscr{F}$, according to \autoref{degree0}. Moreover, consider any nontrivial element
$$a_1\{L/K_1\}+a_2\{L/K_2\}+...+a_n\{L/K_n\}\in M_\mathscr{F}(G/L)$$
with $a_1,a_2,...,a_n\neq 0$, $K_1,K_2,...,K_n$ in different conjugacy classes in $\mathscr{F}$, and $|K_1|\leq |K_2|\leq...\leq|K_n|$. The map $\chi_{K_n}^L$ sends $\{L/K_1\},...,\{L/K_{n-1}\}$ to zero but $\{L/K_n\}$ to a positive value. Thus this element is not in $\ker\chi_{K_n}^L$, and hence $M_\mathscr{F}(G/L)\cap N(L)=\emptyset$.
\medskip

In conclusion, we have $A_G=M_\mathscr{F}\oplus N_\mathscr{F}$ and $M_\mathscr{F}$ appears as a direct summand.. $\Box$

\begin{remark} \label{injectivevar}
When the multiplicative structure is added into consideration, the maps $\chi_H^L$ appear as the components of the ring isomorphism from $A(L)$ to several copies of $\Z[|G|^{-1}]$, with one $H$ chosen from each conjugacy class inside $L$. The generators $e_H^L$ become idempotent elements. So we can view $M_\mathscr{F}$ as a direct summand of $A_G$ as a Green functor.
\end{remark}

\subsection{Homology of $E\mathscr{F}$ in positive degrees}

In positive degrees, we have:

\begin{theorem} \label{torsion}

For any family $\mathscr{F}$, the $\Z$-graded, Mackey functor valued homology
$$\underline{H}_*^G(E\mathscr{F};A_G)$$
contains only torsion when the degree is positive. Moreover, the torsion only has prime factors which divide $|G|$.

\end{theorem}

This is the most important property we want about $E\mathscr{F}$. We will use the rest of this section to prove it.
\bigskip

The main idea is an induction on the size of $\mathscr{F}$. For the base case, we have

\begin{lemma} \label{orduni}

Let $BG=EG/G$ be the classifying space of principal $G$-bundles. Then $H_*(BG;\Z)$ contains only torsion when the degree is positive. Moreover, the torsion only has prime factors which divide $|G|$.

\end{lemma}

This is a standard result about $EG$. We give one possible proof below.

\paragraph{Proof:} Give $EG$ the standard $G$-CW structure which only contains $G$-free cells. Consider the map between cellular chain complexes induced by the projection $EG\rightarrow EG/G=BG$:
$$C_*(EG)\rightarrow C_*(BG).$$
Define another map between chain complexes
$$C_*(BG)\rightarrow C_*(EG)$$
which sends each cell $e$ in $BG$ to the sum of all cells in $EG$ which are sent to $e$ under $EG\rightarrow BG$. Then the composition
$$C_*(BG)\rightarrow C_*(EG)\rightarrow C_*(BG)$$
is the multiplication by $|G|$ since all cells in $EG$ are $G$-free. The same $|G|$-multiplication passes to homology:
$$H_*(BG)\rightarrow H_*(EG)\rightarrow H_*(BG).$$

Since $EG$ is contractible, $H_*(EG)=\Z$, which is concentrated in degree $0$. Thus when the degree is positive, $H_*(BG)$ contains only torsion which divides $|G|$. $\Box$
\bigskip

In the case when $\mathscr{F}$ only contains the trivial subgroup, $E\mathscr{F}$ becomes $EG$ and only has orbit spaces homotopic to $BH$ for $H\subset G$. \autoref{orduni} and \autoref{burnsidehom} imply that it satisfies \autoref{torsion}.
\bigskip

For general $\mathscr{F}$, we can construct $E\mathscr{F}$ as a join. More details about joins and their properties are given in Appendix A.

\begin{proposition} \label{efasjoin}
The join of $G\times_{N_GH}EW_GH$, with one $H$ chosen from each conjugacy class in $\mathscr{F}$, is a valid construction for $E\mathscr{F}$. Here $N_GH$ is the normalizer of $H$ in $G$. We view $EW_GH$ as a left $N_GH$-space and $G$ as a right $N_GH$-set.
\end{proposition}

\paragraph{Proof:} Notice that $E\mathscr{F}$ is characterized by its fixed point subspaces:
$$(E\mathscr{F})^K\simeq
\begin{cases}
*,\text{ if }K\in\mathscr{F},\\
0,\text{ if }K\notin\mathscr{F}
\end{cases}$$
According to \autoref{join}, it suffices to show:

\textbf{(a)} If $K\in\mathscr{F}$, $(G\times_{N_GH}EW_GH)^K$ is contractible for some $H\in\mathscr{F}$.

\textbf{(b)} If $K\notin\mathscr{F}$, $(G\times_{N_GH}EW_GH)^K$ is empty for all $H\in\mathscr{F}$.
\medskip

Consider any point in $G\times_{N_GH}EW_GH$ expressed by $a\times x$ for $a\in G$ and $x\in EW_GH$. For any $g\in G$, if the action of $g$ fixes the point $a\times x$:
$$g(a\times x)=ga\times x=a\times x,$$
then there is $h\in N_GH$, such that $hx=x$, $gah^{-1}=a$. Since all points in $EW_GH$ have isotropy group $H$, we have $h\in H$. Then $g=aha^{-1}\in aHa^{-1}$.

Therefore, the isotropy group of $a\times x$ is $aHa^{-1}$. Thus the $K$-fixed point subspace of $G\times_{N_GH}EW_GH$ is non-empty if and only if $K$ is sub-conjugate to $H$. Moreover, when $K=gHg^{-1}$, we have
$$(G\times_{N_GH}EW_GH)^K=\{g\}\times EW_GH,$$
which is contractible. Thus both conditions (a), (b) are satisfied. $\Box$
\bigskip

\begin{remark} \label{isotropytypes}
Since all points in $G\times_{N_GH}EW_GH$ have isotropy groups conjugate to $H$, according to \autoref{joindef}, the collection of isotropy groups of all points in $E\mathscr{F}$ agrees with $\mathscr{F}$.
\end{remark}
\medskip

\paragraph{Proof of \autoref{torsion}:} Assume that \autoref{torsion} is true for all families smaller than $\mathscr{F}$. Let $\mathscr{F}^\prime$ be a smaller family which is obtained by removing the conjugacy class of one largest subgroup $H\in\mathscr{F}$. According to \autoref{efasjoin}, we have
$$E\mathscr{F}\simeq E\mathscr{F}^\prime\ast(G\times_{N_GH}EW_GH).$$
\medskip

For any $K\lhd L\subset G$, the discussion in the proof of \autoref{efasjoin} shows that the space $(G\times_{N_GH}EW_GH)^K$ is either empty or several copies of $EW_GH$. Since all points in each copy of $EW_GH$ have the same isotropy group, $(G\times_{N_GH}W_GH)^K/L$ is either empty or the disjoint union of $BN$ for some subgroups $N\subset N_GH\subset G$. \autoref{orduni} shows that its positive degree homology only contains torsion dividing $|G|$. Thus \autoref{burnsidehom} shows that $G\times_{N_GH}EW_GH$ satisfies \autoref{torsion}.
\medskip

Since $E\mathscr{F}^\prime$ satisfies \autoref{torsion}, according to \autoref{joinhom}, now it suffices to prove the theorem for $E\mathscr{F}^\prime\times(G\times_{N_GH}EW_GH)$. Using $X(H)$ to denote $G\times_{N_GH}EW_GH$, the homology of $E\mathscr{F}^\prime\times X(H)$ can be computed by the equivariant Künneth spectral sequence \cite{LM}:
$$E_{p,q}^2=\underline{Tor}^{A_G}_{p,q}(\underline{H}_*^G(E\mathscr{F}^\prime;A_G),\underline{H}_*^G(X(H);A_G))\Rightarrow\underline{H}_*^G(E\mathscr{F}^\prime\times X(H);A_G).$$

When we invert $|G|$, the homology of both $E\mathscr{F}^\prime$ and $X(H)$ is concentrated in degree 0. Moreover, $\underline{H}_0^G(E\mathscr{F}^\prime;A_G)=M_{\mathscr{F}^\prime}$ appears as a direct summand of $A_G$ by \autoref{homsplit}. Thus the $E_2$-page collapses into a single box product:
$$E^2=E^2_{0,0}=M_\mathscr{F}\text{ }\Box\text{ } \underline{H}_*^G(X(H);A_G).$$

Therefore, $E\mathscr{F}^\prime\times X(H)$ has trivial homology in positive degrees when $|G|$ is inverted. $\Box$
\bigskip

For special choices of $\mathscr{F}$, the torsion in $\underline{H}^G_*(E\mathscr{F};A_G)$ may not see all prime factors of $|G|$. We emphasize the following criterion which is implied by \autoref{burnsidehom}:

\begin{proposition} \label{torsiontype}
The torsion types in $\underline{H}_*^G(E\mathscr{F},A_G)$ come from $H_*(E\mathscr{F}^K/L)$ for all $K\lhd L\subset G$.
\end{proposition}

\section{Splitting cofiber sequences}

Let $\mathscr{F}$ be a family of subgroups of $G$. Consider the cofiber sequence
$$E\mathscr{F}_+\rightarrow S^0\rightarrow\widetilde{E\mathscr{F}}$$
where the map $E\mathscr{F}_+\rightarrow S^0$ sends $E\mathscr{F}$ into the non-basepoint of $S^0$.

\subsection{A categorical splitting}

\begin{theorem} \label{ringsplit}
When certain prime factors of $|G|$ are inverted, the cofiber sequence above splits the category of $G$-spectra, in the sense that
$$X\simeq (E\mathscr{F}_+\wedge X)\vee(\widetilde{E\mathscr{F}}\wedge X),$$
$$[E\mathscr{F}_+\wedge X,\widetilde{E\mathscr{F}}\wedge Y]^G=[\widetilde{E\mathscr{F}}\wedge Y,E\mathscr{F}_+\wedge X]^G=0$$
for any $G$-spectra $X,Y$. Moreover, the functors
$$X\mapsto E\mathscr{F}_+\wedge X, \text{ }X\mapsto\widetilde{E\mathscr{F}}\wedge X$$
are symmetric monoidal up to homotopy.
\end{theorem}

The proof will give information about which primes must be inverted; see \autoref{invertprime}.
\medskip

Consider the suspension of the cofiber sequence:
$$\Sigma^\infty E\mathscr{F}_+\rightarrow S\rightarrow\Sigma^\infty\widetilde{E\mathscr{F}}.$$

\begin{theorem} \label{cofibsplit}
When $|G|$ is inverted, there exists a left inverse of $\Sigma^\infty E\mathscr{F}_+\rightarrow S$.
\end{theorem}

\paragraph{Proof:} According to \autoref{degree0}, the map $\Sigma^\infty E\mathscr{F}_+\rightarrow S$ induces an inclusion in homology:
$$\underline{HA_G}_0(\Sigma^\infty E\mathscr{F}_+)=:M_\mathscr{F}\hookrightarrow\underline{HA_G}_0S=A_G.$$

When $|G|$ is inverted, \autoref{torsion} shows that the $HA_G$-homology is concentrated in degree 0. Thus it suffices to find a map $S\rightarrow\Sigma^\infty E\mathscr{F}_+$ which induces a projection from $A_G$ to $M_\mathscr{F}$ on their $HA_G$-homology. After composition with $\Sigma^\infty E\mathscr{F}_+\rightarrow S$, we get a self-map of $\Sigma^\infty E\mathscr{F}_+$ that induces an isomorphism on $HA_G$-homology. Since $\Sigma^\infty E\mathscr{F}_+$ is a connective spectrum, this self map must be a weak equivalence.
\smallskip

For each $H\in\mathscr{F}$, the transfer and restriction maps in Mackey functors between images of $G/G$ and $G/H$ are induced by stable maps

$$T: S=\Sigma^\infty S^0=\Sigma^\infty  G/G_+\rightarrow\Sigma^\infty G/H_+$$
$$R: \Sigma^\infty G/H_+\rightarrow \Sigma^\infty G/G_+=S.$$

$R$ is induced by the space-level map $G/H_+\rightarrow S^0$, which sends $G/H$ into the non-basepoint of $S^0$. According to \autoref{isotropytypes}, $E\mathscr{F}$ contains points with any $H\in\mathscr{F}$ as isotropy groups. Thus we can find a map from $G/H$ to $E\mathscr{F}$ and factorize $R$ as
$$\Sigma^\infty G/H_+\rightarrow\Sigma^\infty E\mathscr{F}_+\rightarrow S.$$

Let $\iota_H$ be the composition of the first map above and $T$:
$$\iota_H:S\xrightarrow{T}\Sigma^\infty G/H_+\rightarrow\Sigma^\infty E\mathscr{F}_+.$$

The composition
$$S\xrightarrow{\iota_H}\Sigma^\infty E\mathscr{F}_+\rightarrow S$$
agrees with $R\circ T$, which is the equivariant Euler characteristic of $G/H$. We can fully describe the self-map on
$$\underline{HA_G}_0S\cong A_G\cong\underline{\pi}_0S$$
induced by the Euler characteristic as follows:

For each $L\subset G$, the induced self-map on $A_G(G/L)=A(L)$ is the multiplication by the underlying $L$-set of $\{G/H\}$ (which is defined in \autoref{ag}). We refer to \cite{LMS}, V.1 and V.2] for more details.
\medskip

Now we choose one $H$ from each conjugacy class in $\mathscr{F}$. Assign a random number $c_H\in\Z[|G|^{-1}]$ for each such $H$. Consider the map
$$\sum_H c_H\iota_H:S\rightarrow\Sigma^\infty E\mathscr{F}_+.$$
The main strategy is to choose $c_H$ properly such that the map above induces an isomorphism on equivariant homology with coefficient in $A_G$ after composing with $\Sigma^\infty E\mathscr{F}_+\rightarrow S$, hence appears to be the required left inverse.
\smallskip

The composition
$$S\xrightarrow{\sum_Hc_H\iota_H}\Sigma^\infty E\mathscr{F}_+\rightarrow S$$
induces a self-map on $A_G(G/L)=A(L)$ as the multiplication by $\sum_H c_H\{G/H\}$.

Recall that $\underline{HA_G}_0(\Sigma^\infty E\mathscr{F}_+)(G/L)=M_\mathscr{F}(G/L)\subset A(L)$ is generated by all $\{L/J\}$ with $J\in\mathscr{F}$. The statement that $\sum_Hc_H\iota_H$ induces the projection $A_G\rightarrow M_\mathscr{F}$ is equivalent to
$$\{L/J\}\cdot\sum_Hc_H\{G/H\}=\{L/J\}$$
for any $J\in\mathscr{F}$, $J\subset L$.

Recall the ring map $\chi_K^L:A(L)\rightarrow\Z$ for any $K\subset L$ which sends each $L$-set to the number of $K$-fixed points. According to \autoref{injectivevar}, $\prod_{K\subset L}\chi_K^L$ is an injective ring homomorphism from $A(L)$ to copies of $\Z$. Thus it suffices to check the equation above on the images under each $\chi_K^L$. The equation above under $\chi_K^L$ becomes
$$|(L/J)^K|\cdot\sum_Hc_H|(G/H)^K|=|(L/J)^K|.$$

Notice that $|(L/J)^K|=0$ for any $K\notin\mathscr{F}$ (since $K$ is not sub-conjugate to $J\in\mathscr{F}$). Thus it suffices to choose $c_H$ such that
$$\sum_Hc_H|(G/H)^K|=\sum_H c_Hs_{(K,H)}=1, \text{ } \forall K\in\mathscr{F}$$
with one $H$ chosen from each conjugacy class in $\mathscr{F}$. (Recall that $s_{(K,H)}:=|(G/H)^K|$ is defined in the proof of \autoref{homsplit}.)

Notice that $s_{(K,H)}$ is non-zero if and only if $K$ is sub-congujate to $H$. Denote this relation by $[K]\leq[H]$. We have

$$\sum_{[K]\leq[H]\subset\mathscr{F}} c_Hs_{(K,H)}=1.$$
Thus
$$c_K=s_{(K,K)}^{-1}\left(1-\sum_{[K]\lneqq[H]\subset\mathscr{F}}c_Hs_{(K,H)}\right).$$

Since $s_{(K,K)}=|W_GK|$ is invertible when $|G|$ is inverted, $c_K$ can be chosen inductively from larger subgroups to smaller ones. Moreover, it's clear that the denominator of any $c_K$ only contains prime factors dividing $|G|$.

Now the map $\sum_Hc_H\iota_H$ gives us the required projection on the 0th homology, hence becomes the left inverse of $\Sigma^\infty E\mathscr{F}_+\rightarrow S$. $\Box$
\bigskip

\paragraph{Proof of \autoref{ringsplit}:} Assume that $|G|$ is inverted. \autoref{cofibsplit} implies that the cofiber sequence
$$\Sigma^\infty E\mathscr{F}_+\rightarrow S\rightarrow\Sigma^\infty\widetilde{E\mathscr{F}}$$
splits. Thus we have
$$X\simeq(E\mathscr{F}_+\wedge X)\vee(\widetilde{E\mathscr{F}}\wedge X).$$

In fact, both $\Sigma^\infty E\mathscr{F}_+$ and $\Sigma^\infty\widetilde{E\mathscr{F}}$ become ring spectra (up to homotopy). The unit map on $\Sigma^\infty E\mathscr{F}_+$ is given by the left inverse in \autoref{cofibsplit}, and the product is induced by the space level map $E\mathscr{F}_+\wedge E\mathscr{F}_+\rightarrow E\mathscr{F}_+$. For $\Sigma^\infty\widetilde{E\mathscr{F}}$, the unit and product maps are both induced by the space level maps $S^0\rightarrow\widetilde{E\mathscr{F}}$ and $\widetilde{E\mathscr{F}}\wedge\widetilde{E\mathscr{F}}\rightarrow\widetilde{E\mathscr{F}}$. \autoref{basedmonoid} and the fact that 
$$E\mathscr{F}_+\wedge E\mathscr{F}_+\simeq E\mathscr{F}_+\text{ , }\widetilde{E\mathscr{F}}\wedge\widetilde{E\mathscr{F}}\simeq\widetilde{E\mathscr{F}}$$
guarantee all commutative diagrams in the definition of a ring spectrum.
\smallskip

Since both
$$S\rightarrow\Sigma^\infty E\mathscr{F}_+\text{ , }S\rightarrow\Sigma^\infty\widetilde{E\mathscr{F}}$$
are maps between ring spectra, the functors
$$X\mapsto E\mathscr{F}_+\wedge X,\text{ }X\mapsto\widetilde{E\mathscr{F}}\wedge X$$
are symmetric monoidal.

Moreover, $E\mathscr{F}_+\wedge\widetilde{E\mathscr{F}}\simeq *$ since all its fixed point subspaces are contractible. Use $A,B$ to denote the suspensions of $E\mathscr{F}_+$ and $\widetilde{E\mathscr{F}}$ (in either order). The universal coefficient spectral sequence \cite{LM} tells us that for any $G$-spectra $X,Y$, we have
$$\underline{Ext}^{*,*}_{\underline{B}_*}(\underline{B}_*(A\wedge X),\underline{(B\wedge Y)}_*)\Rightarrow\underline{(B\wedge Y)}^*(A\wedge X)=[A\wedge X,\Sigma^* B\wedge Y]^G.$$
Since
$$\underline{B}_*(A\wedge X)=\underline{\pi}_*(B\wedge A\wedge X)=0,$$
the spectral sequence has trivial $E_2$-page. Thus we have

$$[E\mathscr{F}_+\wedge X,\widetilde{E\mathscr{F}}\wedge Y]^G=[\widetilde{E\mathscr{F}}\wedge Y,E\mathscr{F}_+\wedge X]^G=0.$$
$\Box$
\medskip

Now we can get a criterion about which prime factors of $|G|$ need to be inverted according to the proofs above:

\begin{remark} \label{invertprime}
It suffices to invert the prime factors of the denominators of all $c_H$, and the torsion in $\underline{H}_*^G(E\mathscr{F};A_G)$.

The numbers $c_H$ for all $H\in\mathscr{F}$ are defined in an inductive way by the equations
$$\sum_{[K]\leq[H]\subset\mathscr{F}}c_Hs_{(K,H)}=1,\text{ for any fixed }K\in\mathscr{F}.$$

According to \autoref{torsiontype}, the torsion in $\underline{H}_*^G(E\mathscr{F};A_G)$ comes from the torsion in $H_*(E\mathscr{F}^K/L)$ for all $K\lhd L\subset G$.
\end{remark}

\subsection{Splitting criterion}

In this subsection, we will improve \autoref{invertprime} and make a better criterion, which plays the key role in our splitting method:

\begin{theorem} \label{splitcrit}
The cofiber sequence in \autoref{ringsplit} splits after taking $\Sigma^\infty$ if and only if we invert all prime $p$ such that there exists $H\lhd J\subset G$ with $H\in\mathscr{F}$, $J\notin\mathscr{F}$, and $|J/H|$ is a power of $p$.
\end{theorem}

We prove the easier half first:
\medskip

\paragraph{Proof of ``only if" part:} Consider any prime $p$ such that the corresponding $H,J$ in the theorem exist. We can assume that $J/H\cong C_p$, otherwise we replace $H$ by larger proper normal subgroups of $J$. If the cofiber sequence
$$E\mathscr{F}_+\rightarrow S^0\rightarrow\widetilde{E\mathscr{F}}$$
splits after $\Sigma^\infty$, we have 
$$S\simeq\Sigma^\infty E\mathscr{F}_+\vee\Sigma^\infty\widetilde{E\mathscr{F}}.$$
Take equivariant homology with coefficients in $A_G$, we have
$$\widetilde{H}_*^G(S^0;A_G)=\widetilde{H}_*^G(E\mathscr{F}_+;A_G)\oplus\widetilde{H}_*^G(\widetilde{E\mathscr{F}};A_G).$$
Since $\widetilde{H}_*^G(S^0;A_G)\cong A_G$ concentrates in degree $0$, $\widetilde{H}_*^G(E\mathscr{F}_+;A_G)$ must concentrate in degree $0$. According to \autoref{burnsidehom}, $H_*(E\mathscr{F}^H/W_JH;\Z)$ is a direct summand of $\widetilde{H}_*^G(E\mathscr{F}_+;A_G)$, hence also concentrates in degree $0$. 

Since $H\in\mathscr{F}$, $J\notin\mathscr{F}$, $J/H\simeq C_p$, we have $H_*(E\mathscr{F}^H/W_JH;\Z)\simeq H_*(BC_p;\Z)$, which concentrates in degree $0$ only if $p$ is inverted. $\Box$
\bigskip

Now we consider the ``if" part. Assume that all primes described in \autoref{splitcrit} are inverted. According to \autoref{invertprime}, it suffices to prove:

\begin{proposition} \label{homtor}
For any $K\lhd L\subset G$ and $m>0$, $H_m(E\mathscr{F}^K/L)=0$.
\end{proposition}

\begin{proposition} \label{chcoe}
Consider numbers $c_H$ for all $H\in\mathscr{F}$ determined by the following equations: For any $K\in\mathscr{F}$, we have
$$\sum_H c_H\cdot|(G/H)^K|=1.$$
Here we choose one $H$ from each conjugacy class inside $\mathscr{F}$. Then the denominators of all $c_H$'s become invertible after we invert all prime factors of $|G|$ not in $T$.
\end{proposition}
\bigskip

\paragraph{Proof of \autoref{homtor}:} We divide the proof into three steps:
\medskip

\paragraph{Step 1:} Let $p$ be a prime factor of $|G|$ which is not inverted and $P$ be any Sylow $p$-subgroup of $G$. Then $|G/P|\cdot H_m(E\mathscr{F}/G)=0$ for all $m>0$.
\medskip

Fix any $G$-CW structure on $E\mathscr{F}$, which induces CW structures on both $E\mathscr{F}/P$ and $E\mathscr{F}/G$. The natural map $E\mathscr{F}/P\rightarrow E\mathscr{F}/G$ induces a map between cellular chain complexes:
$$\phi:C_*(E\mathscr{F}/P)\rightarrow C_*(E\mathscr{F}/G).$$

Define a map between abelian groups in the opposite direction
$$\psi:C_*(E\mathscr{F}/G)\rightarrow C_*(E\mathscr{F}/P)$$
as follows:

For any $G$-cell $G/H_+\wedge e$ of $E\mathscr{F}$, use $e_P$ and $e_G$ to denote the corresponding cells in $E\mathscr{F}/P$ and $E\mathscr{F}/G$. Let $Pg_1, Pg_2,...,Pg_k$ be all right cosets of $P$ in $G$. Define $\psi$ by sending each $e_G\in C_*(E\mathscr{F}/G)$ to $\sum_i (g_ie)_P\in C_*(E\mathscr{F}/P)$. 
\medskip

We can show that $\psi$ is well-defined by the following two facts:

\textbf{(1)} If $g,g^\prime$ are in the same right coset of $P$, then $(ge)_P=(g^\prime e)_P$. Thus $\psi$ does not depend on the choice of $g_1,g_2,...,g_k$.

\textbf{(2)} For any $e,e^\prime$, $e_G=e^\prime_G$ if and only if $e^\prime=ge$ for some $g\in G$. Since $g_1g,g_2g,...,g_kg$ also cover all right cosets of $P$, we have $\sum_i(g_i(ge))_P=\sum_i((g_ig)e)_P=\sum_i(g_ie)_P$. Thus $\psi$ is uniquely defined on each element of $C_*(E\mathscr{F}/G)$.
\bigskip

Moreover, $\psi$ commutes with the boundary maps:

We use $d,d_P,d_G$ to denote the cellular boundary maps for $E\mathscr{F},E\mathscr{F}/P,E\mathscr{F}/G$ respectively. For any $G$-cell $G/H_+\wedge e$ of $E\mathscr{F}$, write
$$de=\sum_t\left(\sum_{g\in G}n_gg\right)e_t=\sum_t\left[\sum_i\left(\sum_{a\in P} n_{ag_i}a\right)g_i\right]e_t$$
where the cells $e_t$ are in different $G$-orbits and $n_g\in\Z$. Then we have
$$d_Ge_G=\sum_t\left(\sum_{g\in P}n_g\right)(e_t)_G$$
$$d_Pe_P=\sum_{t,i}\left(\sum_{a\in P}n_{ag_i}\right)(g_ie_t)_P$$

Notice that
$$d(\sum_{g\in G}ge)=\sum_{g\in G}g(de)=\sum_t\left(\sum_{g\in G}n_g\right)\left(\sum_{g\in G}g\right)e_t=\sum_t\left[\sum_i\left(\sum_{a\in P}(\sum_{g\in G} n_g)a\right)g_i\right]e_t.$$
Thus we have
$$d_P(\sum_{g\in G}(ge)_P)=|P|\sum_{t,i}(\sum_{g\in G}n_g)(g_ie_t)_P.$$
Since $\sum_{g\in G}(ge)_P=|P|\sum_i(g_ie)_P$, we have
$$d_P(\sum_i(g_ie)_P)=\frac{1}{|P|}d_P(\sum_{g\in G}(ge)_P)=\sum_{t,i}(\sum_{g\in G}n_g)(g_ie_t)_P.$$
Therefore, we get
$$d_P(\psi e_G)=d_P(\sum_i(g_ie)_P)=\sum_{t,i}(\sum_{g\in G}n_g)(g_ie_t)_P$$
$$=\sum_t(\sum_{g\in G}n_g)(\sum_i(g_ie_t)_P)=\sum_t(\sum_{g\in G}n_g)\psi((e_t)_G)=\psi(d_Ge_G).$$
\medskip

Now $\psi$ becomes a morphism between chain complexes. Notice that $\phi\circ\psi$ is multiplication by $|G/P|$. Thus the induced map on homology
$$H_*(E\mathscr{F}/G)\xrightarrow{\psi_*}H_*(E\mathscr{F}/P)\xrightarrow{\phi_*}H_*(E\mathscr{F}/G)$$
is also multiplication by $|G/P|$.

Notice that since $p$ is not inverted, $\mathscr{F}$ must contain all $p$-subgroups of $G$; otherwise we can choose $H$ to be the trivial subgroup and $J$ to be the $p$-subgroup not in $\mathscr{F}$, which does not satisfy the condition in \autoref{splitcrit}. Thus $E\mathscr{F}$ is contractible as a $P$-space, and hence $H_m(E\mathscr{F}/P)=0$ for all $m>0$. The maps above imply that $|G/P|\cdot H_m(E\mathscr{F}/G)=0$ for all $m>0$.
\bigskip

\paragraph{Step 2:} Now we list all prime factors of $|G|$ which are not inverted: $p_1,p_2,...,p_n$.

For any $i\in\{1,2,...,n\}$, let $P_i$ be a Sylow-$p_i$ subgroup of $G$. By step 1, we have $|G/P_i|\cdot H_m(E\mathscr{F}/G)=0$ for all $m>0$.
\medskip

Letting $i$ vary among $1,2,...,n$, we get
$$(|G/P_1|,|G/P_2|,...,|G/P_n|)\cdot H_m(E\mathscr{F}/G)=0$$
for $m>0$. Here $(|G/P_1|,|G/P_2|,...,|G/P_n|)$ is the greatest common divisor of all $|G/P_i|$. Since all prime factors of $G$ except $p_1,p_2,...,p_n$ are inverted, $(|G/P_1|,|G/P_2|,...,|G/P_n|)$ is inverted. Thus $H_m(E\mathscr{F}/G)=0$ for any $m>0$.
\bigskip

\paragraph{Step 3:} Consider the general case of any $K\lhd L\subset G$. By comparing the fixed point subspaces, $E\mathscr{F}^K$ is $L/K$-equivalent to $E\mathscr{G}$, where $\mathscr{G}$ is the following family of subgroups of $L/K$:
$$\mathscr{G}:=\{H/K\text{ }|\text{ }K\subset H\subset L, H\in\mathscr{F}\}.$$
Consider any $H_1,H_2\subset G$ such that $H_1/K\in\mathscr{G}$, $H_2/K\notin\mathscr{G}$, $(H_1/K)\lhd (H_2/K)$ and $|(H_2/K)/(H_1/K)|$ is a power of prime $p$. We must have $H_1\in\mathscr{F}$, $H_2\notin\mathscr{F}$, $H_1\lhd H_2$ and $|H_2/H_1|$ is a power of $p$. Thus $p$ must be inverted. By induction on the group size and the previous two steps, we have $H_m(E\mathscr{F}^K/L)\cong H_m(E\mathscr{G}/(L/K))=0$ for all $m>0$. $\Box$
\bigskip

We now prove \autoref{chcoe}. Let $T=\{p_1,...,p_n\}$ be the set of all prime factors of $|G|$ which are not inverted.

For any $H\subset G$, recall the ring map $\chi_H: A(G)\rightarrow\Z$ defined in the proof of \autoref{homsplit} which sends each $G$-set $S$ to $|S^H|$. Let $C(G):=\prod_{[H]}\Z$. Then we have a ring map
$$\chi=\prod_{[H]}\chi_H:A(G)\rightarrow C(G).$$
The proof of \autoref{homsplit} implies that $\chi$ is a ring isomorphism when $|G|$ is inverted. If we only invert the prime factors not in $T$, $\chi$ becomes a monomorphism and $C(G)/Im(\chi)$ only contains torsion of $p_1,...,p_n$.
\medskip

We want to show that the denominators of all $c_H$'s are invertible, which is equivalent to the existence of the element $\sum_H c_H\{G/H\}$ in $A(G)$, when all prime factors of $|G|$ not in $T$ are inverted. Since $\chi$ is a monomorphism, this existence can be detected by the image of $\chi$:
$$\chi_K(\sum_H c_H\{G/H\})=\sum_H c_H|(G/H)^K|=
\begin{cases}
1,\text{ if }K\in\mathscr{F}\\
0,\text{ otherwise}
\end{cases}$$

In other words, it suffices to show that the element $\alpha=(\alpha_H)_{[H]}\in C(G)$, which is defined as
$$\alpha_H=
\begin{cases}
1,\text{ if }K\in\mathscr{F},\\
0,\text{ otherwise,}
\end{cases}$$
is inside the image of $\chi$.
\bigskip

We will need the following lemma:

\begin{lemma} \label{keylemma}
For any $H\in\mathscr{F}$ and $p\in T=\{p_1,p_2,...,p_n\}$, there exists $L\in\mathscr{F}$, such that $\chi_H(G/L)$ is not divided by $p$.
\end{lemma}

\paragraph{Proof:} Let $M$ be the subgroup of $H$ generated by all elements whose order is not divided by $p$.

List all prime factors of $|G|$ other than $p$ as $q_1,q_2,...,q_m$. Notice that $M$ is generated by elements $h\in G$ whose order has the form $|h|=q_1^{a_1}q_2^{a_2}...q_m^{a_m}$. Since $h$ is contained in the group generated by $h^{|h|/q_i^{a_i}}$ for all $i=1,2,...,m$ and each $h^{|h|/q_i^{a_i}}$ has a $q_i$-power order, $M$ is also generated by its Sylow $q_i$-subgroups, for $i=1,2,...,m$.
\smallskip

We observe that $M$ is a normal subgroup of $H$ since conjugations do not change the order of an element. Thus $H\subset N_GM$.

Write $H=Q\ltimes M$ and $N_GM=W_GM\ltimes M$ with $Q\subset W_GM$.
\medskip

For any $h\in H$, write the order of $h$ as $|h|=p^ac$ with $c$ not divided by $p$. Then $|h^{p^a}|=c$ and hence $h^{p^a}\in M$. Therefore, $H/M=Q$ is a $p$-group.

Choose $P\in W_GM$ to be a Sylow-$p$ group containing $Q$. Let $L=P\ltimes M\supset H$. Since $M\in\mathscr{F}$ as a subgroup of $H$ and $L/M\simeq P$ has a $p$-power order, we have $L\in\mathscr{F}$. We will show that $\chi_H(G/L)$ is not divided by $p_1$.
\bigskip

Recall that $\chi_H(G/L)=|(G/L)^H|$. For any $g\in G$, $gL\in G/L$ is fixed by $H$ if and only if $g^{-1}Hg\subset L$. If this happens, then $g^{-1}Mg\subset L=P\ltimes M$.

Let $N$ be any Sylow $q_i$-subgroup of $M$ for some $i\in\{1,2,...,m\}$. Then $g^{-1}Ng\subset g^{-1}Mg\subset L$. Since $|L|=|P|\cdot|M|$ and $|P|$ is a $p$-power, both $N$ and $g^{-1}Ng$ are Sylow-$q_i$ subgroups of $L$, hence are conjugate to each other inside $L$. Since $M$ is a normal subgroup of $L$ and $N\subset M$, we get $g^{-1}Ng\subset M$. Since $M$ is generated by its Sylow-$q_1,...,q_m$ subgroups, we have $g^{-1}Mg=M$, and hence $g\in N_GM$.
\medskip

When $gL\in(G/L)^H$, we have $g\in N_GM$ and $gL\in(N_GM/L)^H$. Thus
$$\chi_H(G/L)=|(G/L)^H|=|(N_GM/L)^H|=|(W_GM/P)^Q|.$$
Recall that $P$ is a Sylow-$p$ group of $W_GM$ containing $Q$. Thus $|W_GM/P|$ is not divisible by $p$. Since the sizes of all $Q$-orbits except $Q/Q$ are divided by $p$, the number of $Q/Q$ orbits in $W_GM/P$ is not divided by $p$. Thus $|(W_GM/P)^Q|$ is not divided by $p$.
\medskip

In conclusion, $\chi_H(G/L)$ is not divided by $p$. $\Box$
\bigskip

\paragraph{Proof of \autoref{chcoe}:} Since $C(G)$ has finite rank and $C(G)/Im(\chi)$ only contains torsion of $p_1,...,p_n$, there exist $k_1,k_2,...,k_n\in\Z^+$, such that $p_1^{k_1}...p_n^{k_n}C(G)\subset Im(\chi)$.
\medskip

We divide our proof into three steps:
\medskip

\textbf{Step 1:} For any $H\in\mathscr{F}$ and $i\in\{1,...,n\}$, there exists an element $\beta(H,i)\in Im(\chi)$, such that, modulo $p_i^{k_i}$, we have
$$\beta(H,i)_K\equiv
\begin{cases}
1,\text{ if }K\in[H]\\
1\text{ or }0,\text{ if }K\in\mathscr{F}\\
0,\text{ if }K\notin\mathscr{F}
\end{cases}$$
According to \autoref{keylemma}, there exists $L\in\mathscr{F}$ such that $\chi_H(G/L)$ is not divided by $p_i$. Notice that $\chi_K(G/L)=0$ for any $K\notin\mathscr{F}$. Since $|(\Z/p_i^{k_i})^\times|=p_i^{k_i-1}(p_i-1)$, the image of $\{G/L\}^{p_i^{k_i-1}(p_i-1)}$ is a valid choice for $\beta(H,i)$.
\medskip

\textbf{Step 2:} For any $H\in\mathscr{F}$, there exists an element $\beta(H)\in Im(\chi)$, such that, modulo $p_1^{k_1}...p_n^{k_n}$, we have
$$\beta(H)_K\equiv
\begin{cases}
1,\text{ if }K\in[H]\\
1\text{ or }0,\text{ if }K\in\mathscr{F}\\
0,\text{ if }K\notin\mathscr{F}
\end{cases}$$
For $i=1,2,...,n$, choose $c_i\in\Z$ such that $c_1\equiv 1$ (mod $p_i^{k_i}$) and $c_1\equiv 0$ (mod $p_j^{k_j}$) for all $j\neq i$. We only need to choose
$$\beta(H)=c_1\beta(H,1)+...+c_n\beta(H,n).$$
\smallskip

\textbf{Step 3:} Define $a+_*b=a+b-ab$ for any $a,b$ in the same ring. Notice that for any $a,b,c\in\Z$ such that $a,b\equiv 0$ or $1$ (mod $c$), we have
$$a+_*b\equiv
\begin{cases}
0,\text{ if }a\equiv b\equiv 0\\
1,\text{ otherwise}
\end{cases}$$
mod $c$.

Choose one subgroup from each conjugacy class inside $\mathscr{F}$ and list them as $H_1,H_2,...,H_m$. Let
$$\beta=\beta(H_1)+_*\beta(H_2)+_*...+_*\beta(H_m)$$
Then modulo $p_1^{k_1}...p_n^{k_n}$, we have
$$\beta_K\equiv
\begin{cases}
1,\textbf{ if }K\in\mathscr{F}\\
0,\textbf{ if }K\notin\mathscr{F}
\end{cases}$$

Since $\beta\in Im(\chi)$ and $p_1^{k_1}...p_n^{k_n}\cdot C(G)\subset Im(\chi)$, the special element $\alpha=(\alpha_H)_{[H]}$, which is defined as
$$\alpha_H=
\begin{cases}
1,\text{ if }H\in\mathscr{F},\\
0,\text{ if }H\notin\mathscr{F},
\end{cases}$$
must be contained in $Im(\chi)$. $\Box$

\section{Splitting method}

In this section, we will set up our splitting method in the computation for any equivariant spectrum $X$.
\medskip

The procedure of the splitting method can be explained as follows:
\bigskip

\textbf{(I)} Fix a prime factor $p$ of $|G|$. Invert all other prime factors.
\medskip

\textbf{(II)} Choose a specific family $\mathscr{F}$ such that the criterion in \autoref{splitcrit} is satisfied.
\medskip

\textbf{(III)} Apply \autoref{ringsplit} to split $X$ into multiple pieces. Compute the homotopy of each piece and glue them together. Now we have computed the localized $\pi_\bigstar^G(X)$ when all prime factors of $|G|$ except $p$ are inverted.
\medskip

\textbf{(IV)} List all prime factors of $|G|$ as $p_1,p_2,...,p_n$. For $i=1,2,...,n$, we can compute
$$R_i:=\pi_\bigstar^G(X)[p_1^{-1},...,\widehat{p_i^{-1}},...,p_n^{-1}]$$
by the three steps above with $p=p_i$. 

Consider the diagram with objects $R_1,R_2,...,R_n, \pi_\bigstar^G(HM)[|G|^{-1}]$ and morphisms
$$R_i\rightarrow \pi_\bigstar^G(X)[|G|^{-1}],\text{ }i=1,2,...,n$$
as localizations. This is a diagram inside the category of $RO(G)$-graded rings. The limit of the diagram is exactly $\pi_\bigstar^G(X)$.
\bigskip

The hardest part is the computation of each piece in step (III), which will be discussed in \autoref{moregeneralefsplit}.
\medskip

First we explain the choice of family in step (II):
\medskip

Let $H$ be the subgroup of $G$ generated by all elements whose orders are not divisible by $p$. Choose $\mathscr{F}$ to be the family of all subgroups of $G$ not containing $H$. The criterion in \autoref{splitcrit} is satisfied because of the following lemma:

\begin{lemma} \label{familychoice}
$H$ is a normal subgroup of $G$ with $p$-power index and $H$ is generated by Sylow $q$-subgroups of $G$ for all $q\neq p$. Moreover, $H$ itself does not have any proper normal subgroups with $p$-power index.
\end{lemma}

\paragraph{Proof:} The first half can be proved in a similar way as in the proof of \autoref{keylemma}.

For the second half, assume that $M$ is a normal subgroup of $H$ with $p$-power index. Then for any $q\neq p$, the Sylow $q$-subgroup of $M$ is also a Sylow $q$-subgroup of $H$. Since $M$ is normal in $H$, it contains all Sylow $q$-subgroups of $H$ for any $q\neq p$. Thus $M=H$ since $H$ is generated by these Sylow subgroups. $\Box$
\bigskip

We have $X\simeq (E\mathscr{F}_+\wedge X)\vee (\widetilde{E\mathscr{F}}\wedge X)$ when \autoref{ringsplit} holds. With this specific choice of $\mathscr{F}$, we have
\begin{equation} \label{geosplitpart}
\pi_V^G(\widetilde{E\mathscr{F}}\wedge X)\cong\pi_{V^H}^H(\Phi^HX)
\end{equation}
where $\Phi^HX$ is the $H$-geometric fixed point spectrum of $X$.

\subsection{Decomposition of $E\mathscr{F}_+\wedge X$}

When the splitting happens, we can compute the homotopy of $E\mathscr{F}_+\wedge X$ in a systematic way. Recall that $M_\mathscr{F}$ is defined in \autoref{mfnf} as the homology of $E\mathscr{F}_+$ in degree $0$.

\begin{proposition} \label{eftomap}
Let $\mathscr{F}$ be an arbitrary family. Invert proper prime factors of $|G|$ such that \autoref{ringsplit} holds. Then for any $G$-spectrum $X$, we have
$$\pi_V^G(E\mathscr{F}_+\wedge X)\cong Map(M_\mathscr{F},\underline{\pi}_V^G(X)).$$
Here $\underline{\pi}_V^G(X)$ is the Mackey functor valued homotopy which sends each $G/H$ to $\pi_{V|_H}^H(X)$.
\end{proposition}

Together with the expression of $M_\mathscr{F}$ in \autoref{homsplit}, we have:

\begin{theorem} \label{moregeneralefsplit}
Let $\mathscr{F}$ be an arbitrary family of subgroups of $G$. Invert proper prime factors of $|G|$ such that \autoref{ringsplit} holds.

For any virtual $G$-representation $V$ and $G$-spectrum $X$, $\pi_V^G(E\mathscr{F}_+\wedge X)$ is isomorphic to the limit of the diagram with objects
$$\pi_{V|_H}^H(X),\text{ }H\in\mathscr{F}$$
and morphisms 
$$res_{H_1}^{H_2}:\pi_{V|_{H_2}}^{H_2}(X)\rightarrow \pi_{V|_{H_1}}^{H_1}(X),\text{ }H_1\subset H_2\in\mathscr{F},$$
$$c_g:\pi_{V|_H}^H(X)\rightarrow \pi_{V|_{gHg^{-1}}}^{gHg^{-1}}(X),\text{ }H\in\mathscr{F},g\in G.$$

The map from $\pi_V^G(E\mathscr{F}_+\wedge X)$ to each object $\pi_{V|_H}^H(X)$ can be factorized as
$$\pi_V^G(E\mathscr{F}_+\wedge X)\rightarrow \pi_V^G(X)\xrightarrow{res_H^G}\pi_{V|_H}^H(X).$$

Moreover, when $X$ is a ring spectrum, since taking limits preserves multiplicative structures, this decomposition can also be used to compute $\pi_\bigstar^G(E\mathscr{F}_+\wedge X)$ as a graded ring.
\end{theorem}

\begin{remark}
The image of $res_P^G$ in $\pi_{V|_P}^P(X)$ does not depend only on $V|_P$ since the choice of $g$ such that $g^{-1}Hg\subset P$ may not be contained in $P$, not even in $N_GP$. Thus it's possible that for two different $G$-representations $V,W$, we have $V|_P\cong W|_P$ and yet $\pi_V^G(E\mathscr{F}_+\wedge X)$ and $\pi_W^G(E\mathscr{F}_+\wedge X)$ are different.
\end{remark}
\medskip

\paragraph{Proof of \autoref{moregeneralefsplit}:} According to \autoref{eftomap}, we have
$$\pi_V^G(E\mathscr{F}_+\wedge X)\cong Map(M_{\mathscr{F}},\underline{\pi}_V^GX).$$

According to \autoref{homsplit}, as a sub-Mackey functor of $A_G$, $M_\mathscr{F}$ is generated by $\{H/H\}$ for all $H\in\mathscr{F}$, with relations generated by
$$res_{H_1}^{H_2}(\{H_2/H_2\})=\{H_1/H_1\},\text{ }H_1\subset H_2\in\mathscr{F},$$
$$c_g(\{H/H\})=\{gHg^{-1}/gHg^{-1}\},\text{ }H\in\mathscr{F},g\in G.$$

Each map from $M_\mathscr{F}$ to $\pi_V^G(E\mathscr{F}_+\wedge X)$ is determined by the images of $\{H/H\}$ in
$$\underline{\pi}_V^G(X)(G/H)=\pi_{V|_H}^H(X)$$
for all $H\in\mathscr{F}$ which are compatible with maps $res_{H_1}^{H_2}$ and $c_g$. Thus $\pi_V^G(E\mathscr{F}_+\wedge X)$ is expressed as the limit described in the theorem. $\Box$
\bigskip

In order to prove \autoref{eftomap}, we need a small lemma on the homotopy of $E\mathscr{F}$ first. Recall that $N_\mathscr{F}$ is defined in \autoref{mfnf} as the complement of $M_\mathscr{F}$ in $A_G$.

\begin{lemma} \label{efhomotopy}
Assume that proper primes are inverted such that \autoref{ringsplit} holds. Then
$$\underline{\pi}_\bigstar^G(\Sigma^\infty E\mathscr{F}_+)\cong\underline{\pi}_\bigstar^G(S)\text{ }\Box\text{ }M_\mathscr{F},$$
$$\underline{\pi}_\bigstar^G(\Sigma^\infty\widetilde{E\mathscr{F}})\cong\underline{\pi}_\bigstar^G(S)\text{ }\Box\text{ }N_\mathscr{F}.$$
Here $S=\Sigma^\infty S^0$ is the sphere spectrum.
\end{lemma}

\paragraph{Proof:} According to \autoref{ringsplit}, we have an equivalence of ring spectra:
$$S\simeq \Sigma^\infty E\mathscr{F}_+\vee\Sigma^\infty\widetilde{E\mathscr{F}}.$$
Thus the product in $\underline{\pi}_\bigstar^G(S)$ is induced by the products in $\underline{\pi}_\bigstar^G(\Sigma^\infty E\mathscr{F}_+)$ and $\underline{\pi}_\bigstar^G(\Sigma^\infty\widetilde{E\mathscr{F}})$.

In particular, the map
$$\underline{\pi}_\bigstar^G(S)\text{ }\Box\text{ }\underline{\pi}_0^G(S)\rightarrow \underline{\pi}_\bigstar^G(S)$$
is given by
$$\underline{\pi}_\bigstar^G(S)\text{ }\Box\text{ }\underline{\pi}_0^G(S)\cong\left(\underline{\pi}_\bigstar^G(\Sigma^\infty E\mathscr{F}_+)\text{ }\Box\text{ }\underline{\pi}_0^G(\Sigma^\infty E\mathscr{F}_+)\right)\oplus\left(\underline{\pi}_\bigstar^G(\Sigma^\infty \widetilde{E\mathscr{F}})\text{ }\Box\text{ }\underline{\pi}_0^G(\Sigma^\infty \widetilde{E\mathscr{F}})\right)$$
$$\oplus\left(\underline{\pi}_\bigstar^G(\Sigma^\infty E\mathscr{F}_+)\text{ }\Box\text{ }\underline{\pi}_0^G(\Sigma^\infty \widetilde{E\mathscr{F}})\right)\oplus\left(\underline{\pi}_\bigstar^G(\Sigma^\infty \widetilde{E\mathscr{F}})\text{ }\Box\text{ }\underline{\pi}_0^G(\Sigma^\infty E\mathscr{F}_+)\right)\xrightarrow{f}$$
$$\left(\underline{\pi}_\bigstar^G(\Sigma^\infty E\mathscr{F}_+)\text{ }\Box\text{ }\underline{\pi}_0^G(\Sigma^\infty E\mathscr{F}_+)\right)\oplus\left(\underline{\pi}_\bigstar^G(\Sigma^\infty \widetilde{E\mathscr{F}})\text{ }\Box\text{ }\underline{\pi}_0^G(\Sigma^\infty \widetilde{E\mathscr{F}})\right)\rightarrow$$
$$\underline{\pi}_\bigstar^G(\Sigma^\infty E\mathscr{F}_+)\oplus\underline{\pi}_\bigstar^G(\Sigma^\infty \widetilde{E\mathscr{F}})\cong\underline{\pi}_\bigstar^G(S).$$
Here the map $f$ is the projection onto the first two summands.
\smallskip

On the other hand, the map
$$\underline{\pi}_\bigstar^G(S)\text{ }\Box\text{ }\underline{\pi}_0^G(S)\rightarrow \underline{\pi}_\bigstar^G(S)$$
is an isomorphism since $\underline{\pi}_0^G(S)\cong A_G$ is the unit Mackey functor. Thus the projection $f$ must be an isomorphism and hence
$$\underline{\pi}_\bigstar^G(\Sigma^\infty E\mathscr{F}_+)\text{ }\Box\text{ }\underline{\pi}_0^G(\Sigma^\infty \widetilde{E\mathscr{F}})=\underline{\pi}_\bigstar^G(\Sigma^\infty \widetilde{E\mathscr{F}})\text{ }\Box\text{ }\underline{\pi}_0^G(\Sigma^\infty E\mathscr{F}_+)=0.$$

Since $\Sigma^\infty E\mathscr{F}_+$ and $\Sigma^\infty \widetilde{E\mathscr{F}}$ are connective, we have
$$\underline{\pi}_0^G(\Sigma^\infty  E\mathscr{F}_+)\cong\widetilde{\underline{H}}_0^G(E\mathscr{F}_+;A_G)=M_\mathscr{F},$$
$$\underline{\pi}_0^G(\Sigma^\infty \widetilde{E\mathscr{F}})\cong\widetilde{\underline{H}}_0^G(\widetilde{E\mathscr{F}};A_G)=N_\mathscr{F}.$$
Thus
$$\underline{\pi}_\bigstar^G(\Sigma^\infty E\mathscr{F}_+)\text{ }\Box\text{ }N_\mathscr{F}=\underline{\pi}_\bigstar^G(\Sigma^\infty \widetilde{E\mathscr{F}})\text{ }\Box\text{ }M_\mathscr{F}=0$$
and hence
$$\underline{\pi}_\bigstar^G(\Sigma^\infty E\mathscr{F}_+)\cong\underline{\pi}_\bigstar^G(S)\text{ }\Box\text{ }M_\mathscr{F},$$
$$\underline{\pi}_\bigstar^G(\Sigma^\infty\widetilde{E\mathscr{F}})\cong\underline{\pi}_\bigstar^G(S)\text{ }\Box\text{ }N_\mathscr{F}.$$
$\Box$

\paragraph{Proof of \autoref{eftomap}:} According to \autoref{ringsplit}, we have
$$\pi_V^G(E\mathscr{F}_+\wedge X)=[S^V,E\mathscr{F}_+\wedge X]^G\cong[E\mathscr{F}_+,\Sigma^{-V}X]^G=X^{-V}E\mathscr{F}_+.$$
Here we express the homotopy of $E\mathscr{F}_+\wedge X$ as the $X$-cohomology of the space $E\mathscr{F}_+$.

View $X$ as an $S$-module and apply the equivariant universal coefficient spectral sequence in \cite{LM}:
$$Ext^{*,*}_{\underline{\pi}_*^G(S)}(\underline{\pi}_*^G(\Sigma^\infty E\mathscr{F}_+),\underline{X}^{*-V})\Rightarrow X^{*-V}E\mathscr{F}_+.$$

According to \autoref{efhomotopy} and \autoref{homsplit}, we have
$$Ext^{*,*}_{\underline{\pi}_*^G(S)}(\underline{\pi}_*^G(\Sigma^\infty E\mathscr{F}_+),\underline{X}^{*-V})\cong Ext^{*,*}_{\underline{\pi}_*^G(S)}(\underline{\pi}_*^G(S)\text{ }\Box\text{ }M_\mathscr{F},\underline{X}^{*-V})$$
$$\cong Ext^{*,*}_{A_G}(M_\mathscr{F},\underline{X}^{*-V})\cong Map(M_\mathscr{F},\underline{\pi}_V^GX)$$
and the $E_2$-page is trivial except on the bottom line. Thus we have
$$\pi_V^G(E\mathscr{F}_+\wedge X)=[S^V,E\mathscr{F}_+\wedge X]^G\cong
\smallskip Map(M_\mathscr{F},\underline{\pi}_V^G(X)).$$
$\Box$

\subsection{Input of the computation}

When we apply the decomposition into smaller subgroups by \autoref{moregeneralefsplit} and \autoref{geosplitpart}, if the subgroups still contain multiple prime factors, we can apply the splitting method again. Finally everything will be decomposed into equivariant homotopy for subgroups with prime power orders. We will give a precise list of input of this computational machine. We need the following notation:
\smallskip

For any finite group $H$, let $H_p$ be the subgroup of $H$ generated by all elements whose orders are not divisible by $p$. For any $K\subset H$, let $\mathscr{F}_H[K]$ be the family of subgroups of $H$ not containing $K$. Then our splitting method uses the family $\mathscr{F}_H[H_p]$ to split $H$-spectra when localized at $p$.

\begin{theorem} \label{machineinput}
Let $X$ be a $G$-spectrum and $V$ be a virtual $G$-representation. Assume that all prime factors of $|G|$ except $p$ are inverted. Then $\pi_V^G(X)$ is the limit of the diagram consisting of:

\textbf{(i)} $\pi_{V|_H}^H(X)$ for any $p$-subgroup $H\subset G$;

\textbf{(ii)} $\pi_{(V|_K)^{K_p}}^{K/K_p}(\Phi^{K_p}X)$ for any $K\in\mathscr{F}_G[G_p]\cup\{G\}$ such that $K$ is not a $p$-group;

\textbf{(iii)} Restriction map $\pi_{V|_{H_2}}^{H_2}(X)\rightarrow\pi_{V|_{H_1}}^{H_1}(X)$ for any $H_1\subset H_2$;

\textbf{(iv)} Restriction map $\pi_{(V|_L)^{L_p}}^{L/L_p}(\Phi^{L_p}X)\rightarrow\pi_{(V|_K)^{K_p}}^{K/K_p}(\Phi^{K_p}X)$ for any $K\subset L$ such that $K_p=L_p$;

\textbf{(v)} Conjugacy map on $\pi_{V|_H}^H(X)$ by any $g\in W_GH$;

\textbf{(vi)} Conjugacy map on $\pi_{(V|_K)^{K_p}}^{K/K_p}(\Phi^{K_p}X)$ by any $g\in W_GK$ (notice that $g^{-1}K_pg=K_p$ according to the definition of $K_p$).
\end{theorem}
\medskip

We will prove this theorem in the rest of this subsection. Notice that \autoref{catmain} is a consequence of this theorem.
\bigskip

Before we start the proof, we need a more precise explanation about how to repeat the splitting method. After we apply the splitting method on the $G$-spectrum $X$ with family $\mathscr{F}_G[G_p]$, we obtain a diagram with objects $\pi_{V|_H}^H(X)$, for all $H$ in $\mathscr{F}_G[G_p]$ and $\pi_{V^{G_p}}^{G/G_p}(\Phi^{G_p}X)$. We don't have to further split the geometric fixed point spectrum since it is acted on by the $p$-group $G/G_p$. But we can apply the splitting method again on each $\pi_{V|_H}^H(X)$ with family $\mathscr{F}_H[H_p]$ to make it into the limit of another diagram.

The conjugacy map on $\pi_{V|_H}^H(X)$ splits in the same way since the splitting in \autoref{ringsplit} is categorical. The following lemma implies that the restriction maps are also compatible with the splittings. Thus repeating the splitting method not only expresses each object as a limit, but also enlarge the whole diagram.

By applying \autoref{moregeneralefsplit} and tracing the restriction maps, we have

\begin{lemma} \label{splitres}
Consider any $H\subset L\subset G$. Assume that all prime factors of $|G|$ except $p$ are inverted. The restriction map $res_H^L:\pi_{V|_L}^L(X)\rightarrow\pi_{V|_H}^H(X)$ can be expressed as follows:

\textbf{(1)} If $H\in\mathscr{F}_L[L_p]$, or equivalently, $H$ does not contain $L_p$, then after applying the splitting on $\pi_{V|_L}^L(X)$ with family $\mathscr{F}_L[L_p]$, $res_H^L$ can be factorized as
$$\pi_{V|_L}^L(X)\rightarrow\pi_{V|_L}^L(E\mathscr{F}_L[L_p]_+\wedge X)\rightarrow\pi_{V|_H}^H(X)$$
since $\pi_{V|_H}^H(X)$ is an object in the diagram computing $\pi_{V|_L}^L(E\mathscr{F}_L[L_p]_+\wedge X)$ in \autoref{moregeneralefsplit}.

\textbf{(2)} If $H\notin\mathscr{F}_L[L_p]$, or equivalently, $H$ contains $L_p$ and hence $H_p=L_p$, we can apply the splitting on both $\pi_{V|_L}^L(X)$ and $\pi_{V|_H}^H(X)$ with families $\mathscr{F}_L[L_p]$ and $\mathscr{F}_H[H_p]$. Then $res_H^L$ is split into two pieces:
$$\pi_{V|_L}^L(E\mathscr{F}_L[L_p]_+\wedge X)\rightarrow\pi_{V|_H}^H(E\mathscr{F}_H[H_p]_+\wedge X)$$
$$\pi_{(V|_L)^{L_p}}^{L/L_p}(\Phi^{L_p}X)\rightarrow\pi_{(V|_H)^{H_p}}^{H/H_p}(\Phi^{H_p}X)$$
The second map is obtained by taking the underlying $H$-spectrum. The first map is obtained by the fact that the diagram computing $\pi_{V|_H}^H(E\mathscr{F}_H[H_p]_+\wedge X)$ in \autoref{moregeneralefsplit} is a sub-diagram of the diagram computing $\pi_{V|_L}^L(E\mathscr{F}_L[L_p]_+\wedge X)$.
\end{lemma}

\begin{remark}
The splittings on the underlying spectra can be applied in the following way to keep enlarging the diagram:

Choose an object $\pi_{V|_H}^H(X)$ such that $H_p$ has the largest size (if the largest size is $1$, no more splittings are required). Apply the splitting on all $\pi_{V|_K}^K(X)$ such that $K_p=H_p$. \autoref{splitres} guarantees that all restriction maps are transferred into the new diagram. Thus the diagram is enlarged in each such round of splitting.
\end{remark}
\medskip

Now we can repeat the spiltting method to enlarge the diagram as much as possible. Notice that when $H$ is a $p$-group, $H_p$ is trivial and hence the splitting does not change the $H$-homotopy. In other words, the diagram will finally stabilize when each summand is acted on by a $p$-group.
\medskip

\paragraph{Proof of \autoref{machineinput}:} After the diagram stabilizes, the argument above and \autoref{splitres} shows that all objects and maps will have the forms in (i)-(vi). On the other hand, all objects and maps in (i)-(vi) appear after the splitting is applied once on all underlying spectra where the acting group is not a $p$-group, and don't change in further splittings, hence survive in the stabilized diagram. $\Box$

\subsection{Computability of $H\underline{\Z}$}

In this subsection, we will discuss the computability of Eilenberg-Maclane spectra with constant Mackey functors, with $H\underline{\Z}$ the universal choice among them. Recall that the splitting method decomposes $G$-spectra into underlying and geometric fixed point spectra. $H\underline{\Z}$ is quite special because the geometric fixed point spectra will be trivial when the splitting happens:

\begin{proposition} \label{trivialtilde}
If \autoref{ringsplit} holds, we have $\widetilde{E\mathscr{F}}\wedge H\underline{\Z}\simeq*$. Thus
$$H\underline{\Z}\simeq E\mathscr{F}_+\wedge H\underline{\Z}.$$
\end{proposition}

\paragraph{Proof:} Consider the equivariant universal coefficient spectral sequence \cite{LM}:
$$E^2_{*,*}=\underline{Tor}_{*,*}^{A_G}(\underline{HA_G}_*\widetilde{E\mathscr{F}},\underline{\Z})\Rightarrow \underline{H\underline{\Z}}_*\widetilde{E\mathscr{F}}.$$
According to the proof of \autoref{homsplit}, $\underline{HA_G}_*\widetilde{E\mathscr{F}}$ is concentrated in degree $0$ and appears as the direct summand $N_\mathscr{F}$ of $A_G$. Thus the $E^2$-page is trivial except 
$$E^2_{0,0}\cong\underline{\Z}\text{ }\Box\text{ }N_\mathscr{F}.$$

Since
$$\underline{\Z}=\underline{\Z}\text{ }\Box\text{ }A_G=\underline{\Z}\text{ }\Box\text{ }M_\mathscr{F}\oplus\underline{\Z}\text{ }\Box\text{ }N_\mathscr{F},$$
there must be one $\underline{\Z}$ and one $0$ in $\{\underline{\Z}\text{ }\Box\text{ }M_\mathscr{F}, \underline{\Z}\text{ }\Box\text{ }N_\mathscr{F}\}$.

Notice that $M_\mathscr{F}(G/e)=H_0(E\mathscr{F};\Z)=\Z$. Thus
$$\underline{\Z}\text{ }\Box\text{ }M_\mathscr{F}(G/e)=\underline{\Z}(G/e)\otimes M_\mathscr{F}(G/e)=\Z.$$
Now we get
$$\underline{\Z}\text{ }\Box\text{ }M_\mathscr{F}=\underline{\Z},\text{ }\underline{\Z}\text{ }\Box\text{ }N_\mathscr{F}=0.$$
Therefore, the $E^2$-page of the spectral sequence is trivial, hence we have $\widetilde{E\mathscr{F}}\wedge H\underline{\Z}\simeq *$. $\Box$
\bigskip

\autoref{trivialtilde} still works when $\underline{\Z}$ is replaced by other constant Mackey functors, as can be proved by considering the induced maps from $\underline{\Z}$.
\bigskip

We now list the currently known computations on $p$-groups that can be fed in as inputs in \autoref{machineinput}:
\smallskip

\textbf{(i)} $\pi_\bigstar^{C_p}(H\underline{\Z})$, which is computed by multiple people with different methods.

\textbf{(ii)} $\pi_\bigstar^{C_{p^2}}(H\underline{\Z})$ by Zeng \cite{Zen} for odd $p$ and Georgakopoulos \cite{Geo} for $p=2$.

\textbf{(iii)} $\pi_\bigstar^{C_2^2}(H\underline{\F_2})$ and partial information about $\pi_\bigstar^{C_2^2}(H\underline{\Z})$ by Ellis-Bloor \cite{Ell}.

\textbf{(iv)} $\pi_\bigstar^{Q_8}(H\underline{\Z})$ by Lu \cite{Lu}, only with additive structure.
\medskip

The computation which works for $H\underline{\Z}$ usually works for other constant Mackey functors as well.
\bigskip

Most of the computations above are done in a cellular way. Thus we can obtain the restriction and conjugacy maps by treating the homotopy of $H\underline{\Z}$ as the equivariant homology of the representation spheres.
\medskip

In conclusion, now we have a computability list with our splitting method:

\begin{theorem} \label{grouplist}
For a finite group $G$, $\pi_\bigstar^G(H\underline{\Z})$ is computable if all Sylow subgroups $P\subset G$ satisfy one of the following conditions:

\textbf{(i)} $P=C_p$ or $C_{p^2}$;

\textbf{(ii)} $P=C_2^2$ and $P\subset A_4\subset W_GP$ (this guarantees that we only need the partial information in \cite{Ell});

\textbf{(iii)} $P=Q_8$ (we can only compute the additive structure in this case).
\end{theorem}
\medskip

\begin{remark}
When $\underline{\Z}$ is replaced by a general Mackey functor $M$, the geometric fixed point part is still computable since we have
$$\widetilde{E\mathscr{F}}\wedge HM\simeq H(M\text{ }\Box\text{ }N_\mathscr{F}).$$
However, we don't know how to compute the homotopy of general $HM$ for $p$-groups other than $C_p$. Thus $\pi_\bigstar^G(HM)$ is only computable when all Sylow subgroups of $G$ has prime order.
\end{remark}

\subsection{Computibility of $MU_G$}

We will discuss the computibility of $MU_G$ in this section. According to the special structure of Thom spectra, $MU_G$ is more computable than $H\underline{\Z}$ in some sense. We will prove the following theorem:

\begin{theorem} \label{mucomp}
The homotopy of $MU_G$ is computable if all Sylow subgroups of $G$ are abelian.
\end{theorem}

We just need to check the list of input in \autoref{machineinput}:
\medskip

The homotopy of $MU_H$ for abelian $H$ is computed by \cite{AK}, which uses the Tate diagram repeatedly with families corresponding to all sequences of subgroups.

The conjugacy maps become isomorphisms of diagrams. For any $H_1\subset H_2$, $res_{H_1}^{H_2}$ can be obtained by the fact that the diagram computing $MU_{H_1}$ is isomorphism to a sub-diagram (after taking underlying $H_1$-spectra) of the diagram computing $MU_{H_2}$.
\medskip

The study of geometric fixed point spectrum $\Phi^KMU_G$ traces back to \cite{Dieck2}. That paper only discusses the case when $G$ is abelian, but some of the results still work for the non-abelian case. We refer to \cite{Gre} for a more modern discussion:

\begin{theorem} \label{mugeo}
Let $J(K)$ be the set of all complex virtual $K$-representations with dimension $0$, and $t(K)$ be the number of isomorphism classes of non-trivial irreducible complex representations of $K$. Then as a $G/K$-spectrum, we have
$$\Phi^KMU_G\simeq\bigvee_{j\in J(K)}\Sigma^{|j^K|}MU_{G/K}\wedge (BU_{G/K})_+^{\wedge t(K)}$$
\end{theorem}

When $G/K$ is abelian, the homotopy of each summand is computed in \cite{CGK}:

\begin{theorem} \label{mubu}
Let $H$ be an abelian group. Then
$$\pi_*^H(MU_H\wedge BU_H^{\wedge n})=(\pi_*^HMU_H)[\gamma_1,\gamma_2,...]^{\otimes n}.$$
\end{theorem}

The conjugacy maps fix the polynomial generators $\gamma_i$. For any $H_1\subset H_2$ such that $(H_1)_p=(H_2)_p$ (written as $H_p$), the $H_p$-geometric fixed point spectra of $MU_{H_1}$ and $MU_{H_2}$ have the same number of $BU$ factors in the expression in \autoref{mugeo}. The restriction map $res_{H_1}^{H_2}$ can be obtained by the restriction map
$$\pi_*^{H_2/H_p}MU_{H_2/H_p}\rightarrow\pi_*^{H_1/H_p}MU_{H_1/H_p}$$
and identity maps on the polynomial generators.
\medskip

Therefore, \autoref{mucomp} is proved since all inputs in \autoref{machineinput} are computable.
\bigskip

\autoref{machineinput} also implies the following approach on the evenness conjecture of $MU_G$:

\begin{theorem} \label{mueven}
If all Sylow subgroups of $G$ are abelian, then the homotopy of $MU_G$ concentrates in even degrees.

In general, if for any $p$-subquotient $P$ of $G$ and $n\geq 0$, the homotopy of $MU_P\wedge(BU_P)_+^{\wedge n}$ is even, then the homotopy of $MU_G$ is even.
\end{theorem}

\subsection{Some Examples}

In this subsection, we will apply \autoref{machineinput} to compute general $G$-spectra with $G=D_{2p}$ and $A_5$. These results will be used in later computations.
\bigskip

For $G=D_{2p}$, the subgroups of $G$ (up to conjugacy) are listed as $e,C_2,C_p,G$.

\begin{theorem} \label{d2pgeneral}
Let $G=D_{2p}$ and $V$ be a $G$-virtual representation. For any $G$-spectrum $X$, $\pi_V^G(X)$ is isomorphic to:

\textbf{(i)} The direct sum of $\pi_{V^G}(\Phi^GX)$, $\pi_{(V|_{C_2})^{C_2}}^{C_2}(\Phi^{C_2}X)$ and $\pi_{V|_{C_p}}^{C_p}(X)^{G/C_p}$ when $2$ is inverted;

\textbf{(ii)} The direct sum of $\pi_{V^{C_p}}^{G/C_p}(\Phi^{C_p}X)$ and the pullback of
$$\pi_{V|_{C_2}}^{C_2}(X)\xrightarrow{res_e^{C_2}}\pi_{|V|}(X)\hookleftarrow\pi_{|V|}(X)^G$$
when $p$ is inverted.
\end{theorem}

\paragraph{Proof:} When $2$ is inverted, according to \autoref{machineinput}, $\pi_V^G(X)$ is the limit of the diagram with objects $\pi_{|V|}(X)$ (with $G$-action), $\pi_{V|_{C_p}}^{C_p}(X)$ (with $G/C_p$-action), $\pi_{V^G}(\Phi^GX)$, $\pi_{(V|_{C_2})^{C_2}}^{C_2}(\Phi^{C_2}X)$, and the map 
$$res_e^{C_p}:\pi_{V|_{C_p}}^{C_p}(X)\rightarrow\pi_{|V|}(X).$$
Since the image of $res_e^{C_p}$ is fixed by $C_p$, the image of $\pi_{V|_{C_p}}^{C_p}(X)^{G/C_p}$ under $res_e^{C_p}$ is fixed by $G$. Thus all information related to $\pi_{|V|}(X)$ and $\pi_{V|_{C_p}}^{C_p}(X)$ is encoded in $\pi_{V|_{C_p}}^{C_p}(X)^{G/C_p}$.

When $p$ is inverted, part (ii) is a direct consequence of \autoref{machineinput}. $\Box$
\bigskip

For $G=A_5$, the subgroups of $G$ (up to conjugacy) are listed as $e,C_2,C_3,K_4,C_5,D_6,D_{10},A_4,G$. Here $K_4=C_2\times C_2$ is the Sylow $2$-subgroup of $G$.
\smallskip

The case of $G=A_5$ is a little special. Notice that $G$ has no proper normal subgroups and we have $G=G_p$ for all prime $p$. When we apply the splitting method described at the begining of section 4, the family $\mathscr{F}$ will always be the family of all proper subgroups no matter which primes we invert. Thus we can apply \autoref{moregeneralefsplit} and \autoref{geosplitpart} directly without any localizations:

\begin{theorem} \label{a5generalpre}
Let $G=A_5$ and $V$ be a $G$-virtual representation. For any $G$-spectrum $X$, $\pi_V^G(X)$ is isomorphic to the direct sum of $\pi_{V^G}(\Phi^GX)$ and the limit of the following diagram:
$$\xymatrix{ & \pi_{V|_{A_4}}^{A_4}(X) \ar[d]^{res_e^{A_4}} & \\
\pi_{V|_{D_6}}^{A_4}(X) \ar[r]^{res_e^{D_6}} & \pi_{|V|}(X) & \pi_{V|_{D_{10}}}^{D_{10}}(X) \ar[l]^{res_e^{D_{10}}}
}$$
\end{theorem}

The diagram above comes from \autoref{moregeneralefsplit} and is simplified by similar arguments as in the proof of \autoref{d2pgeneral} part (i).
\medskip

On the other hand, we can also apply \autoref{machineinput} after inverting proper primes to get a complete decomposition. For simplicity, we will use the same notation to denote the underlying representations of any $G$-virtual representation.

\begin{theorem} \label{a5general}
Let $G=A_5$ and $V$ be a $G$-virtual representation. For any $G$-spectrum $X$, $\pi_V^G(X)$ is isomorphic to:
\smallskip

\textbf{(i)} The direct sum of $\pi_{V^G}(\Phi^GX)$, $\pi_{V^{D_6}}(\Phi^{D_6}X)$, $\pi_{V^{A_4}}(\Phi^{A_4}X)$, $\pi_{V^{D_{10}}}(\Phi^{D_{10}}X)$, $\pi_{V^{K_4}}(\Phi^{K_4}X)^{A_4/K_4}$, $\pi_{V^{C_3}}(\Phi^{C_3}X)^{D_6/C_3}$, $\pi_{V^{C_2}}(\Phi^{C_2}X)^{K_4/C_2}$ and the limit of
$$\pi_V^{C_5}(X)^{D_{10}/C_5}\xrightarrow{res_e^{C_5}}\pi_V(X)\hookleftarrow\pi_V(X)^G$$
when $2,3$ are inverted.
\smallskip

\textbf{(ii)} The direct sum of $\pi_{V^G}(\Phi^GX)$, $\pi_{V^{D_6}}(\Phi^{D_6}X)$, $\pi_{V^{K_4}}^{A_4/K_4}(\Phi^{K_4}X)$, $\pi_{V^{D_{10}}}(\Phi^{D_{10}}X)$, $\pi_{V^{K_4}}(\Phi^{K_4}X)^{A_4/K_4}$, $\pi_{V^{C_5}}(\Phi^{C_5}X)^{D_{10}/C_5}$, $\pi_{V^{C_2}}(\Phi^{C_2}X)^{K_4/C_2}$, and the limit of
$$\pi_V^{C_3}(X)^{D_6/C_3}\xrightarrow{res_e^{C_3}}\pi_V(X)\hookleftarrow\pi_V(X)^G$$
when $2,5$ are inverted.

\textbf{(iii)} The direct sum of $\pi_{V^G}(\Phi^GX)$, $\pi_{V^{C_3}}^{D_6/C_3}(\Phi^{C_3}X)$, $\pi_{V^{A_4}}(\Phi^{A_4}X)$, $\pi_{V^{C_5}}^{D_{10}/C_5}(\Phi^{C_5}X)$, $\pi_{V^{C_3}}(\Phi^{C_3}X)^{D_6/C_3}$, $\pi_{V^{C_5}}(\Phi^{C_5}X)^{D_{10}/C_5}$, and the limit of
$$\pi_{V|_{K_4}}^{K_4}(X)^{A_4/K_4}\xrightarrow{res_e^{K_4}}\pi_{|V|}(X)\hookleftarrow\pi_{|V|}(X)^G$$
when $3,5$ are inverted.
\end{theorem}

\subsection{An algebraic point of view}

We want to mention a special case of \autoref{moregeneralefsplit} and its relation with a recent paper by Angeltveit \cite{Ang}:

\begin{theorem} \label{generalsplit}
Let $\mathscr{F}=\mathscr{F}_p$ be the family of all $p$-subgroups of $G$. Let $P\subset G$ be a Sylow $p$-subgroup. Assume that all prime factors of $|G|$ except $p$ are inverted. Then for any virtual $G$-representation $V$, The composition
$$\pi_V^G(E\mathscr{F}_+\wedge X)\rightarrow\pi_V^G(X)\xrightarrow{res_P^G}\pi_{V|_P}^P(X)$$
is injective. Its image consists of all elements $x\in\pi_{V|_P}^P(X)$ such that
$$res^P_H(x)=c_g(res^P_{g^{-1}Hg}(x))$$
for any $H\subset P$ and $g\in G$ such that $g^{-1}Hg\subset P$. Here $c_g$ is the conjugation map induced by $g\in G$.
\end{theorem}

\paragraph{Proof:} It suffices to check the criterion in \autoref{splitcrit}: For any $H\in\mathscr{F}$, $H\lhd J$ such that $|J/H|$ is a power of $p$, since $H$ itself is a $p$-group, $J$ is also a $p$-group. Thus $J$ cannot be outside $\mathscr{F}$ according to the construction of $\mathscr{F}$. $\Box$
\bigskip

When we choose $X=HM$ for some cohomological $M$, the same decomposition problem is studied in a more algebraic way:

\begin{definition}
A Mackey functor $M$ is called cohomological if $tr_K^H\circ res_K^H$ is the multiplication by $|H/K|$ for any $K\subset H\subset G$.
\end{definition}

\begin{theorem} \label{angeltveit}
\cite[Theorem 3.1]{Ang} Let $P_1,...,P_n$ be Sylow subgroups of $G$, with one for each prime factor. There is an isomorphism
$$\oplus_i tr_{P_i}^G:\bigoplus_iM(G/P_i)/\sim\rightarrow M(G/G),$$
where $\sim$ is generated by
$$tr_H^{P_i}(y)\sim tr_{g^{-1}Hg}^{P_j}(c_{g^{-1}}(y))$$
for any $y\in M(G/H)$ and $g\in G$ such that $H\subset P_i$, $g^{-1}Hg\subset P_j$.
\end{theorem}

This theorem can be used to compute equivariant homology with coefficients in a cohomological Mackey functor (like $\underline{\Z}$) since

\begin{proposition}
\cite[Proposition 2.5]{Ang} Let $X$ be a $G$-CW spectrum. If $M$ is a cohomological Mackey functor, then $\underline{\pi}_\bigstar^G(X\wedge HM)$ is also cohomological.
\end{proposition}
\bigskip

Notice that, in \autoref{angeltveit}, the equivalence relation is generated by the following two types:
\smallskip

\textbf{(1)} The equivalence relation on one single summand $M(G/P_i)$:
$$tr_H^{P_i}(y)\sim tr_{g^{-1}Hg}^{P_i}(c_{g^{-1}}(y))$$
for $y\in M(G/H)$ such that $H,g^{-1}Hg\subset P_i$;

\textbf{(2)} For $i\neq j$, we have
$$tr_e^{P_i}(y)\sim tr_e^{P_j}(y)$$
for any $y\in M(G/e)$. This is equivalent to $tr_{P_i}^G\circ tr_e^{P_i}=tr_{P_j}^G\circ tr_e^{P_j}$ as part of the Mackey functor structure.
\bigskip

Since $tr_{P_i}^G\circ res_G^{P_i}$ and $tr_e^{P_i}\circ res_e^{P_i}$ are multiplications by $|G/P_i|$ and $|P_i|$, it's not hard to check that \autoref{angeltveit} is implied by (in fact equivalent to) its localized version:

\begin{proposition} \label{cohmackey}
Assume that all prime factors of $|G|$ except $p$ are inverted. Let $P$ be a Sylow $p$-subgroup of $G$. For any cohomological Mackey functor $M$, we have an isomorphism
$$tr_P^G: M(G/P)/\sim\rightarrow M(G/G),$$
where $\sim$ is generated by
$$tr_H^P(y)\sim tr_{g^{-1}Hg}^P(c_{g^{-1}}(y))$$
for any $y\in M(G/H)$ and $g\in G$ such that $H,g^{-1}Hg\subset P$.
\end{proposition}
\bigskip

This proposition is related to \autoref{generalsplit} by the following two lemmas, which we will prove later:

\begin{lemma} \label{samerestr}
With the same assumptions as in \autoref{cohmackey}, let $R$ be the subgroup of $M(G/P)$ which consists of all elements $x$ such that
$$res_H^P(x)=c_g(res_{g^{-1}Hg}^P(x))$$
for any $H\subset P$ and $g\in G$ such that $g^{-1}Hg\subset P$. Then the composition
$$R\hookrightarrow M(G/P)\xrightarrow{tr_P^G}M(G/G)$$
is an isomorphism.
\end{lemma}

\begin{lemma} \label{cohtilde}
Assume that all prime factors of $|G|$ except $p$ are inverted. Let $\mathscr{F}=\mathscr{F}_p$. Then for any cohomological Mackey functor $M$, we have
$$\widetilde{E\mathscr{F}}\wedge HM\simeq *.$$
\end{lemma}
\medskip

\begin{remark}
When we are computing equivariant homology with coefficients in a cohomological Mackey functor, \autoref{cohtilde} shows that smashing with $E\mathscr{F}_+$ makes no changes. Then $\autoref{samerestr}$ shows that the methods in \cite{Ang} and our paper provide the same result.

In fact, our method also works for some non-cohomological $M$. We only need to guarantee
$$\widetilde{E\mathscr{F}}\wedge HM\simeq *.$$
This condition holds if and only if $M$ is an $M_\mathscr{F}$-module.
\end{remark}

\begin{remark}
From an algebraic point of view, both our method and the method in \cite{Ang} share the same idea: Expressing the top level $M(G/G)$ of the Mackey functor $M$ by the lower levels. 

The only difference is whether we use transfer maps or restriction maps, each of which has its own advantage: The expression with transfer maps does not require any localizations, while the expression with restriction maps preserves the multiplicative structure.
\end{remark}
\medskip

\paragraph{Proof of \autoref{samerestr}:} Since $M$ is cohomological, $tr_P^G\circ res_P^G$ is multiplication by $|G/P|$, hence is an isomorphism from $M(G/G)$ to itself. Thus $res_P^G$ is injective and $tr_P^G$ is surjective.

Fix $z\in M(G/G)$. For any $H\subset P$ and $g\in G$ such that $g^{-1}Hg\subset P$, we have
$$c_g(res_{g^{-1}Hg}^P(res_P^G(z)))=c_g(res_{g^{-1}Hg}^G(z))=res_H^G(c_g(z))=res_H^G(z)=res_H^P(res_P^G(z)).$$
Since $H,g$ are arbitrary, $res_P^G(z)$ is contained in $R$. Thus $res_P^G$ can be viewed as an injective map from $M(G/G)$ to $R$, which is the right inverse of the composition 
$$R\hookrightarrow M(G/P)\xrightarrow{tr_P^G}M(G/G).$$

On the other hand, for any $x\in R$, choose $H=P$ and $g\in N_GP$. We get $x=c_g(x)$. Thus
$$res_P^G(tr_P^G(x))=\sum_{gP\in N_GP/P}c_g(x)=|W_GP|x.$$
Since $|W_GP|$ is inverted, $res_P^G\circ tr_P^G$ becomes a self-isomorphism of $R$. Therefore, the composition
$$R\hookrightarrow M(G/P)\xrightarrow{tr_P^G}M(G/G)$$
is an isomorphism. $\Box$
\bigskip

\paragraph{Proof of \autoref{cohtilde}:} With the same argument as in the proof of \autoref{trivialtilde}, it suffices to prove $M=M\text{ }\Box\text{ }M_\mathscr{F}$ for any cohomological Mackey functor $M$.
\medskip

For any $H\subset G$, let $P$ be a Sylow $p$-subgroup of $H$. Then $tr_P^H\circ res_P^H$ is multiplication by $|H/P|$, hence is an isomorphism. Thus $res_P^H$ is an injection.

Notice that for any $p$-subgroup $K$ of $G$, according to the definition of $M_\mathscr{F}$, we have $M_\mathscr{F}(G/K)=A(K)$. Thus $M\text{ }\Box\text{ }M_\mathscr{F}(G/P)\cong M(G/P)$.

Consider the following commutative diagram.
$$\xymatrix{
M(G/H) \ar[r] \ar[d]^{res_P^H} & M\text{ }\Box\text{ }M_\mathscr{F}(G/H) \ar[d]^{res_P^H}\\
M(G/P) \ar[r] & M\text{ }\Box\text{ }M_\mathscr{F}(G/P)}$$

We have shown that the left vertical map and the bottom horizonal map are isomorphisms. Moreover, the top horizontal map is surjective, since
$$M=(M\text{ }\Box\text{ }M_\mathscr{F})\oplus(M\text{ }\Box\text{ }N_\mathscr{F}).$$
Thus all maps in the diagram above are isomorphisms. Since $H\subset G$ is arbitrary, we have $M=M\text{ }\Box\text{ }M_\mathscr{F}$. $\Box$

\section{Computation of $\pi_\bigstar^{D_{2p}}(H\underline{\Z})$}

Starting from this section, we will apply our splitting method to compute $\pi_\bigstar^G(H\underline{\Z})$ for several choices of $G$ as $RO(G)$-graded rings.
\medskip

We will use the same way to define generators as in \cite{HHR} Definition 3.4:

\begin{definition} \label{generator}

For any non-virtual $G$-representation $V$ with $V^G=0$, let $a_V\in\pi_{-V}(S^0)$ be the map $S^0\rightarrow S^V$ embedding $S^0$ to $0$ and $\infty$. We also use $a_V$ to denote the Hurewicz image of $a_V\in\pi_{-V}(S^0)$ in $\pi_{-V}(H\underline{\Z})$.

For any non-virtual orientable representation $V$ of dimension $n$, let $u_V$ be the generator of $\pi_{n-V}(H\underline{\Z})=H_n^G(S^V;\underline{\Z})$ which restricts to the choice of orientation in
$$\underline{H}_n^G(S^V;\underline{\Z})(G/e)\cong H_n(S^n;\Z).$$

\end{definition}

Some important relations on these generators are given below:

\begin{proposition} \label{genrelation}
\textbf{(a)} For any $V_1,V_2$,
$$a_{V_1+V_2}=a_{V_1}a_{V_2},\text{ }u_{V_1+V_2}=u_{V_1}u_{V_2}.$$

\textbf{(b)} Let $G_V$ be the isotropy subgroup of $V$. Then $|G/G_V|a_V=0$.

\textbf{(c)} For $V,W$ both oriented with dimension $2$, with $G_V\subset G_W$, we have
$$a_Wu_V=|G_W/G_V|a_Vu_W.$$
\end{proposition}
\medskip

\paragraph{Convention:} For any $H\subset G$ and $H$-representations $V_1,V_2$, we will identify them in our computation if
$$S^{V_1}\wedge HA_H\simeq S^{V_2}\wedge HA_H.$$
We may not have $V_1\simeq V_2$. But in the computation of homology, it's not necessary to distinguish these two representations.
\bigskip

For the rest of this section, let $G=D_{2p}$.
\medskip

Up to conjugacy, the nontrivial subgroups of $G$ consist of $C_2$, $C_p$ and $G$ itself. According to the splitting method in section 4.1, the computation of $D_{2p}$-homotopy can be decomposed into computations of $C_2$ and $C_p$-homotopy.
\medskip

We first list all irreducible representations for these subgroups except the trivial representation:

For $C_2$, we have the 1-dimensional sign representation $\sigma$.

For $C_p$, we have $p-1$ different 2-dimensional rotation representations, which are identified according to the convention above. We denote them by $\lambda$.

For $D_{2p}$, we have the sign representation $\sigma$ for $D_{2p}/C_p$, and the 2-dimensional dihedral representation $\gamma$, where all elements with order $p$ act as rotations and all elements with order $2$ act as reflections.

\subsection{Decomposition into $C_2$ and $C_p$}

We apply the splitting method and invert $2$ and $p$ separately. Notice that $\pi_*(H\underline{\Z})=\Z$ which is $G$-fixed. According to \autoref{d2pgeneral} and \autoref{trivialtilde}, we get:

\textbf{(i)} $\pi_V^{D_{2p}}(H\underline{\Z})[1/p]\cong\pi^{C_2}_{V|_{C_2}}(H\underline{\Z})[1/p]$.

\textbf{(ii)} $\pi_V^{D_{2p}}(H\underline{\Z})[1/2]\cong\pi^{C_p}_{V|_{C_p}}(H\underline{\Z})^{D_{2p}/C_p}[1/2]$.
\bigskip

It suffices to compute the action of $D_{2p}/C_p$ on $\pi_{V|_{C_p}}^{C_p}(H\underline{\Z})$. Express the generators of $D_{2p}$ by $\zeta$ and $\tau$ such that $\zeta^p=\tau^2=1$, $\zeta\tau=\tau\zeta^{-1}$.
\bigskip

First we consider the case that $V$ has no copies of $\sigma$:
\medskip

\begin{proposition} \label{tauactionp}
Let $V=k+n\gamma$ for some integers $k,n$. Then the action of $\tau$ on $[S^V,\Sigma^t H\underline{\Z}]^{C_p}$ is multiplication by $-1$ if $|k-t|\equiv 2$ or $3$ (mod $4$). Otherwise, the action is multiplication by $1$.
\end{proposition}

\paragraph{Proof:} Notice that $[S^V,\Sigma^t H\underline{\Z}]^{C_p}$ can be viewed as the equivariant cohomology or homology of $S^{|n|\gamma}$ with coefficients in $\underline{\Z}$. The computation can be done in a cellular way. We first give $S^{|n|\gamma}$ a $C_p$-CW structure.
\smallskip

We view $S^\gamma$ as $\gamma$ compactified at $\infty$. The action of $\zeta$ is the counter-clockwise rotation by $2\pi/p$ and the action of $\tau$ is the reflection by the $x$-axis. The standard $CW$ structure of $S^\gamma$ can be described as follows:
\medskip

As a based non-equivariant space, $S^\gamma\cong S^2$ can be constructed by one single $2$-cell and the base point at the origin. Denote that $2$-cell by $a$.
\smallskip

As a based $C_p$-space, $S^\gamma\cong S^\lambda$ can be constructed by:

One $C_p$-free $2$-cell denoted by $\{b_1,b_2,...,b_p\}$, such that each $b_i$ has image
$$\left\{(r\cos\theta,r\sin\theta):\text{ }0\leq r\leq\infty, \frac{2\pi(i-1)}{p}\leq\theta\leq\frac{2\pi i}{p}\right\};$$

One $C_p$-free $1$-cell denoted by $\{c_1,c_2,...,c_p\}$, such that each $c_i$ has image
$$\left\{\left(r\cos\frac{2\pi i}{p},r\sin\frac{2\pi i}{p}\right):\text{ }0\leq r\leq\infty\right\};$$

One $C_p$-fixed $0$-cell $d$ at $\infty$ and the base point at $0$.

Here the subscripts of $b_i$ and $c_i$ are defined modulo $p$.
\medskip

For both $CW$ structures, the action of $\tau$ on $S^\gamma$ can be made into a cellular map:
\smallskip

The action of $\tau$ on the whole $S^\gamma$ is a reflection. Thus $\tau a=-a$.

The action of $\tau$ gives a permutation among the images of $b_1,b_2,...,b_p$ and reverses the orientations of these $2$-cells. Thus $\tau b_i=-b_{p+1-i}$.

The action of $\tau$ gives a permutation among the images of $c_1,c_2,...,c_p$ and keeps the orientations of these $1$-cells. Thus $\tau c_i=c_{p+1-i}$.

Finally, $\tau d=d$ since the point $\infty$ is fixed.
\medskip

Since $S^{|n|\gamma}=S^\gamma\wedge S^\gamma\wedge...\wedge S^\gamma$, we can construct a CW structure as follows:

For $i\in\Z$ and $0\leq j\leq |n|-1$, let
$$e_{2j+2,i}:=(\zeta^ia,...,\zeta^ia,b_i,d,...,d)$$
$$e_{2j+1,i}:=(\zeta^ia,...,\zeta^ia,c_i,d,...,d)$$
with $j$ copies of $\zeta^ia$ and $|n|-1-j$ copies of $d$.

We also have one $0$-cell $e_0:=(d,...,d)$ and the base point.

Notice that $\zeta b_i=b_{i+1}$, $\zeta c_i=c_{i+1}$. This construction makes $S^{|n|\gamma}$ into a $C_p$-CW complex with one $C_p$-cell (except the base point) in each degree between $0$ and $2|n|$. The action of $\tau$ on $S^{|n|\gamma}$ is induced by the action on each copy of $S^\gamma$, hence is also made cellular:
$$\tau e_{2j+2,i}=(-1)^{j+1}e_{2j+2,p+1-i},\text{ }\tau e_{2j+1,i}=(-1)^je_{2j+1,p+1-i},\text{ }\tau e_0=e_0.$$
\medskip

Recall that $H_*^{C_p}(S^{|n|\gamma};\underline{\Z})$ and $H^*_{C_p}(S^{|n|\gamma};\underline{\Z})$ can be computed by the following chain and cochain complexes:
$$C_*^{C_p}(S^{|n|\gamma};\underline{\Z}):=\underline{C}_*(S^{|n|\gamma})\otimes_{\mathcal{O}_{C_p}}\underline{\Z};$$
$$C^*_{C_p}(S^{|n|\gamma};\underline{\Z}):=Hom_{\mathcal{O}_{C_p}}(\underline{C}_*(S^{|n|\gamma});\underline{\Z}),$$
where we view $\underline{\Z}$ as a covariant coefficient system in the first equation and a contravariant one in the second equation.

The actions of $\tau$ on the homology and cohomology are obtained by applying $\tau$ on both $\underline{C}_*(S^{|n|\gamma})$ and $\underline{\Z}$.
\smallskip

Since $S^{|n|\gamma}$ has only $C_p$-free cells in positive degrees and only $C_p$-fixed cells in degree $0$, and $\underline{\Z}$ is the constant coefficient system, we have
$$C_*^{C_p}(S^{|n|\gamma};\underline{\Z}):=\underline{C}_*(S^{|n|\gamma})\otimes_{\mathcal{O}_{C_p}}\underline{\Z}\cong C_*(S^{|n|\gamma})/C_p;$$
$$C^*_{C_p}(S^{|n|\gamma};\underline{\Z}):=Hom_{\mathcal{O}_{C_p}}(\underline{C}_*(S^{|n|\gamma});\underline{\Z})\cong Hom_\Z(C_*(S^{|n|\gamma})/C_p;\Z).$$

Let $e_{2j}$ and $e_{2j-1}$ be the orbits of cells $e_{2j,i}$ and $e_{2j-1,i}$. The induced $\tau$-action on $C_*(S^{|n|\gamma})/C_p$ is expressed as
$$\tau e_{2j}=(-1)^je_{2j},\text{ }\tau e_{2j-1}=(-1)^{j-1}e_{2j-1}.$$

Therefore, the $\tau$-action on the homology and cohomology is multiplication by $-1$ when the degree is congruent to $2$ or $3$ mod $4$. Otherwise the $\tau$-action is multiplication by $1$. $\Box$
\bigskip

Now we consider the case that $V$ also contains copies of $\sigma$. For any $D_{2p}$-spectra $X,Y$, let $\tau$ act on $[X,Y]^{C_p}$ by conjugation. Since $S^\sigma\cong S^1$ as $C_p$-spaces, we have
$$[X,Y]^{C_p}\cong [X,S^{m(\sigma-1)}\wedge Y]^{C_p}.$$
Since $\tau$ acts on $S^\sigma$ as a reflection, the composition
$$[X,Y]^{C_p}\cong [X,S^{m(\sigma-1)}\wedge Y]^{C_p}\xrightarrow{\tau}[X,S^{m(\sigma-1)}\wedge Y]^{C_p}\cong [X,Y]^{C_p}$$
agrees with the $\tau$-action on $[X,Y]^{C_p}$ multiplied by $(-1)^m$.
\medskip

Together with \autoref{tauactionp}, we get

\begin{proposition} \label{tauaction}
Let $V=k+m\sigma+n\gamma$. The action of $\tau$ on $[S^V,\Sigma^tH\underline{\Z}]^{C_p}$ is multiplication by $-1$ if $[|k+m-t|/2]+m$ is odd. Otherwise, the action is multiplication by $1$.
\end{proposition}
\bigskip

In conclusion, we have

\begin{proposition} \label{d2pdecom}
For any integers $k,m,n$,
$$[S^{k+m\sigma+n\gamma},H\underline{\Z}]^G[1/p]\cong[S^{(k+n)+(m+n)\sigma},H\underline{\Z}]^{C_2}[1/p],$$
$$[S^{k+m\sigma+n\gamma},H\underline{\Z}]^G[1/2]\cong
\begin{cases}
0,\text{ if }[|k+m|/2]+m\text{ is odd,}\\
[S^{(k+m)+n\lambda},H\underline{\Z}]^{C_p}[1/2],\text{ if }[|k+m|/2]+m\text{ is even.}
\end{cases}$$
\end{proposition}

\subsection{A complete expression}

Now we combine \autoref{d2pdecom} with the computations of $\pi_\bigstar^{C_2}(H\underline{\Z})$ and $\pi_\bigstar^{C_p}(H\underline{\Z})$ to get an explicit expression for $\pi_\bigstar^{D_{2p}}(H\underline{\Z})$. We refer to \cite{Zen}, which provides a modern method to compute $C_2$ and $C_p$-homotopy of $H\underline{\Z}$ and gives the following expression:

\begin{theorem} \label{primecase}
The $C_2$ and $C_p$-homology of a point are given below:
$$\pi^{C_2}_\bigstar(H\underline{\Z})=\Z[u_{2\sigma},a_\sigma]/(2a_\sigma)\oplus \left(\bigoplus_{i>0}2\Z\langle u_{2\sigma}^{-i}\rangle\right)\oplus\left(\bigoplus_{j,k>0}\Z/2\langle\Sigma^{-1}u_{2\sigma}^{-j}a_\sigma^{-k}\rangle\right).$$
$$\pi^{C_p}_\bigstar(H\underline{\Z})=\Z[u_\lambda,a_\lambda]/(pa_\lambda)\oplus\left(\bigoplus_{i>0}p\Z\langle u_\lambda^{-i}\rangle\right)\oplus\left(\bigoplus_{j,k>0}\Z/p\langle\Sigma^{-1}u_\lambda^{-j}a_\lambda^{-k}\rangle\right).$$
Here the generators $u_{2\sigma}, a_\sigma, u_\lambda, a_\lambda$ are defined in \autoref{generator}.
\end{theorem}
\bigskip

We can obtain generating elements in $\pi_\bigstar^{D_{2p}}(H\underline{\Z})$ by tracing the pre-image of generators in $\pi_\bigstar^{C_2}(H\underline{\Z})$ and $\pi_\bigstar^{C_p}(H\underline{\Z})$.
\bigskip

First we consider the case when $p$ is inverted. According to the first equation in \autoref{d2pdecom}, we have

\textbf{(1)}
$$\pi_0^{C_2}(H\underline{\Z})\cong\pi_{c(1+\sigma-\gamma)}^{D_{2p}}(H\underline{\Z}),\text{ }\forall c\in\Z.$$
We use $u_{\gamma-\sigma}$ to denote the generator of $\pi_{1+\sigma-\gamma}^{D_{2p}}(H\underline{\Z})$. Then $u_{\gamma-\sigma}$ is invertible.
\smallskip

\textbf{(2)}
$$\pi_{k+m\sigma}^{C_2}(H\underline{\Z})\cong\pi_{k+m\sigma+c(1+\sigma-\gamma)}^{D_{2p}}(H\underline{\Z}),\text{ }\forall c\in\Z.$$
Thus the pre-images of $u_{2\sigma},a_\sigma\in\pi_\bigstar^{C_2}(H\underline{\Z})$ are $u_{2\sigma}u_{\gamma-\sigma}^c,a_\sigma u_{\gamma-\sigma}^c\in \pi_\bigstar^{D_{2p}}(H\underline{\Z})$, for any $c\in\Z$.
\medskip

Therefore, $\pi_\bigstar^{D_{2p}}(H\underline{\Z})$ can be obtained by adding the invertible $u_{\gamma-\sigma}$ to $\pi_\bigstar^{C_2}(H\underline{\Z})$:
$$\pi_\bigstar^{D_{2p}}(H\underline{\Z})[\frac{1}{p}]=\Z[\frac{1}{p}][u_{2\sigma},a_\sigma,u_{\gamma-\sigma}^\pm]/(2a_\sigma)$$
$$\oplus\left(\bigoplus_{i>0}2\Z[\frac{1}{p}][u_{\gamma-\sigma}^\pm]\langle u_{2\sigma}^{-i}\rangle\right)\oplus\left(\bigoplus_{j,k>0}\Z/2[u_{\gamma-\sigma}^\pm]\langle\Sigma^{-1}u_{2\sigma}^{-j}a_\sigma^{-k}\rangle\right).$$
\bigskip

The case when $2$ is inverted is similar. According to the second equation in \autoref{d2pdecom}, we have:

\textbf{(1)}
$$\pi_0^{C_p}(H\underline{\Z})\cong\pi_{2c(1-\sigma)}^{D_{2p}}(H\underline{\Z}),\text{ }\forall c\in\Z.$$
Thus $u_{2\sigma}$, as the generator of $\pi_{2-2\sigma}^{D_{2p}}(H\underline{\Z})$, becomes invertible.
\smallskip

\textbf{(2)}
$$\pi_{2-\lambda}^{C_p}(H\underline{\Z})\cong\pi_{1+\sigma-\gamma+2c(1-\sigma)}^{D_{2p}}(H\underline{\Z}),\text{ }\forall c\in\Z.$$
Thus the pre-image of $u_\lambda\in\pi_\bigstar^{C_p}(H\underline{\Z})$ is $u_{\gamma-\sigma}u_{2\sigma}^c$, for any $c\in\Z$.
\smallskip

\textbf{(3)}
$$\pi_{-\lambda}^{C_p}(H\underline{\Z})\cong\pi_{-\gamma+2c(1-\sigma)}^{D_{2p}}(H\underline{\Z}),\text{ }\forall c\in\Z.$$
Thus the pre-image of $a_\lambda\in\pi_\bigstar^{C_p}(H\underline{\Z})$ is $a_\gamma u_{2\sigma}^c$, for any $c\in\Z$.
\medskip

Therefore, we have
$$\pi_\bigstar^{D_{2p}}(H\underline{\Z})[\frac{1}{2}]=\Z[\frac{1}{2}][u_{\gamma-\sigma},a_\gamma,u_{2\sigma}^\pm]/(pa_\gamma)$$
$$\oplus\left(\bigoplus_{i>0}p\Z[\frac{1}{2}][u_{2\sigma}^\pm]\langle u_{\gamma-\sigma}^{-i}\rangle\right)\oplus\left(\bigoplus_{j,k>0}\Z/p[u_{2\sigma}^\pm]\langle\Sigma^{-1}u_{\gamma-\sigma}^{-j}a_\gamma^{-k}\rangle\right).$$
\bigskip

Finally, we just need to glue $\pi_\bigstar^{D_{2p}}(H\underline{\Z})[1/p]$ and $\pi_\bigstar^{D_{2p}}(H\underline{\Z})[1/2]$ together.

\begin{theorem} \label{d2phz}
Recall that we obtain $a_\sigma$, $a_\gamma$ and $u_{2\sigma}$ from \autoref{generator} and choose $u_{\gamma-\sigma}$ to be the generator of
$$\pi_{1+\sigma-\gamma}^{D_{2p}}(H\underline{\Z})\cong\Z.$$
Then we have
$$\pi_\bigstar^{D_{2p}}(H\underline{\Z})=\Z[u_{\gamma-\sigma},u_{2\sigma},a_\sigma,a_\gamma]/(2a_\sigma,pa_\gamma)$$
$$\oplus\left(\bigoplus_{i>0}2\Z[u_{\gamma-\sigma}]\langle u_{2\sigma}^{-i}\rangle\right)\oplus\left(\bigoplus_{i>0}p\Z[u_{2\sigma}]\langle u_{\gamma-\sigma}^{-i}\rangle\right)\oplus\left(\bigoplus_{j,k>0}2p\Z\langle u_{2\sigma}^{-j}u_{\gamma-\sigma}^{-k}\rangle\right)$$
$$\oplus\left(\bigoplus_{j,k>0}\Z/p[u_{\gamma-\sigma}]\langle a_\gamma^j u_{2\sigma}^{-k}\rangle\right)\oplus\left(\bigoplus_{j,k>0}\Z/2[u_{2\sigma}]\langle a_\sigma^j u_{\gamma-\sigma}^{-k}\rangle\right)$$
$$\oplus\left(\bigoplus_{j,k>0}\Z/2[u_{\gamma-\sigma}^\pm]\langle\Sigma^{-1}u_{2\sigma}^{-j}a_\sigma^{-k}\rangle\right)\oplus\left(\bigoplus_{j,k>0}\Z/p[u_{2\sigma}^\pm]\langle\Sigma^{-1}u_{\gamma-\sigma}^{-j}a_\gamma^{-k}\rangle\right).$$
\end{theorem}

\section{Computation of $\pi_\bigstar^{A_5}(H\underline{\Z})$}

In this section, we will use the splitting method to compute $\pi_\bigstar^G(H\underline{\Z})$ for $G=A_5$. First we list all subgroups of $A_5$.
\bigskip

\textbf{Subgroups of $A_5$ up to conjugacy:}
\smallskip

Sylow $2$-subgroup $K_4=\langle (12)(34),(13)(24)\rangle$, which is isomorphic to $C_2\times C_2$;

Sylow $3$-subgroup $C_3=\langle (123)\rangle$;

Sylow $5$-subgroup $C_5=\langle (12345)\rangle$;
\smallskip

Normalizers of Sylow subgroups:

$A_4=N_GK_4=\langle (12)(34),(123)\rangle$;

$D_6=N_GC_3=\langle (123),(12)(45)\rangle$;

$D_{10}=N_GC_5=\langle (12345),(15)(24)\rangle$;
\smallskip

Other subgroups: $C_2=\langle (12)(34)\rangle$ and the trivial subgroup.
\bigskip

We also list all irreducible real representations of $A_5$ and its subgroups. Since there are lots of groups, we express the group and dimension of a representation in the subscript. For example, $V_{H,n}$ is an $n$-dimension $H$-representation. The trivial representations of all groups will be written as $1$.
\bigskip

\textbf{Irreducible representations of $A_5$:}

$V_{A_5,3}^+$: This is obtained by viewing $A_5$ as the group of rotations of the regular icosahedron in the  3 dimensional space.

$V_{A_5,3}^-$: This is the same as $V_{A_5,3}^+$ except applying the conjugation by $(12)$ (inside $S_5$) first.

For the two representations $V_{A_5,3}^+$ and $V_{A_5,3}^-$, we will see that all our computations which work for one of them also work for the other. This further implies that the two representation spheres are $HA_G$-equivalent. We use $V_{A_5,3}$ to denote both representations.
\smallskip

$V_{A_5,4}$: Consider the space of functions from $\{1,2,3,4,5\}$ to $\R$ whose images have sum $0$. The action of $A_5$ is the natural permutation.

$V_{A_5,5}$: Consider the space of functions from antipodal pairs of vertices of the regular icosahedron, such that the sum of all images of the function is $0$. The action of $A_5$ is induced by the rotations of the icosahedron.
\medskip

\textbf{Irreducible representations of $A_4$:}

$V_{A_4,2}$: All elements in $A_4$ with order $3$ become rotations of either $2\pi/3$ or $\pi/3$. All elements with order $2$ act trivially.

$V_{A_4,3}$: Consider the space of functions from $\{1,2,3,4\}$ to $\R$ whose images have sum $0$. The action of $A_4$ is the natural permutation.
\medskip

\textbf{Irreducible representations of $D_{2p}$ ($p=3$ or $5$):}

$V_{D_{2p},1}$: The sign representation of $D_{2p}/C_p$, with $C_p$ acting trivially.

$V_{D_{2p},2}$: Dihedral representations, where all elements with odd orders act as rotations, and all elements with order 2 act as reflections.
\medskip

\textbf{Irreducible representations of $K_4$:}

Three different sign representations $V_{K_4,1,1}, V_{K_4,1,2}, V_{K_4,1,3}$ which correspond to the three proper subgroups of $K_4$.

In addition, we use $V_{K_4,3}$ to denote $V_{K_4,1,1}+V_{K_4,1,2}+V_{K_4,1,3}$ (which is not irreducible).
\medskip

\textbf{Irreducible representations of $C_p$:} ($p=3$ or $5$)

$V_{C_p,2}$: All elements act as rotations.
\medskip

\textbf{Irreducible representations of $C_2$:}

$V_{C_2,1}$: Sign representation.
\bigskip

The representations of these groups are related by restrictions:
\smallskip

The restrictions of $V_{A_5,3}$ in $A_4,D_6,D_{10}$ are $V_{A_4,3}$, $V_{D_6,1}+V_{D_6,2}$ and $V_{D_{10},1}+V_{D_{10},2}$.

The restrictions of $V_{A_5,4}$ in $A_4,D_6,D_{10}$ are $1+V_{A_4,3}$, $1+V_{D_6,1}+V_{D_6,2}$ and $2V_{D_{10},2}$.

The restrictions of $V_{A_5,5}$ in $A_4,D_6,D_{10}$ are $V_{A_4,2}+V_{A_4,3}$, $1+2V_{D_6,2}$ and $1+2V_{D_{10},2}$.

The restrictions of $V_{A_4,2}, V_{A_4,3}$ in $K_4$ are $2$ and $V_{K_4,3}$.

The restrictions of $V_{D_{2p},1}, V_{D_{2p},2}$ in $C_p$ are $1$ and $V_{C_p,2}$ for $p=3,5$.
\bigskip

Applying \autoref{a5general}, eliminating the geometric fixed point spectra with \autoref{trivialtilde}, and simplifying the diagrams by the fact that $\pi_*(H\underline{\Z})=\Z$ is fixed by $G$, we get

\begin{proposition} \label{a5decom}
For $p=2,3,5$, let $P$ be a Sylow $p$-subgroup of $G=A_5$. Assume that all prime factors of $|G|$ except $p$ are inverted. Then
$$\pi_V^G(H\underline{\Z})\cong\pi_{V|_P}^P(H\underline{\Z})^{W_GP}.$$
\end{proposition}
\bigskip

Now we can use the same strategy as in section 5.2 to obtain explicit expressions of localized homotopy:

\begin{theorem} \label{3local}
When $2,5$ are inverted, we have
$$\pi_{\bigstar}^{A_5}(H\underline{\Z})=\Z[u_{V_4-V_3}^\pm,u_{2V_3-V_5}^\pm]\otimes$$
$$\left[\Z[u_{V_3},a_{V_5-1}]/(3a_{V_5-1})\oplus\left(\bigoplus_{i>0}3\Z\langle u_{V_3}^{-i}\rangle\right)\oplus\left(\bigoplus_{j,k>0}\Z/3\langle\Sigma^{-1}u_{V_3}^{-j}a_{V_5-1}^{-k}\rangle\right)\right].$$
Here $V_i$ denotes $V_{A_5,i}$, $i=3,4,5$. For any virtual representation $V$, we have $|a_V|=-V$, $|u_V|=|V|-V$.
\end{theorem}

\begin{theorem} \label{5local}
When $2,3$ are inverted, we have
$$\pi_{\bigstar}^{A_5}(H\underline{\Z})=\Z[u_{V_5-V_4}^\pm,u_{2V_3-V_4}^\pm]\otimes$$
$$\left[\Z[u_{V_3},a_{V_4}]/(5a_{V_4})\oplus\left(\bigoplus_{i>0}5\Z\langle u_{V_3}^{-i}\rangle\right)\oplus\left(\bigoplus_{j,k>0}\Z/5\langle\Sigma^{-1}u_{V_3}^{-j}a_{V_4}^{-k}\rangle\right)\right].$$
\end{theorem}
\bigskip

For $p=2$, we have a similar result:

\begin{theorem} \label{2local}
Let $n_1,n_3,n_4,n_5$ be arbitrary integers. Assume that $3,5$ are inverted. Then
$$\pi_{n_1+n_3V_{A_5,3}+n_4V_{A_5,4}+n_5V_{A_5,5}}^{A_5}(H\underline{\Z})\cong\pi_{(n_1+n_4+2n_5)+(n_3+n_4+n_5)V_{K_4,3}}^{K_4}(H\underline{\Z})^{A_4/K_4}.$$
\end{theorem}

The computation of $\pi_{*+*V_{K_4,3}}^{K_4}(H\underline{\Z})$ is given in Appendix B. Unlike Theorems \ref{3local} and \ref{5local}, it is quite hard to give an explicit expression for this ring. We will provide a less explicit, but still computable expression in Theorems \ref{positivez} and \ref{negativez}.
\bigskip

The unlocalized $\pi_\bigstar^{A_5}(H\underline{\Z})$ can be recovered from Theorems \ref{3local}, \ref{5local}, \ref{2local}, \ref{positivez}, and \ref{negativez}.

\section{Mackey functor valued homology}

In this section, we will discuss the computation of the Mackey functor valued homotopy $\underline{\pi}_\bigstar^G(H\underline{\Z})$. All arguments and theorems in this section still work when $\underline{\Z}$ is replaced by other constant Mackey functors.
\medskip

\paragraph{Notation:} We use $\Pi$ to denote the Mackey functor $\underline{\pi}_V^G(H\underline{\Z})$ when the virtual $G$-representation $V$ is fixed. Then we have
$$\Pi(G/H)=[G/H_+\wedge S^V,H\underline{\Z}]^G\cong [S^V,H\underline{\Z}]^H=\pi_{V|_H}^H(H\underline{\Z}),$$
$$\Pi|_H=\underline{\pi}_{V|_H}^H(H\underline{\Z}).$$
\bigskip

When all prime factors of $|G|$ except $p$ are inverted, since $H\underline{\Z}=H\underline{\Z}\wedge E\mathscr{F}_+$, we can use \autoref{generalsplit} and \autoref{trivialtilde} to compute $[S^V,H\underline{\Z}]^H$ as a subgroup of $[S^V,H\underline{\Z}]^P$
where $P$ is a Sylow $p$-subgroup of $H$. In order to get the Mackey functor structure on $\Pi$, we still need to compute the restrction and transfer maps.
\bigskip

Consider any $H_1\subset H_2\subset G$. Choose Sylow $p$-subgroups $P_1,P_2$ of $H_1,H_2$ such that $P_1\subset P_2$. We will compute $tr_{H_1}^{H_2}$ and $res_{H_1}^{H_2}$ with the Mackey functor structure of $\Pi|_{P_2}$ and conjugation maps.
\smallskip

The restriction map is simple:

\begin{proposition} \label{restrictionmap}
When computing $\Pi(G/H_1)$ and $\Pi(G/H_2)$ as subgroups of $\Pi(G/P_1)$ and $\Pi(G/P_2)$, $res_{H_1}^{H_2}$ is determined by $res_{P_1}^{P_2}$.
\end{proposition}

\paragraph{Proof:} The inclusion $\Pi(G/H_i)\hookrightarrow\Pi(G/P_i)$ is given by $res_{P_i}^{H_i}$, $i=1,2$. The proposition follows from the fact that
$$res_{P_1}^{P_2}\circ res_{P_2}^{H_2}=res_{P_1}^{H_1}\circ res_{H_1}^{H_2}.$$
$\Box$
\bigskip

The transfer map is more complicated:

Write the underlying $P_2$-set of $H_2/H_1$ as the disjoint union of $P_2$-orbits:
$$H_2/H_1=\bigsqcup_{i=1}^t P_2/Q_i.$$
Choose $g_1,g_2,...,g_t\in H_2$ such that each $Q_i$ corresponds to $g_iH_1$ inside $H_2/H_1$. Then we have $g_i^{-1}Q_ig_i\subset H_1$, $i=1,2,...,t$.

\begin{proposition} \label{transfer}
When computing $\Pi(G/H_1)$ and $\Pi(G/H_2)$ as subgroups of $\Pi(G/P_1)$ and $\Pi(G/P_2)$, $tr_{H_1}^{H_2}$ is determined by 
$$\bigoplus_i tr_{Q_i}^{P_2}\circ c_{g_i}\circ res_{g_i^{-1}Q_ig_i}^{P_1}.$$
\end{proposition}

\paragraph{Proof:} Notice that
$$tr_{H_1}^{H_2}: [(H_2/H_1)_+\wedge S^V,H\underline{\Z}]^{H_2}\rightarrow [S^V,H\underline{\Z}]^{H_2}$$
is induced by the $H_2$-stable map
$$S^0=\Sigma^{-W}S^W\rightarrow\Sigma^{-W}(H_2/H_1)_+\wedge S^W=(H_2/H_1)_+.$$
Here $W$ is an $H_2$-representation with an embedding $H_2/H_1\hookrightarrow W$. Taking a tubular neighborhood of $H_2/H_1$ inside $S^W$ and collapsing the complement to a point, we get the map $S^W\rightarrow (H_2/H_1)_+\wedge S^W$ which defines the middle map above.
\medskip

Now consider the commutative diagram
$$\xymatrix{
[(H_2/H_1)_+\wedge S^V,H\underline{\Z}]^{H_2} \ar[d] \ar[r]^-{tr_{H_1}^{H_2}} & 
[S^V,H\underline{\Z}]^{H_2} \ar[d] \\
[(H_2/H_1)_+\wedge S^V,H\underline{\Z}]^{P_2} \ar[r] & 
[S^V,H\underline{\Z}]^{P_2}
}$$
where the vertical maps come from taking the underlying $P_2$-maps, hence are inclusions by \autoref{moregeneralefsplit}. The bottom horizontal map is induced by the underlying $P_2$-stable map of $S^0\rightarrow(H_2/H_1)_+$.

For each $i=1,2,...,t$, compose the $P_2$-stable map $S^0\rightarrow(H_2/H_1)_+$ with the projection $H_2/H_1\rightarrow P_2/Q_i$. The composition can be expressed by
$$S^0=\Sigma^{-W}S^W\rightarrow\Sigma^{-W}(P_2/Q_i)_+\wedge S^W=(P_2/Q_i)_+,$$
which exactly induces $tr_{Q_i}^{P_2}$.

Thus the bottom map
$$[(H_2/H_1)_+\wedge S^V,H\underline{\Z}]^{P_2}\rightarrow [S^V,H\underline{\Z}]^{P_2}$$
is given by
$$[(H_2/H_1)_+\wedge S^V,H\underline{\Z}]^{P_2}\cong\bigoplus_{i=1}^t [(P_2/Q_i)_+\wedge S^V,H\underline{\Z}]^{P_2}\xrightarrow{\oplus tr_{Q_i}^{P_2}}[S^V,H\underline{\Z}]^{P_2}.$$

Moreover, the map
$$[(H_2/H_1)_+\wedge S^V,H\underline{\Z}]^{H_2}\rightarrow [(P_2/Q_i)_+\wedge S^V,H\underline{\Z}]^{P_2}$$
is obtained by taking projection and the underlying $Q_i$-map. Thus we can write it as
$$[(H_2/H_1)_+\wedge S^V,H\underline{\Z}]^{H_2}\cong [S^V,H\underline{\Z}]^{g_iH_1g_i^{-1}}\xrightarrow{res_{g_iP_1g_i^{-1}}^{g_iH_1g_i^{-1}}}[S^V,H\underline{\Z}]^{g_iP_1g_i^{-1}}\xrightarrow{res_{Q_i}^{g_iP_1g_i^{-1}}}[S^V,H\underline{\Z}]^{Q_i}.$$
Commute the conjugation by $g_i$ with the restriction map. The composition above becomes
$$c_{g_i}\circ res_{g_i^{-1}Q_ig_i}^{P_1}\circ res_{P_1}^{H_1}.$$
\smallskip

In conclusion, the left vertical map of the commutative diagram above can be factored as
$$(c_{g_i}\circ res_{g_i^{-1}Q_ig_i}^{P_1})_i\circ res_{P_1}^{H_1},$$
while the bottom horizontal map is given by $\oplus_itr_{Q_i}^{P_2}$. Notice that $res_{P_1}^{H_1}$ is the map we used to compute $\Pi(G/H_1)$ by \autoref{generalsplit}. Thus $tr_{H_1}^{H_2}$ is determined by
$$\bigoplus_i tr_{Q_i}^{P_2}\circ c_{g_i}\circ res_{g_i^{-1}Q_ig_i}^{P_1}.$$
$\Box$
\bigskip

\begin{remark}
In practice, usually we do not need to compute both restriction and transfer maps. It's possible that one of them is induced by the other together with the Mackey functor structure, especially when the Mackey functor is cohomological. Recall that the definition and properties of a cohomological Mackey functor are given in section 4.4.
\end{remark}
\bigskip

Theoretically, whenever our splitting method in section 4 applies, we can always compute the homotopy of $H\underline{\Z}$ as an $RO(G)$-graded Green functor instead of an $RO(G)$-graded ring since \autoref{generalsplit} also preserves the multiplicative structure. But it would be too complicated to finish the complete computation by hand.
\medskip

For the rest of this section, we will give two examples to show how \autoref{restrictionmap} and \autoref{transfer} work in practice.

\subsection{One example}

We will focus on $G=A_5$ since it is the most complicated group for which $\underline{\pi}_\bigstar^G(H\underline{\Z})$ is computable by our splitting methods. The computation for other smaller groups (like $D_{2p}$) will be similar and more simple.
\bigskip

When $H$ is a subgroup of $A_5$, a similar argument as in \autoref{a5decom} gives us

\begin{proposition} \label{a5levels}
Assume that all prime factors of $|A_5|$ except $p$ are inverted. For any $H\subset A_5$, let $P$ be a Sylow $p$-subgroup of $H$. Then we have
$$res_{P}^H:\Pi(A_5/H)\xrightarrow{\sim}\Pi(A_5/P)^{W_HP}.$$
\end{proposition}
\bigskip

We choose a specific virtual representation $V=V_{A_5,3}+V_{A_5,4}-V_{A_5,5}-2$ and compute $\underline{\pi}_V^G(H\underline{\Z})$. We will provide details in the application of \autoref{restrictionmap} and how to use the Mackey functor structure to simplify our computation.
\medskip

We have $V|_{C_3}=0$, $V|_{K_4}=V_{K_4,3}-3$, $V|_{C_5}=V_{C_5,2}-2$. The restrictions of $V$ as Mackey functors over $C_3,K_4,C_5$ are given below:
$$\xymatrix{
\Z \ar@/^/[d]^{1} & \\
\Z \ar@/^/[u]^{3} &
}\text{ }
\xymatrix{
 & \Z \ar@/^/[d]^{2} & \\
 & \Z \ar@/^/[u]^{1} \ar@/^/[d]^{2} & \\
 & \Z \ar@/^/[u]^{1} & 
}\text{ }
\xymatrix{
 & \Z \ar@/^/[d]^{5} \\
 & \Z \ar@/^/[u]^{1}
}$$
The Weyl group actions are all trivial in this case.
\medskip

Using \autoref{a5levels} and \autoref{restrictionmap}, we can obtain the values of all levels in $\Pi$ and the restriction maps when two of the three prime factors $2,3,5$ are inverted.

The first diagram below is the standard expression of an $A_5$-Mackey functor. The other three are the expressions of $\Pi$ (with only restriction maps) when 6, 10, or 15 is inverted respectively.

$$
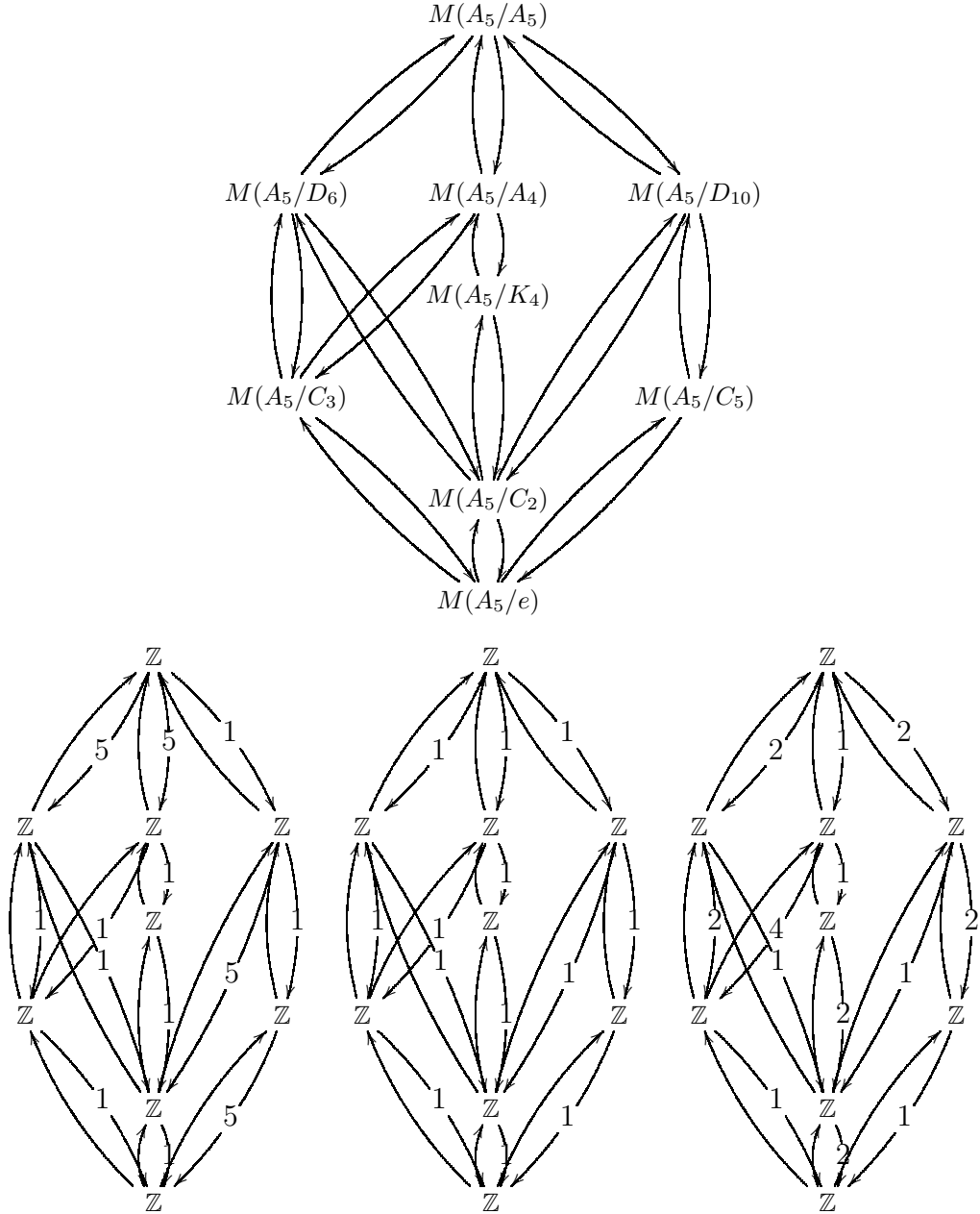
\xymatrix{
  & M(A_5/A_5) \ar@/^/[ddl] \ar@/^/[dd] \ar@/^/[ddr] & \\
  & & \\
  M(A_5/D_6) \ar@/^/[ruu] \ar@/^/[dd] \ar@/^/[dddr] &
  M(A_5/A_4) \ar@/^/[uu] \ar@/^/[d] \ar@/^/[ddl] &
  M(A_5/D_{10}) \ar@/^/[uul] \ar@/^/[dd] \ar@/^/[dddl] \\
  & M(A_5/K_4) \ar@/^/[u] \ar@/^/[dd] & \\
  M(A_5/C_3) \ar@/^/[uu] \ar@/^/[uur] \ar@/^/[ddr] & &
  M(A_5/C_5) \ar@/^/[uu] \ar@/^/[ddl] \\
  & M(A_5/C_2) \ar@/^/[uuul] \ar@/^/[uu] \ar@/^/[uuur] \ar@/^/[d] & \\
  & M(A_5/e) \ar@/^/[uul] \ar@/^/[u] \ar@/^/[uur] & \\
  }$$
$$\xymatrix{
  & \text{ }\text{ }\Z\text{ }\text{ } \ar@/^/[ddl]|{\hole \text{\large{5}} \hole} \ar@/^/[dd]|{\hole \text{\large{5}} \hole} \ar@/^/[ddr]|{\hole \text{\large{1}} \hole} & \\
  & & \\
  \text{ }\text{ }\Z\text{ }\text{ } \ar@/^/[ruu] \ar@/^/[dd]|{\hole \text{\large{1}} \hole} \ar@/^/[dddr]|{\hole \text{\large{1}} \hole} &
  \text{ }\text{ }\Z\text{ }\text{ } \ar@/^/[uu] \ar@/^/[d]|{\hole \text{\large{1}} \hole} \ar@/^/[ddl]|{\hole \text{\large{1}} \hole} &
  \text{ }\text{ }\Z\text{ }\text{ } \ar@/^/[uul] \ar@/^/[dd]|{\hole \text{\large{1}} \hole} \ar@/^/[dddl]|{\hole \text{\large{5}} \hole} \\
  & \text{ }\text{ }\Z\text{ }\text{ } \ar@/^/[u] \ar@/^/[dd]|{\hole \text{\large{1}} \hole} & \\
  \text{ }\text{ }\Z\text{ }\text{ } \ar@/^/[uu] \ar@/^/[uur] \ar@/^/[ddr]|{\hole \text{\large{1}} \hole} & &
  \text{ }\text{ }\Z\text{ }\text{ } \ar@/^/[uu] \ar@/^/[ddl]|{\hole \text{\large{5}} \hole} \\
  & \text{ }\text{ }\Z\text{ }\text{ } \ar@/^/[uuul] \ar@/^/[uu] \ar@/^/[uuur] \ar@/^/[d]|{\hole \text{\large{1}} \hole} & \\
  & \text{ }\text{ }\Z\text{ }\text{ } \ar@/^/[uul] \ar@/^/[u] \ar@/^/[uur] & \\
  }\text{ }
\xymatrix{
  & \text{ }\text{ }\Z\text{ }\text{ } \ar@/^/[ddl]|{\hole \text{\large{1}} \hole} \ar@/^/[dd]|{\hole \text{\large{1}} \hole} \ar@/^/[ddr]|{\hole \text{\large{1}} \hole} & \\
  & & \\
  \text{ }\text{ }\Z\text{ }\text{ } \ar@/^/[ruu] \ar@/^/[dd]|{\hole \text{\large{1}} \hole} \ar@/^/[dddr]|{\hole \text{\large{1}} \hole} &
  \text{ }\text{ }\Z\text{ }\text{ } \ar@/^/[uu] \ar@/^/[d]|{\hole \text{\large{1}} \hole} \ar@/^/[ddl]|{\hole \text{\large{1}} \hole} &
  \text{ }\text{ }\Z\text{ }\text{ } \ar@/^/[uul] \ar@/^/[dd]|{\hole \text{\large{1}} \hole} \ar@/^/[dddl]|{\hole \text{\large{1}} \hole} \\
  & \text{ }\text{ }\Z\text{ }\text{ } \ar@/^/[u] \ar@/^/[dd]|{\hole \text{\large{1}} \hole} & \\
  \text{ }\text{ }\Z\text{ }\text{ } \ar@/^/[uu] \ar@/^/[uur] \ar@/^/[ddr]|{\hole \text{\large{1}} \hole} & &
  \text{ }\text{ }\Z\text{ }\text{ } \ar@/^/[uu] \ar@/^/[ddl]|{\hole \text{\large{1}} \hole} \\
  & \text{ }\text{ }\Z\text{ }\text{ } \ar@/^/[uuul] \ar@/^/[uu] \ar@/^/[uuur] \ar@/^/[d]|{\hole \text{\large{1}} \hole} & \\
  & \text{ }\text{ }\Z\text{ }\text{ } \ar@/^/[uul] \ar@/^/[u] \ar@/^/[uur] & \\
  }\text{ }
\xymatrix{
  & \text{ }\text{ }\Z\text{ }\text{ } \ar@/^/[ddl]|{\hole \text{\large{2}} \hole} \ar@/^/[dd]|{\hole \text{\large{1}} \hole} \ar@/^/[ddr]|{\hole \text{\large{2}} \hole} & \\
  & & \\
  \text{ }\text{ }\Z\text{ }\text{ } \ar@/^/[ruu] \ar@/^/[dd]|{\hole \text{\large{2}} \hole} \ar@/^/[dddr]|{\hole \text{\large{1}} \hole} &
  \text{ }\text{ }\Z\text{ }\text{ } \ar@/^/[uu] \ar@/^/[d]|{\hole \text{\large{1}} \hole} \ar@/^/[ddl]|{\hole \text{\large{4}} \hole} &
  \text{ }\text{ }\Z\text{ }\text{ } \ar@/^/[uul] \ar@/^/[dd]|{\hole \text{\large{2}} \hole} \ar@/^/[dddl]|{\hole \text{\large{1}} \hole} \\
  & \text{ }\text{ }\Z\text{ }\text{ } \ar@/^/[u] \ar@/^/[dd]|{\hole \text{\large{2}} \hole} & \\
  \text{ }\text{ }\Z\text{ }\text{ } \ar@/^/[uu] \ar@/^/[uur] \ar@/^/[ddr]|{\hole \text{\large{1}} \hole} & &
  \text{ }\text{ }\Z\text{ }\text{ } \ar@/^/[uu] \ar@/^/[ddl]|{\hole \text{\large{1}} \hole} \\
  & \text{ }\text{ }\Z\text{ }\text{ } \ar@/^/[uuul] \ar@/^/[uu] \ar@/^/[uuur] \ar@/^/[d]|{\hole \text{\large{2}} \hole} & \\
  & \text{ }\text{ }\Z\text{ }\text{ } \ar@/^/[uul] \ar@/^/[u] \ar@/^/[uur] & \\
  }$$

Now we can recover $\Pi$ (without transfer maps) when no prime factors are inverted, which is given in the left graph below. According to the fact that $\Pi$ is cohomological, the transfer maps are determined. The complete $\Pi$ is given in the right graph below.
$$
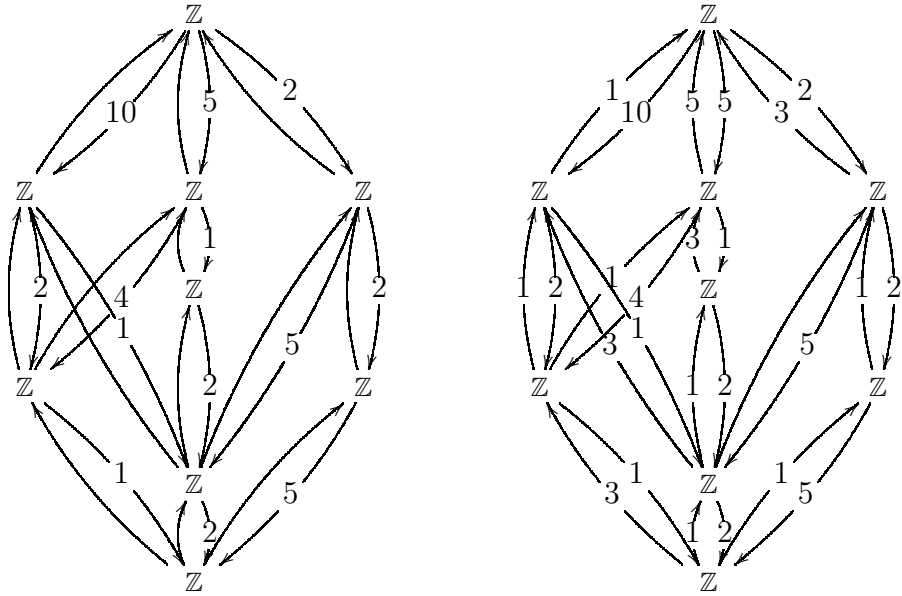
\xymatrix{
  & \text{ }\text{ }\text{ }\Z\text{ }\text{ }\text{ } \ar@/^/[ddl]|{\hole \text{\large{10}} \hole} \ar@/^/[dd]|{\hole \text{\large{5}} \hole} \ar@/^/[ddr]|{\hole \text{\large{2}} \hole} & \\
  & & \\
  \text{ }\text{ }\text{ }\text{ }\Z\text{ }\text{ }\text{ }\text{ } \ar@/^/[ruu] \ar@/^/[dd]|{\hole \text{\large{2}} \hole} \ar@/^/[dddr]|{\hole \text{\large{1}} \hole} &
  \text{ }\text{ }\text{ }\text{ }\Z\text{ }\text{ }\text{ }\text{ } \ar@/^/[uu] \ar@/^/[d]|{\hole \text{\large{1}} \hole} \ar@/^/[ddl]|{\hole \text{\large{4}} \hole} &
  \text{ }\text{ }\text{ }\text{ }\Z\text{ }\text{ }\text{ }\text{ } \ar@/^/[uul] \ar@/^/[dd]|{\hole \text{\large{2}} \hole} \ar@/^/[dddl]|{\hole \text{\large{5}} \hole} \\
  & \text{ }\text{ }\text{ }\text{ }\Z\text{ }\text{ }\text{ }\text{ } \ar@/^/[u] \ar@/^/[dd]|{\hole \text{\large{2}} \hole} & \\
  \text{ }\text{ }\text{ }\text{ }\Z\text{ }\text{ }\text{ }\text{ } \ar@/^/[uu] \ar@/^/[uur] \ar@/^/[ddr]|{\hole \text{\large{1}} \hole}& &
  \text{ }\text{ }\text{ }\text{ }\Z\text{ }\text{ }\text{ }\text{ } \ar@/^/[uu] \ar@/^/[ddl]|{\hole \text{\large{5}} \hole} \\
  & \text{ }\text{ }\text{ }\text{ }\Z\text{ }\text{ }\text{ }\text{ } \ar@/^/[uuul] \ar@/^/[uu] \ar@/^/[uuur] \ar@/^/[d]|{\hole \text{\large{2}} \hole} & \\
  & \text{ }\text{ }\text{ }\text{ }\Z\text{ }\text{ }\text{ }\text{ } \ar@/^/[uul] \ar@/^/[u] \ar@/^/[uur] & \\
  }\text{ }\text{ }\text{ }\text{ }\text{ }\text{ }\text{ }\text{ }
\xymatrix{
  & \text{ }\text{ }\text{ }\Z\text{ }\text{ }\text{ } \ar@/^/[ddl]|{\hole \text{\large{10}} \hole} \ar@/^/[dd]|{\hole \text{\large{5}} \hole} \ar@/^/[ddr]|{\hole \text{\large{2}} \hole} & \\
  & & \\
  \text{ }\text{ }\text{ }\text{ }\Z\text{ }\text{ }\text{ }\text{ } \ar@/^/[ruu]|{\hole \text{\large{1}} \hole} \ar@/^/[dd]|{\hole \text{\large{2}} \hole} \ar@/^/[dddr]|{\hole \text{\large{1}} \hole} &
  \text{ }\text{ }\text{ }\text{ }\Z\text{ }\text{ }\text{ }\text{ } \ar@/^/[uu]|{\hole \text{\large{5}} \hole} \ar@/^/[d]|{\hole \text{\large{1}} \hole} \ar@/^/[ddl]|{\hole \text{\large{4}} \hole} &
  \text{ }\text{ }\text{ }\text{ }\Z\text{ }\text{ }\text{ }\text{ } \ar@/^/[uul]|{\hole \text{\large{3}} \hole} \ar@/^/[dd]|{\hole \text{\large{2}} \hole} \ar@/^/[dddl]|{\hole \text{\large{5}} \hole} \\
  & \text{ }\text{ }\text{ }\text{ }\Z\text{ }\text{ }\text{ }\text{ } \ar@/^/[u]|{\hole \text{\large{3}} \hole} \ar@/^/[dd]|{\hole \text{\large{2}} \hole} & \\
  \text{ }\text{ }\text{ }\text{ }\Z\text{ }\text{ }\text{ }\text{ } \ar@/^/[uu]|{\hole \text{\large{1}} \hole} \ar@/^/[uur]|{\hole \text{\large{1}} \hole} \ar@/^/[ddr]|{\hole \text{\large{1}} \hole}& &
  \text{ }\text{ }\text{ }\text{ }\Z\text{ }\text{ }\text{ }\text{ } \ar@/^/[uu]|{\hole \text{\large{1}} \hole} \ar@/^/[ddl]|{\hole \text{\large{5}} \hole} \\
  & \text{ }\text{ }\text{ }\text{ }\Z\text{ }\text{ }\text{ }\text{ } \ar@/^/[uuul]|{\hole \text{\large{3}} \hole} \ar@/^/[uu]|{\hole \text{\large{1}} \hole} \ar@/^/[uuur] \ar@/^/[d]|{\hole \text{\large{2}} \hole} & \\
  & \text{ }\text{ }\text{ }\text{ }\Z\text{ }\text{ }\text{ }\text{ } \ar@/^/[uul]|{\hole \text{\large{3}} \hole} \ar@/^/[u]|{\hole \text{\large{1}} \hole} \ar@/^/[uur]|{\hole \text{\large{1}} \hole} & \\
  }$$

\subsection{Another example}

Sometimes we may not be able to obtain all transfer maps from the restriction maps and Mackey functor structure, in which case we need to apply \autoref{transfer}. Here is an example on its application.
\smallskip

We will still focus on $G=A_5$. We choose another specific virtual representation $V=3-V_{A_5,3}-V_{A_5,4}$ and compute $\underline{\pi}_V^G(H\underline{\Z})$.
\medskip

We have $V|_{C_3}=-2V_{C_3,2}$, $V|_{K_4}=2-2V_{K_4,3}$, $V|_{C_5}=2-3V_{C_5,2}$. The restrictions of $V$ as Mackey functors over $C_3,K_4,C_5$ are given below:
$$\xymatrix{
\F_3 \ar@/^/[d]^{0} & \\
0 \ar@/^/[u]^{0} &
}\text{ }
\xymatrix{
 & \F_2^3 \ar@/^/[d] & \\
 & \F_2 \ar@/^/[u]^0 \ar@/^/[d]^{0} & \\
 & 0 \ar@/^/[u]^{0} & 
}\text{ }
\xymatrix{
 & \F_5 \ar@/^/[d]^{0} \\
 & 0 \ar@/^/[u]^{0}
}$$
Here the Weyl group action is trivial, except that $A_4/K_4$ acts on $\F_2^3$ by permuting the 3 copies of $\F_2$. The restriction map $res_{C_2}^{K_4}:\F_2^3\rightarrow\F_2$ is the projection onto one copy of $\F_2$, depending on which $C_2\subset K_4$ we choose. We will assume the restriction map to be $1\oplus 0\oplus 0$.
\medskip

Now we can still combine \autoref{a5levels} and \autoref{restrictionmap} to compute $\Pi$ without transfer maps when different prime factors are inverted, and then recover the unlocalized version of $\Pi$. The answer is given in the left graph below. Notice that the values of $\Pi(A_5/H)$ and restriction maps are more complicated this time. We use $P$ to denote the natural projection.

$$
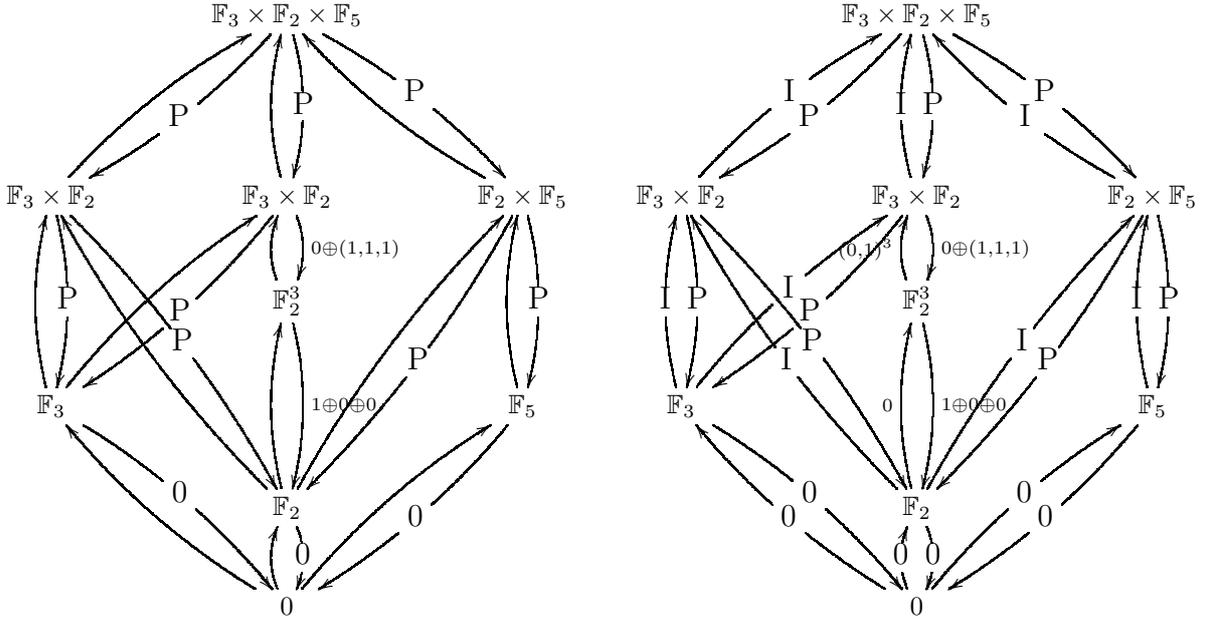
\xymatrix{
  & \text{ }\F_3\times\F_2\times\F_5\text{ } \ar@/^/[ddl]|{\hole \text{\large{P}} \hole} \ar@/^/[dd]|{\hole \text{\large{P}} \hole} \ar@/^/[ddr]|{\hole \text{\large{P}} \hole} & \\
  & & \\
  \text{ }\text{ }\text{ }\F_3\times\F_2\text{ }\text{ }\text{ } \ar@/^/[ruu] \ar@/^/[dd]|{\hole \text{\large{P}} \hole} \ar@/^/[dddr]|{\hole \text{\large{P}} \hole} &
  \text{ }\text{ }\text{ }\F_3\times\F_2\text{ }\text{ }\text{ } \ar@/^/[uu] \ar@/^/[d]^{0\oplus(1,1,1)} \ar@/^/[ddl]|{\hole \text{\large{P}} \hole} &
  \text{ }\text{ }\text{ }\F_2\times\F_5\text{ }\text{ }\text{ } \ar@/^/[uul] \ar@/^/[dd]|{\hole \text{\large{P}} \hole} \ar@/^/[dddl]|{\hole \text{\large{P}} \hole} \\
  & \text{ }\text{ }\text{ }\text{ }\F_2^3\text{ }\text{ }\text{ }\text{ } \ar@/^/[u] \ar@/^/[dd]^{1\oplus 0\oplus 0} & \\
  \text{ }\text{ }\text{ }\text{ }\F_3\text{ }\text{ }\text{ }\text{ } \ar@/^/[uu] \ar@/^/[uur] \ar@/^/[ddr]|{\hole \text{\large{0}} \hole} & &
  \text{ }\text{ }\text{ }\text{ }\F_5\text{ }\text{ }\text{ }\text{ } \ar@/^/[uu] \ar@/^/[ddl]|{\hole \text{\large{0}} \hole} \\
  & \text{ }\text{ }\text{ }\text{ }\F_2\text{ }\text{ }\text{ }\text{ } \ar@/^/[uuul] \ar@/^/[uu] \ar@/^/[uuur] \ar@/^/[d]|{\hole \text{\large{0}} \hole} & \\
  & \text{ }\text{ }\text{ }\text{ }0\text{ }\text{ }\text{ }\text{ } \ar@/^/[uul] \ar@/^/[u] \ar@/^/[uur] & \\
  }
\xymatrix{
  & \text{ }\F_3\times\F_2\times\F_5\text{ } \ar@/^/[ddl]|{\hole \text{\large{P}} \hole} \ar@/^/[dd]|{\hole \text{\large{P}} \hole} \ar@/^/[ddr]|{\hole \text{\large{P}} \hole} & \\
  & & \\
  \text{ }\text{ }\text{ }\F_3\times\F_2\text{ }\text{ }\text{ } \ar@/^/[ruu]|{\hole \text{\large{I}} \hole} \ar@/^/[dd]|{\hole \text{\large{P}} \hole} \ar@/^/[dddr]|{\hole \text{\large{P}} \hole} &
  \text{ }\text{ }\text{ }\F_3\times\F_2\text{ }\text{ }\text{ } \ar@/^/[uu]|{\hole \text{\large{I}} \hole} \ar@/^/[d]^{0\oplus(1,1,1)} \ar@/^/[ddl]|{\hole \text{\large{P}} \hole} &
  \text{ }\text{ }\text{ }\F_2\times\F_5\text{ }\text{ }\text{ } \ar@/^/[uul]|{\hole \text{\large{I}} \hole} \ar@/^/[dd]|{\hole \text{\large{P}} \hole} \ar@/^/[dddl]|{\hole \text{\large{P}} \hole} \\
  & \text{ }\text{ }\text{ }\text{ }\F_2^3\text{ }\text{ }\text{ }\text{ } \ar@/^/[u]^{(0,1)^3} \ar@/^/[dd]^{1\oplus 0\oplus 0} & \\
  \text{ }\text{ }\text{ }\text{ }\F_3\text{ }\text{ }\text{ }\text{ } \ar@/^/[uu]|{\hole \text{\large{I}} \hole} \ar@/^/[uur]|{\hole \text{\large{I}} \hole} \ar@/^/[ddr]|{\hole \text{\large{0}} \hole} & &
  \text{ }\text{ }\text{ }\text{ }\F_5\text{ }\text{ }\text{ }\text{ } \ar@/^/[uu]|{\hole \text{\large{I}} \hole} \ar@/^/[ddl]|{\hole \text{\large{0}} \hole} \\
  & \text{ }\text{ }\text{ }\text{ }\F_2\text{ }\text{ }\text{ }\text{ } \ar@/^/[uuul]|{\hole \text{\large{I}} \hole} \ar@/^/[uu]^0 \ar@/^/[uuur]|{\hole \text{\large{I}} \hole} \ar@/^/[d]|{\hole \text{\large{0}} \hole} & \\
  & \text{ }\text{ }\text{ }\text{ }0\text{ }\text{ }\text{ }\text{ } \ar@/^/[uul]|{\hole \text{\large{0}} \hole} \ar@/^/[u]|{\hole \text{\large{0}} \hole} \ar@/^/[uur]|{\hole \text{\large{0}} \hole} & \\
  }$$

The complete expression of $\Pi$ is given in the right graph above, where $I$ denotes the natural inclusion. Most of the transfer maps can be obtained from the fact that $\Pi$ is cohomological, except $tr_{K_4}^{A_4}$, for which we have to use \autoref{transfer}:
\smallskip

Assume that $3,5$ are inverted. We have $P_1=P_2=H_1=K_4$, $H_2=A_4$. The underlying $K_4$-space of $A_4/K_4$ is three copies of $K_4/K_4$ corresponding to the three $K_4$-cosets in $A_4/K_4$. Write $A_4/K_4=\{g_1K_4,g_2K_4,g_3K_4\}$. Then we have $Q_1=Q_2=Q_3=K_4$ corresponding to $g_1,g_2,g_3$.

Applying \autoref{transfer}, we have
$$tr_{K_4}^{A_4}=\bigoplus_{i=1}^3 tr_{A_4}^{A_4}\circ c_{g_i}\circ res_{K_4}^{K_4}=c_{g_1}+c_{g_2}+c_{g_3}.$$
Thus the map
$$tr_{K_4}^{A_4}: \F_2^3=\Pi(A_5/K_4)\rightarrow\Pi(A_5/A_4)\cong\Pi(A_5/K_4)^{A_4/K_4}=\F_2$$
sends each copy of $\F_2$ identically into $\F_2$. Adding the $3$-torsion into consideration, we get
$$tr_{K_4}^{A_4}:\F_2^3\xrightarrow{(0,1)^3}\F_3\times\F_2.$$

\newpage
\begin{appendices}

\section{Homology of a join}

In section 2.4, we construct the universal space $E\mathscr{F}$ as a join. We will provide more details and properties of a join in this appendix.

\begin{definition} \label{joindef}
The \textbf{join} $X\ast Y$ of two $G$-spaces $X,Y$ is defined as
$$X\ast Y:= (X\times Y\times[0,1]\sqcup X\sqcup Y)/\sim,$$
with equivalence relation given by projections
$$X\times Y\times\{0\}\rightarrow X, \text{ }X\times Y\times\{1\}\rightarrow Y.$$
The $G$-action on $X\ast Y$ is induced by the $G$-actions on $X,Y,$ and $X\times Y$.
\end{definition}

It's not hard to check:

\begin{lemma} \label{join}
For $G$-CW complexes $X,Y$, $X\ast Y$ has a natural $G$-CW structure. For any subgroup $H$, $(X\ast Y)^H$ is

\textbf{(a)} contractible if and only if at least one of $X^H, Y^H$ is contractible.

\textbf{(b)} empty if and only if both $X^H, Y^H$ are empty.
\end{lemma}
\smallskip

We can control the torsion types of $\underline{H}^G_*(X\ast Y; A_G)$ as follows:

\begin{lemma} \label{joinhom}
Assume that the Mackey functor valued homology of $X,Y, X\times Y$ with coefficients in $A_G$ only contains torsion in positive degrees. Then so does the homology of $X\ast Y$. Moreover, for any prime $p$, if $\underline{H}^G_*(X\ast Y; A_G)$ contains $p$-torsion, then at least one of $X,Y,X\times Y$ has $p$-torsion in its homology with coefficients in $A_G$.
\end{lemma}

\paragraph{Proof:} The homology of the join can be computed by the Mayer-Vietoris sequence:
$$X\ast Y=(X\times Y\times[0,1]\sqcup X\sqcup Y)/\sim=(X\times Y\times\left[0,\frac{1}{2}\right]\sqcup X/\sim)\cup(X\times Y\times\left[\frac{1}{2},1\right]\sqcup Y/\sim).$$

We write $X\ast Y$ as the union of two mapping cylinders, which are homotopic to $X,Y$ respectively. The intersection of these two cylinders is $X\times Y\times\{1/2\}\simeq X\times Y$.

Thus we have a short exact sequence
$$0\rightarrow C_*(X\times Y)\rightarrow C_*(X)\oplus C_*(Y)\rightarrow C_*(X\ast Y)\rightarrow 0$$
which passes to fixed point subspaces. So we get a long exact sequence on equivariant homology with coefficients in $A_G$:
$$...\rightarrow \underline{H}_{n+1}^G(X\ast Y;A_G)\rightarrow \underline{H}_n^G(X\times Y;A_G)\rightarrow \underline{H}_n^G(X;A_G)\oplus \underline{H}_n^G(Y;A_G)\rightarrow \underline{H}_n^G(X\ast Y;A_G)\rightarrow...$$

On degree $0$, the map
$$\underline{H}^G_0(X\times Y;A_G)\rightarrow \underline{H}^G_0(X;A_G)\oplus \underline{H}^G_0(Y;A_G)$$
is always an inclusion since the 0th homology is determined by the number of connected components for each fixed point subspace. Thus the long exact sequence above implies that the homology of $X\ast Y$ in positive degrees only contains torsion in the homology of $X,Y,X\times Y$. $\Box$

\newpage
\section{Computations on $K_4$-homology}

The computation of $\pi_\bigstar^{A_5}(H\underline{\Z})$ in section 7 is not complete since $\pi_{*+*V_{K_4,3}}^{K_4}(H\underline{\Z})$ and the $A_4/K_4$-action on it are required in the case when $3,5$ are inverted. We will compute this last missing part in this appendix.
\bigskip

For any $m,n\in\Z$, notice that $\pi^{K_4}_{m+nV_{K_4,3}}(H\underline{\Z})$ can be expressed as the equivariant homology or cohomology of $S^{|n|V_{K_4,3}}$ with coefficients in $\underline{\Z}$, which can be computed by considering an explicit CW structure on $S^{|n|V_{K_4,3}}$. In order to simplify our computation, especially for the multiplicative structure, we compare the homotopy of $H\underline{\Z}$ with the homotopy of $H\underline{\F_2}$, which is already computed in Ellis-Bloor's thesis \cite{Ell}.
\bigskip
\bigskip

Let $G=K_4$ throughout the section. 
\smallskip

We use $V_1,V_2,V_3,V$ to denote the representations $V_{K_4,1,1},V_{K_4,1,2},V_{K_4,1,3},V_{K_4,3}$. Let $H_1,H_2,H_3$ be the three proper subgroups of $K_4$ such that $V_i$ is the sign representation of $G/H_i$, $i=1,2,3$.
\bigskip

Partial computation of $\pi_\bigstar^{K_4}(H\underline{\F_2})$ is given below.

\begin{definition} \label{defplusminuscone}
The \textbf{positive cone} $\bigstar +$ of $RO(K_4)$ consists of all grades with the form $a+bV_2+cV_2+dV_3$ such that $b,c,d\leq 0$.

The \textbf{negative cone} $\bigstar -$ of $RO(K_4)$ consists of all grades with the form $a+bV_2+cV_2+dV_3$ such that $b,c,d>0$.
\end{definition}

\begin{remark} \label{plusminuscone}
Any element in $\bigstar+$ can be written as $a-W$ for some $a\in\Z$ and actual representation $W$. The homology of a point in degree $a-W$ is exactly the homology of the representation sphere $S^W$ in degree $a$. This is why we call $\bigstar+$ the positive cone. 

Similarly, any element in $\bigstar -$ can be written as $a+W$ for some actual representation $W$. The homology of a point in degree $a+W$ can be expressed as the cohomology of $S^W$ in degree $-a$. Thus we call $\bigstar-$ the negative cone.
\end{remark}

\begin{theorem} \label{positivecone}
\cite[Theorem 4.14]{Ell} Over the positive cone, we have
$$\pi_{\bigstar +}^{K_4}(H\underline{\F_2})=\frac{\F_2[x_1,y_1,x_2,y_2,x_3,y_3]}{(x_1y_2y_3+y_2x_2y_3+y_1y_2x_3)}.$$
Here $|x_i|=-V_i$, $|y_i|=1-V_i$, $i=1,2,3$.
\end{theorem}

\begin{theorem} \label{negativecone}
\cite[Proposition 4.27 and Theorem 4.30]{Ell} For the negative cone, consider the $\F_2$-linear span
$$\left\langle\frac{1}{x_1^{i_1}y_1^{j_1}x_2^{i_2}y_2^{j_2}x_3^{i_3}y_3^{j_3}}\text{ : }i_1,j_1,i_2,j_2,i_3,j_3\geq 0\right\rangle.$$
Define a self-map $f$ on the graded $\F_2$-module above by multiplication with $x_1y_2y_3+y_1x_2y_3+y_1y_2x_3$. Then
$$\pi_{\bigstar -}^{K_4}(H\underline{\F_2})=\Theta\cdot\ker(f)$$
with $|\Theta|=V-3$.
\end{theorem}

The ring structure is implied by the generators and $\Theta^2=0$. The action of $A_4/K_4\cong C_3$ is the cyclic permutation on $x_1,x_2,x_3$ and $y_1,y_2,y_3$.

\medskip

The relation between $H\underline{\Z}$ and $H\underline{\F_2}$ is given by

\begin{theorem} \label{bockstein}
\cite[Theorem 4.40]{Ell} The Bockstein spectral sequence computing the $RO(K_4)$-graded homology of a point with constant coefficients corresponding to
$$\Z\xrightarrow{2}\Z\rightarrow\F_2$$
collapses to the $E^2$-page.
\end{theorem}

In other words, $\pi_{a+bV}^{K_4}(H\underline{\Z})$ only contains $\F_2$ components unless $a+3b=0$. When $a+3b=0$, our computations will show that $\pi_{-3b+bV}^{K_4}(H\underline{\Z})=\Z$.

Applying the Bockstein long exact sequence, we know that the map
$$\pi_{a+bV}^{K_4}(H\underline{\Z})\rightarrow\pi_{a+bV}^{K_4}(H\underline{\F_2})$$
is an inclusion when $a+3b\neq 0$, and becomes the projection $\Z\rightarrow\F_2$ when $a+3b=0$. Therefore, we have

\begin{proposition} \label{keycomute}
Both the ring structure and the action of $A_4/K_4$ on $\pi_{*+*V}^{K_4}(H\underline{\Z})$ are determined by their images in $\pi_{*+*V}^{K_4}(H\underline{\F_2})$.
\end{proposition}
\medskip

Notice that $\pi_{*+*V}^{K_4}(H\underline{\Z})$ is the equivariant homology or cohomology of $S^{*V}$ with coefficients in $\underline{\Z}$. We can compute it by explicitly assigning a $G$-CW structure on $S^{*V}$:
\bigskip

For any $n>0$, $S^{nV}=S^{nV_1}\wedge S^{nV_2}\wedge S^{nV_3}$. Each $S^{nV_i}$ has the $G$-CW structure with one fixed $0$-cell $e_{i,0}$, and one cell $(G/H_i)_+\wedge e_{i,j}$ for each positive degree $0<j\leq n$.

Write $(G/H_i)_+\wedge e_{i,j}=e_{i,j}\vee e_{i,j}^\prime$. Then the boundary map is expressed as
$$\partial e_{i,j}=e_{i,j-1}+(-1)^{j-1}e_{i,j-1}^\prime,\text{ }\partial e_{i,j}^\prime=e^\prime_{i,j-1}+(-1)^{j-1}e_{i,j-1},\text{ if }j>1,$$
$$\partial e_{i,1}=\partial e_{i,1}^\prime=e_0.$$
We also formally define $e_0^\prime=e_0$.
\medskip

The $G$-CW structure on $S^{nV}$ can be obtained by smashing the structures on $S^{nV_i}$, $i=1,2,3$. For any $k,l,m\geq 0$, the elements in
$$\{e_{1,k},e_{1,k}^\prime\}\wedge\{e_{2,l},e_{2,l}^\prime\}\wedge\{e_{3,m},e_{3,m}^\prime\}$$
consist of two $G$-cells if $k,l,m>0$, or one $G$-cell otherwise. 

In particular, if exactly one of $k,l,m$ is zero, we obtain a $G$-free cell. If two of $k,l,m$ are zero, the $G$-cell has isotropy group as one of $H_1,H_2,H_3$. When $k=l=m=0$, we get the fixed cell $e_{1,0}\wedge e_{2,0}\wedge e_{3,0}$.

\subsection{The positive cone}

We give the expression of the positive cone of $\pi_{*+*V}^{K_4}(H\underline{\Z})$ first:

\begin{theorem} \label{positivez}
The homotopy $\pi_{a+bV}^{K_4}(H\underline{\Z})$, with $b\leq 0$, is the subring of
$$\frac{\Z[x_1,y_1,x_2,y_2,x_3,y_3]}{(2x_1,2x_2,2x_3,x_1y_2y_3+y_1x_2y_3+y_1y_2x_3)}$$
generated by 

\textbf{(1)} $x_1^{i_1}y_1^{j_1}x_2^{i_2}y_2^{j_2}x_3^{i_3}y_3^{j_3}$ such that
$$i_1+j_1=i_2+j_2=i_3+j_3,$$
$$j_1\equiv j_2\equiv j_3\text{ (mod 2)},$$
$$j_1j_2i_3=0.$$

\textbf{(2)} $x_1^{i_1+1}y_1^{j_1}x_2^{i_2}y_2^{j_2+1}y_3^{j_3}+x_1^{i_1}y_1^{j_1+1}x_2^{i_2+1}y_2^{j_2}y_3^{j_3}$ such that
$$i_1+j_1=i_2+j_2=j_3-1,$$
$$j_1\equiv j_2\equiv j_3\text{ (mod 2)}.$$

\textbf{(3)} $x_1^{i_1+1}y_1^{j_1}x_2^{i_2}x_3^{i_3}y_3^{j_3+1}+x_1^{i_1}y_1^{j_1+1}x_2^{i_2}x_3^{i_3+1}y_3^{j_3}$ such that
$$i_1+j_1=i_3+j_3=i_2-1,$$
$$j_1\equiv j_3\equiv 0\text{ (mod 2)}.$$

\textbf{(4)} $x_1^{i_1}x_2^{i_2+1}y_2^{j_2}x_3^{i_3}y_3^{j_3+1}+x_1^{i_1}x_2^{i_2}y_2^{j_2+1}x_3^{i_3+1}y_3^{j_3}$ such that
$$i_2+j_2=i_3+j_3=i_1-1,$$
$$j_2\equiv j_3\equiv 0\text{ (mod 2)}.$$

The action of $A_4/K_4=C_3$ is the cyclic permutation on $x_1,x_2,x_3$ and $y_1,y_2,y_3$. 
\end{theorem}

\begin{remark}
The generators we give in the theorem above are not symmetric on $x_1,x_2,x_3$ and $y_1,y_2,y_3$. However, the subring itself will be symmetric after $x_1y_2y_3+y_2x_2y_3+y_1y_2x_3$ is quotiented out.
\end{remark}

We will prove \autoref{positivez} in the rest of this section.
\bigskip

For any $n\geq 0$, $\pi_{*-nV}^{K_4}(H\underline{\Z})$ is the equivariant homology of $S^{nV}$ with coefficients in $\underline{\Z}$, which is computed as the homology of the chain complex:
$$C_*^{K_4}(S^{nV};\underline{\Z}):=\underline{C}_*(S^{nV})\otimes_{\mathcal{O}_G}\underline{\Z}.$$

We can express the generators and the boundary map with the following notations:

Use $1_H$ to denote the unit element in $\underline{\Z}(G/H)=\Z$. For any $k,l,m>0$, let $(k,l,m)$ and $(k,l,m)^\prime$ be the equivalence classes of 
$$(e_{1,k}\wedge e_{2,l}\wedge e_{3,m})\otimes 1_{\{e\}}\text{ and }(e^\prime_{1,k}\wedge e_{2,l}\wedge e_{3,m})\otimes 1_{\{e\}}.$$
Recall that $e_{1,k}^\prime$ is defined by $(G/H_1)_+\wedge e_{1,k}=e_{1,k}\vee e_{1,k}^\prime$.

Let $(k,l,0)$ be the equivalence class of
$$(e_{1,k}\wedge e_{2.l}\wedge e_{3,0})\otimes 1_{\{e\}}.$$
Define $(k,0,m)$ and $(0,l,m)$ in the same way. Let $(k,0,0)$ be the equivalence class of
$$(e_{1,k}\wedge e_{2.0}\wedge e_{3,0})\otimes 1_{H_1}.$$
Define $(0,l,0)$ and $(0,0,m)$ in the same way. Finally, $(0,0,0)$ is the equivalence class of
$$(e_{1,0}\wedge e_{2.l}\wedge e_{3,0})\otimes 1_G.$$

\begin{definition} \label{lambda1}
Let $\lambda$ be the linear self-map of $C_*^{K_4}(S^{nV};\underline{\Z})$ which exchanges $(k,l,m)$ and $(k,l,m)^\prime$ if $klm\neq 0$ and fixes $(k,l,m)$ if $klm=0$.
\end{definition}

\begin{remark} \label{lambdacoe}
For further convenience, we will use $\lambda(k,l,m)$ instead of $(k,l,m)^\prime$. Moreover, we treat $\lambda$ as part of the coefficient of $(k,l,m)$. To be precise, for $a,b\in\Z$, $(a+b\lambda)(k,l,m)$ is the image of $(k,l,m)$ under the map $a\cdot id+b\cdot\lambda$, which contains $a$ copies of cell $(k,l,m)$ and $b$ copies of cell $(k,l,m)^\prime$. We will call $a+b\lambda$ the coefficient of $(k,l,m)$ in this element.
\end{remark}
\medskip

For any degree $t$, $C_t^{K_4}(S^{nV};\underline{\Z})$ is generated by all $(k,l,m)$ and $\lambda(k,l,m)$ with $0\leq k,l,m\leq n$, $k+l+m=t$. The boundary map is given below:
$$\partial(k,l,m)=(1+(-1)^{k+1}\lambda)(k-1,l,m)$$
$$+(-1)^k(1+(-1)^{l+1}\lambda)(k,l-1,m)+(-1)^{k+l}(1+(-1)^{m+1}\lambda)(k,l,m-1)$$
if $k,l,m>1$. When $k=1$, we replace $1+(-1)^{k+1}\lambda$ by $1$. The cases $l=1$ and $m=1$ are similar.
$$\partial(k,l,0)=(1+(-1)^{k+1})(k-1,l,0)+(-1)^k(1+(-1)^{l+1})(k,l-1,0)$$
if $k,l>0$. The boundaries of $(k,0,m)$ and $(0,l,m)$ are defined similarly.
$$\partial(k,0,0)=(1+(-1)^{k+1})(k-1,0,0)$$
if $k>0$. The boundaries of $(0,l,0)$ and $(0,0,m)$ are defined similarly.

Moreover, $\partial$ commutes with $\lambda$.
\medskip

The top degree can be computed directly:

\begin{lemma} \label{tophom}
For degree $3n$, $\ker\partial$ is generated by $(1+(-1)^n\lambda)(n,n,n)$. Thus
$$H_{3n}^{K_4}(S^{nV};\underline{\Z})=\Z,\text{ }H_{3n}^{K_4}(S^{nV};\underline{\F_2})=\F_2.$$
\end{lemma}
\bigskip

It is quite hard to compute homology directly by $\ker\partial/Im\,\partial$. However, we can use the Bockstein long exact sequence
$$H_{3n}^{K_4}(S^{nV};\underline{\Z})\xrightarrow{2}H_{3n}^{K_4}(S^{nV};\underline{\Z})\rightarrow H_{3n}^{K_4}(S^{nV};\underline{\F_2})\rightarrow H_{3n-1}^{K_4}(S^{nV};\underline{\Z})\xrightarrow{2}...$$
$$...\rightarrow H_1^{K_4}(S^{nV};\underline{\F_2})\rightarrow \widetilde{H}_0^{K_4}(S^{nV};\underline{\Z})\xrightarrow{2}\widetilde{H}_0^{K_4}(S^{nV};\underline{\Z})\rightarrow \widetilde{H}_0^{K_4}(S^{nV};\underline{\F_2}).$$

Since $H_*^{K_4}(S_{nV};\underline{\Z})$ only contains $\F_2$-components except in the top degree, the long exact sequence is broken into several pieces:
$$0\rightarrow H_{3n}^{K_4}(S^{nV};\underline{\Z})\xrightarrow{2}H_{3n}^{K_4}(S^{nV};\underline{\Z})\rightarrow H_{3n}^{K_4}(S^{nV};\underline{\F_2})\rightarrow 0,$$
$$H_{3n-1}^{K_4}(S^{nV};\underline{\Z})=0,$$
$$H_i^{K_4}(S^{nV};\underline{\F_2})=H_i^{K_4}(S^{nV};\underline{\Z})\oplus \widetilde{H}_{i-1}^{K_4}(S^{nV};\underline{\Z}),\text{ }i=1,2,...,3n-1,$$
$$\widetilde{H}_0^{K_4}(S^{nV};\underline{\Z})=\widetilde{H}_0^{K_4}(S^{nV};\underline{\F_2}).$$
The first short exact sequence is $0\rightarrow\Z\rightarrow\Z\rightarrow\F_2\rightarrow 0$.
\bigskip

We can compute the homology with coefficients in $\underline{\Z}$ by the following strategy:

\textbf{(1)} Compute the dimension of $H_*^{K_4}(S^{nV};\underline{\F_2})$ as an $\F_2$-vector space. This is not hard since the dimensions of $\ker\partial$ and $Im\,\partial$ with $\F_2$-coefficients can be computed explicitly.

\textbf{(2)} Compute the dimension of $H_*^{K_4}(S^{nV};\underline{\Z})$ by the pieces of the Bockstein long exact sequence above.

\textbf{(3)} Guess elements in $\ker\partial$ (with $\Z$-coefficients) whose images are $\F_2$-independent modulo $Im\,\partial$ (with $\F_2$-coefficients). If we can find enough such elements to match the dimension of the homology, they must generate the whole homology group.
\bigskip

The boundary map with $\F_2$-coefficients is simple:
$$\partial(k,l,m)=(1+\lambda)((k-1,l,m)+(k,l-1,m)+(k,l,m-1))$$
if $k,l,m>0$. If $klm=0$, $\partial(k,l,m)=0$.

\begin{lemma} \label{f2cycle}
With $\F_2$-coefficients, $\ker\partial$ is generated by:
$$(k,l,m),\text{ }klm=0,$$
$$(1+\lambda)(k,l,m),\text{ }k,l,m>0.$$
\end{lemma}

\paragraph{Proof:} All elements listed above are inside $\ker\partial$. It suffices to show that those elements generate the whole kernel.

Assume that the coefficient of some $(k,l,m)$ ($k,l,m>0$) is $1$ or $\lambda$ in some element $Z\in\ker\partial$. Further assume that $k$ is the smallest number among all such cells. Consider the coefficient of $(k-1,l,m)$ in $\partial Z$, which should be $0$. 

If $k=1$, the coefficient of $(0,l,m)$ in $\partial(1,l,m)$ is $1$. Since $(0,l,m)$ does not appear in the boundary of any other cells, its coefficient in $\partial Z$ is $1\neq 0$, which is a contradiction.

If $k>1$, the coefficient of $(k-1,l,m)$ in $\partial(k,l,m)$ is $1+\lambda$. The cell $(k-1,l,m)$ appears in the boundaries of $(k,l,m)$, $(k-1,l+1,m)$, $(k-1,l,m+1)$. Since we already assume $k$ to be the minimal such number, the other two cells cannot have $1$ or $\lambda$ as coefficients in $Z$. Thus the coefficient of $(k-1,l,m)$ in $\partial Z$ is $1+\lambda\neq 0$, which is a contradiction.
\smallskip

Therefore, the coefficient of any $(k,l,m)$ ($k,l,m>0$) is either $1+\lambda$ or $0$ in any element in $\ker\partial$. $\Box$
\bigskip

For each dimension $t$, with $\F_2$-coefficients, the dimension of $\ker\partial$ agrees with the number of $(k,l,m)$ such that $0\leq k,l,m\leq n$ and $k+l+m=t$. The dimension of $Im\,\partial$ agrees with the number of $(k,l,m)$ such that $0<k,l,m\leq n$ and $k+l+m=t+1$. After some combinatorial arguments, we get:

\begin{proposition} \label{homdim}
The dimensions of $H_*^{K_4}(S^{nV};\underline{\F_2})$ in degrees $0,1,...,3n$ are
$$1,3,5,...,2n-1,2n+1,2n,2n-1,2n-2,...,2,1.$$
This sequence is obtained by gluing two arithmetic sequences $1,3,5,...,2n+1$ and $2n+1,2n,...,1$.
\smallskip

Thus the dimension of $H_*^{K_4}(S^{nV};\underline{\Z})$ in degrees $0,1,...,3n-1$ are
$$1,2,3,...,n,n+1,n-1,n,n-2,n-1,n-3,...,1,2,0$$
To be precise, the sequence comes from three arithmetic sequences:

$1,2,...,n+1$ for the first $n+1$ terms;

$n+1,n,n-1,...,2$ for the $(n+1)th$, $(n+3)th$, $(n+5)th$,..., $(3n-1)th$ terms;

$n,n-1,n-2,...,1,0$ for the $(n)th$, $(n+2)th$, $(n+4)th$,..., $(3n)th$ terms.
\smallskip

In addition, $H_{3n}^{K_4}(S^{nV};\underline{\Z})=\Z$.
\end{proposition}

Now we just need to guess enough elements in $\ker\partial$ with $\Z$-coefficients:

\begin{proposition} \label{homgen}
If $n$ is odd, the generators of $H_*^{K_4}(S^{nV};\underline{\Z})$ can be represented by
$$(1-\lambda)(2i+1,2j+1,n),\text{ }0\leq i,j\leq\frac{n-1}{2};$$
$$(1-\lambda)[(2i-1,2j,n)+(2i,2j-1,n)],\text{ }0<i,j\leq\frac{n-1}{2};$$
$$(2i,0,2j),(0,2i,2j),\text{ }0\leq i,j\leq\frac{n-1}{2};$$
$$(2i+1,0,2j)-(2i,0,2j+1),(0,2i+1,2j)-(0,2i,2j+1),\text{ }0\leq i,j\leq\frac{n-1}{2}.$$
If $n$ is even, the generators of $H_*^{K_4}(S^{nV};\underline{\Z})$ can be represented by
$$(1+\lambda)(2i,2j,n),\text{ }0\leq i,j\leq\frac{n}{2};$$
$$(1+\lambda)[(2i+1,2j,n)-(2i,2j+1,n)],\text{ }0\leq i,j\leq\frac{n-2}{2};$$
$$(2i,0,2j),(2j,0,2i),\text{ }0\leq i,j\leq\frac{n}{2};$$
$$(2i+1,0,2j)-(2i,0,2j+1),(0,2i+1,2j)-(0,2i,2j+1),\text{ }0\leq i,j\leq\frac{n-2}{2}.$$
Some cells should appear multiple times in different classes, in which case we will only count once.
\end{proposition}

\paragraph{Proof:} By combinatorial arguments, we can show that the number of such elements matches the dimension of the homology (although not simple). It suffices to show that, any non-zero $\F_2$-linear combination of the elements above is not in $Im\,\partial$ with $\F_2$-coefficients.

For any $Z=\sum_{\alpha\in J} (k_\alpha,l_\alpha,m_\alpha)$, assume that $\partial Z$ is a linear combination of the elements in the list above with $\F_2$-coefficients. Further assume that $|J|$ reaches its minimum.

For any $\alpha\in J$, we must have $k_\alpha,l_\alpha,m_\alpha\neq 0$. Otherwise $\partial(k_\alpha,l_\alpha,m_\alpha)=0$ and we can remove this cell to make $J$ smaller.

Choose $(k,l,m)\in\{(k_\alpha,l_\alpha,m_\alpha)\}$ such that $m$ reaches its minimum. Notice that $(k,l,m-1)$ does not appear in the boundary of any other cells inside $Z$. Thus the coefficient of $(k,l,m-1)$ in $\partial Z$ is the same as the coefficient in $\partial (k,l,m)$, which is either $1$ or $1+\lambda$, hence non-zero. But $(k,l,m-1)$ does not appear in the list above, which is a contradiction. $\Box$
\bigskip

Finally, we only need to check the image of $H_*^{K_4}(S^{nV};\underline{\Z})$ in $\pi_{\bigstar +}^{K_4}(H\underline{\F_2})$ in order to obtain the multiplicative structure and the action of $A_4/K_4=C_3$. This is not hard since the generators $x_i,y_i$ in \autoref{positivecone} come from the homology of $S^{*V_i}$, which can be decomposed into a $G/H_i=C_2$ computation. To be precise, in the homology of $S^{nV}$, the cell $(k,l,m)$ or $(1\pm\lambda)(k,l,m)$, when representing an element in homology, corrresponds to $x_1^{n-k}y_1^kx_2^{n-l}y_2^lx_3^{n-m}y_3^m$.
\smallskip

By transferring the generators in \autoref{homgen} to elements in homology, we get \autoref{positivez}.

\subsection{The negative cone}

The expression of the negative cone of $\pi_{*+*V}^{K_4}(H\underline{\Z})$ is given below:

\begin{theorem} \label{negativez}
Define $U$ as the $\F_2$-span of all $x_1^{-i_1}y_1^{-j_1}x_2^{-i_2}y_2^{-j_2}x_3^{-i_3}y_3^{-j_3}$ such that
$$i_1+j_1=i_2+j_2=i_3+j_3,$$
$$j_1,j_2,j_3>0,\text{ }i_1,i_2,i_3\geq 0,$$
$$(i_1,i_2,i_3)\neq (0,0,0).$$

Let $f$ be the self-map on $U$ given by multiplying $x_1y_2y_3+y_1x_2y_3+y_1y_2x_3$.

Let $T$ be the sub-module of $U$ generated by all
$$x_1^{-i_1}y_1^{-j_1-1}x_2^{-i_2-1}y_2^{-j_2}x_3^{-i_3-1}y_3^{-j_3}+x_1^{-i_1-1}y_1^{-j_1}x_2^{-i_2}y_2^{-j_2-1}x_3^{-i_3-1}y_3^{-j_3}$$
$$+x_1^{-i_1-1}y_1^{-j_1}x_2^{-i_2-1}y_2^{-j_2}x_3^{-i_3}y_3^{-j_3-1}$$
such that $j_1\equiv j_2\equiv j_3$ (mod 2).
\medskip

The homotopy $\pi_{a+bV}^{K_4}(H\underline{\Z})$, with $b>0$, is given by
$$4\Z\langle\bigoplus_{n>0}(y_1y_2y_3)^{-n}\rangle\oplus(\ker(f)\cap T).$$
\end{theorem}

We will use the rest of this section to prove \autoref{negativez}.
\bigskip

For any $n>0$, $\pi_{*+nV}^{K_4}(H\underline{\Z})$ is the equivariant cohomology of $S^{nV}$ with coefficients in $\underline{\Z}$, which is computed as the cohomology of the cochain complex
$$C_{K_4}^*(S^{nV};\underline{\Z}):=Hom_{\mathcal{O}_G}(\underline{C}_*(S^{nV}),\underline{\Z}).$$

Since all restriction maps in $\underline{\Z}$ are the identity, the cochain complex above agrees with
$$Hom_\Z(C_*(S^{nV}/G),\Z).$$
Thus we are computing the non-equivariant cohomology of the orbit space $S^{nV}/G$.
\medskip

We can express the cochain complex and the coboundary map with the following notations:

For any $k,l,m$, let $[k,l,m]$ be the function sending the orbit of
$$e_{1,k}\wedge e_{2,l}\wedge e_{3,m}$$
to $1$ and all other cells to $0$. 

When $klm\neq 0$, let $[k,l,m]^\prime$ be the function sending the orbit of 
$$e_{1,k}^\prime\wedge e_{2,l}\wedge e_{3,m}$$
to $1$ and all other cells to $0$.

Moreover, when some of $k,l,m$ are greater than $n$ or less than $0$, we write $[k,l,m]=0$.
\medskip

\begin{definition}
Let $\lambda$ be the linear self-map of $C^*_{K_4}(S^{nV};\underline{\Z})$ which exchanges $[k,l,m]$ and $[k,l,m]^\prime$ if $klm\neq 0$ and fixes $[k,l,m]$ if $klm=0$. We treat $\lambda$ as part of the coefficient of $[k,l,m]$ in the same way as in \autoref{lambdacoe}.
\end{definition}
\medskip

For any degree $t$, $C^t(S^{nV}/G;\Z)$ is generated by all $[k,l,m]$ and $\lambda[k,l,m]$ with $0\leq k,l,m\leq n$, $k+l+m=t$. The coboundary map is given below:
$$\delta[k,l,m]=(1+(-1)^k\lambda)[k+1,l,m]$$
$$+(-1)^k(1+(-1)^l\lambda)[k,l+1,m]+(-1)^{k+l}(1+(-1)^m\lambda)[k,l,m+1]$$
if there is at most one $0$ in $\{k,l,m\}$. In addition, we have
$$\delta[k,0,0]=(1+(-1)^k)[k+1,0,0]+(-1)^k[k,1,0]+(-1)^k[k,0,1].$$
The coboundaries of $[0,l,0]$ and $[0,0,m]$ are defined similarly.
\medskip

Again, we can compute the top degree directly:

\begin{lemma} \label{topcohom}
For degree $3n$, $Im\,\delta$ is generated by $(1+(-1)^{n-1}\lambda)[n,n,n]$. Thus $H^{3n}_{K_4}(S^{nV};\underline{\Z})=\Z$, $H^{3n}_{K_4}(S^{nV};\underline{\F_2})=\F_2$.
\end{lemma}
\medskip

For the remaining degrees, first we describe $Im\,\delta$:

\begin{lemma} \label{f2coboundary}
With $\F_2$-coefficients, $Im\,\delta$ is generated by:
$$(1+\lambda)[k,l,m],\text{ }k,l,m>0,$$
$$\delta[0,0,m],\delta[0,l,0],\delta[k,0,0].$$
\end{lemma}

\paragraph{Proof:} It's not hard to check that, except $\delta[0,0,m],\delta[0,l,0],\delta[k,0,0]$, all other coboundaries are sums of $(1+\lambda)[k,l,m]$ with $k,l,m>0$. It suffices to check that any $(1+\lambda)[k,l,m]$ can be expressed as a coboundary.

Consider an induction on $\min\{k,l,m\}$. Without loss of generality, assume that $k=\min\{k,l,m\}$. The base case $k=1$ is given by $\delta[0,l,m]=(1+\lambda)[1,l,m]$.

From case $i$ to $i+1$: When $l,m\geq k=i+1$, we have
$$\delta[i,l,m]=(1+\lambda)([i+1,l,m]+[i,l+1,m]+[i,l,m+1]).$$
By induction, $(1+\lambda)[i,l+1,m]$ and $(1+\lambda)[i,l,m+1]$ are in $Im\,\delta$. Thus $(1+\lambda)[i+1,l,m]$ is also in $Im\,\delta$. $\Box$
\bigskip

Unlike the computation of the positive cone in section 7.1, we cannot explicitly guess the generators, since they are quite complicated. Instead, we will point out the types of cocyles with $\F_2$-coefficients which can be lifted to cocycles with $\Z$-coefficients.

\begin{definition} \label{triangle}
For $k,l,m\geq 0$ define
$$\langle k,l,m\rangle^+:=[2k+1,2l,2m]+[2k,2l+1,2m]+[2k,2l,2m+1],$$
$$\langle k,l,m\rangle^-:=[2k,2l-1,2m-1]-[2k-1,2l,2m-1]+[2k-1,2l-1,2m].$$
\end{definition}
\medskip

Notice that when $k=0$, $\langle k,l,m\rangle^-$ contains one single term $[0,2l-1,2m-1]$. The cases when $l=0$ or $m=0$ are similar.

\begin{proposition} \label{cohomkey}
\textbf{(a)} Except in the top degree, any cocycle $Z$ with $\Z$-coefficients is a linear combination of $\langle k,l,m\rangle^+$ and $\langle k,l,m\rangle^-$ modulo $(2,1+\lambda)$.

\textbf{(b)} On the other hand, if the sum of some $\langle k,l,m\rangle^+$ and $\langle k,l,m\rangle^-$ is a cocycle with $\F_2$-coefficients, it can be lifted to a cocycle with $\Z$-coefficients.
\end{proposition}

Since $(1+\lambda)C^*(S^{nV},\F_2)$ is contained in $Im\,\delta$ according to \autoref{f2coboundary}, this proposition gives us a complete description of the image of $H^*(S^{nV};\Z)$ inside $H^*(S^{nV};\F_2)$.
\bigskip

\paragraph{Proof of \autoref{cohomkey} part (a):} First, we show that any generators of the cochain complex with forms $[2k,2l,2m]$ and $[2k+1,2l+1,2m+1]$ do not appear in a cocycle modulo $(2,1+\lambda)$.

Without loss of generality, assume that $k<n$. Consider the element $[2k+1,2l,2m]$, which appears in the coboundary of $[2k,2l,2m]$, $[2k+1,2l-1,2m]$, $[2k+1,2l,2m-1]$ (or multiplied by $\lambda$), with coefficients $1+\lambda$, $\pm(1-\lambda)$, $\pm(1-\lambda)$. The coefficients of $[2k+1,2l,2m]$ in $\delta Z$ can never be zero unless the coefficient of $[2k,2l,2m]$ in $Z$ is a multiple of $1-\lambda\in(2,1+\lambda)$.

The case of $[2k+1,2l+1,2m+1]$ can be proved in a similar way.
\medskip

Next, we show that if $[2k+1,2l,2m]$ appears in some cocycle $Z$ modulo $(2,1+\lambda)$, then $[2k,2l+1,2m]$ must also appear.

Consider the element $[2k+1,2l+1,2m]$, which appears in the coboundaries of $[2k+1,2l,2m]$, $[2k,2l+1,2m]$, $[2k+1,2l+1,2m-1]$ (or multiplied by $\lambda$), with coefficients $\pm(1-\lambda),\pm(1-\lambda),1+\lambda$. Since the coefficient of $[2k+1,2l+1,2m]$ is $0$ in $\delta Z$, the only possibility is that $[2k+1,2l+1,2m-1]$ does not appear in $Z$, while $[2k+1,2l,2m]$ and $[2k,2l+1,2m]$ have the same coefficient modulo $(2,1+\lambda)$.

The same argument can be applied to any other pair of components in $\langle k,l,m\rangle^+$ or $\langle k,l,m\rangle^-$.
\smallskip

In conclusion, we proved that any cocycle must be a linear combination of $\langle k,l,m\rangle^+$ and $\langle k,l,m\rangle^-$ modulo $(2,1+\lambda)$. $\Box$
\bigskip

\paragraph{Proof of \autoref{cohomkey} part (b):} We only consider odd degrees, for which we only have $\langle k,l,m\rangle^+$. The case of even degrees can be proved in a similar way.

Notice that
$$\delta\langle k,l,m\rangle^+=(1-\lambda)([2k+2,2l,2m]+[2k,2l+2,2m]+[2k,2l,2m+2]).$$
Moreover, each $[2k,2l,2m]$ only appears in the coboundaries of $\langle k-1,l,m\rangle^+$, $\langle k,l-1,m\rangle^+$, $\langle k,l,m-1\rangle^+$, with the same coefficient $1-\lambda$.

Assume that $\sum_{\alpha\in J}\langle k_\alpha,l_\alpha,m_\alpha\rangle^+$ is a cocycle modulo 2. Then with $\Z$-coefficients, the coboundary consists of $2(1-\lambda)[2k,2l,2m]$ for some $k,l,m$.

Consider any such $[2k,2l,2m]$. Find one $\langle k_\alpha,l_\alpha,m_\alpha\rangle^+$ whose coboundary contains $(1-\lambda)[2k,2l,2m]$. Without loss of generality, assume that we have $\langle k-1,m,n\rangle^+$. We change the component of $[2k-1,2m,2n]$ to $\lambda[2k-1,2m,2n]$. The image of $\langle k-1,m,n\rangle^+$ is unchanged modulo $(2,1+\lambda)$. However, $\delta\langle k-1,m,n\rangle^+$ is changed from
$$(1-\lambda)([2k,2l,2m]+[2k-2,2l+2,2m]+[2k-2,2l,2m+2])$$
to
$$(1-\lambda)(-[2k,2l,2m]+[2k-2,2l+2,2m]+[2k-2,2l,2m+2]).$$
Thus we eliminate $2(1-\lambda)[2k,2l,2m]$ inside $\delta\sum_{\alpha\in J}\langle k_\alpha,l_\alpha,m_\alpha\rangle^+$.
\smallskip

Applying the same procedure for each $[2k,2l,2m]$. Finally we can make $\sum_{\alpha\in J}\langle k_\alpha,l_\alpha,m_\alpha\rangle^+$ into a cocycle with $\Z$-coefficients. $\Box$
\bigskip

Now we still have the last piece in order to finish the ring structure of the negative cone:

\begin{lemma} \label{theta4}
Consider $\Theta$ and $y_1y_2y_3$ as generators of $\pi_{V-3}^{K_4}(H\underline{\Z})$ and $\pi_{3-V}^{K_4}(H\underline{\Z})$. We have $\Theta y_1y_2y_3=4$.
\end{lemma}

\paragraph{Proof:} Consider the Mackey functor valued homotopy:
$$\underline{\pi}_{V-3}^{K_4}(H\underline{\Z})\text{ and }\underline{\pi}_{3-V}^{K_4}(H\underline{\Z}).$$

For any $H\subset G=K_4$, we have
$$\underline{\pi}_{3-V}^{K_4}(H\underline{\Z})(G/H)=H_3^H(S^V;\underline{\Z})=\Z,$$
$$\underline{\pi}_{V-3}^{K_4}(H\underline{\Z})(G/H)=H_H^3(S^V;\underline{\Z})=\Z.$$
The restriction and transfer maps are computed in [\text{Ang22}, Section 8]:
$$\underline{\pi}_{V-3}^{K_4}(H\underline{\Z})=\widehat{\underline{\Z}},\text{ } \underline{\pi}_{3-V}^{K_4}(H\underline{\Z})=\underline{\Z}.$$
Here $\widehat{\underline{\Z}}$ is the $K_4$-Mackey functor with $\Z$-values and identity transfer maps.

The product between $\Theta$ and $y_1y_2y_3$ can be computed from the $G/G$-value of
$$\underline{\pi}_{3-V}^{K_4}(H\underline{\Z})\text{ }\Box\text{ }\underline{\pi}_{V-3}^{K_4}(H\underline{\Z})\rightarrow \underline{\pi}_0^{K_4}(H\underline{\Z}),$$
which is in fact $\widehat{\underline{\Z}}\text{ }\Box\text{ }\underline{\Z}\rightarrow\underline{\Z}$. On the $G/e$-value, we have the common multiplication $\Z\times\Z\rightarrow\Z$. Thus on the $G/G$-value, $(1,1)\in\Z\times\Z$ is sent to $4\in\Z$. $\Box$
\bigskip

\autoref{negativez} can be proved by combining $\autoref{negativecone}$, $\autoref{cohomkey}$, and $\autoref{theta4}$.

\end{appendices}

\bibliographystyle{abbrv}
%references

\end{document}